\documentclass[hidelinks,onefignum,onetabnum]{siamart220329}
\usepackage{amssymb}
\usepackage{amsmath}
\usepackage{enumitem}
\usepackage{bm}
\usepackage{color}
\usepackage{exscale}
\usepackage{relsize}
\usepackage{graphicx}
\DeclareSymbolFont{yhlargesymbols}{OMX}{yhex}{m}{n}
\DeclareSymbolFont{ugmL}{OMX}{mdugm}{m}{n}
\DeclareMathAccent{\wideparen}{\mathord}{yhlargesymbols}{"F3}
\usepackage{wrapfig}%(P227)
\usepackage{multicol}
\usepackage{verbatim}
\theoremstyle{remark}%
\newtheorem{rmk}{Remark}
\usepackage{appendix}
\usepackage{lipsum}
\usepackage{amsfonts}
\usepackage{graphicx}
\usepackage{epstopdf}
\usepackage{algorithmic}
\ifpdf
  \DeclareGraphicsExtensions{.eps,.pdf,.png,.jpg}
\else
  \DeclareGraphicsExtensions{.eps}
\fi

\newsiamremark{remark}{Remark}
\newsiamremark{hypothesis}{Hypothesis}
\crefname{hypothesis}{Hypothesis}{Hypotheses}
\newsiamthm{claim}{Claim}

\headers{Higher-order multiscale method for thermo-electric coupling}{H. Dong, Z. Z. Yang and Y. F. Nie}

\title{Higher-order multiscale method and its convergence analysis for nonlinear thermo-electric coupling problems of composite structures\thanks{Submitted to the editors' DATE.
\funding{This work was supported by the National Natural Science Foundation of China (No.\hspace{1mm}12471387), the Young Talent Fund of Association for Science and Technology in Shaanxi, China (No.\hspace{1mm}20220506), Xidian University Specially Funded Project for Interdisciplinary Exploration (No.\hspace{1mm}TZJH2024008).
}}}

\author{Hao Dong\thanks{Corresponding author. School of Mathematics and Statistics, Xidian University, Xi'an 710071, China
  (\email{donghaoxd@xidian.edu.cn}).}
\and Zongze Yang\thanks{Department of Applied Mathematics, The Hong Kong Polytechnic University, Hong Kong 999077, China (\email{zongze.yang@polyu.edu.hk}).}
\and Yufeng Nie\thanks{School of Mathematics and Statistics, Northwestern Polytechnical University, Xi'an 710129, China (\email{yfnie@nwpu.edu.cn}).}
}

\usepackage{amsopn}

\ifpdf
\hypersetup{
  pdftitle={Higher-order multiscale method for thermo-electric coupling problems},
  pdfauthor={H. Dong, Z. Z. Yang and Y. F. Nie}
}
\fi

\begin{document}

\maketitle

% REQUIRED
\begin{abstract}
This paper proposes a higher-order multiscale computational method for nonlinear thermo-electric coupling problems of composite structures, which possess temperature-dependent material properties and nonlinear Joule heating. The innovative contributions of this work are the novel multiscale formulation with the higher-order correction terms for periodic composite structures and the global error estimation with an explicit rate for higher-order multiscale solutions. By employing the multiscale asymptotic approach and the Taylor series technique, the higher-order multiscale method is established for time-dependent nonlinear thermo-electric coupling problems, which can keep the local balance of heat flux and electric charge for high-accuracy multiscale simulation. Furthermore, an efficient numerical algorithm with off-line and on-line stages is presented in detail, and corresponding convergent analysis is also obtained. Two- and three-dimensional numerical experiments are conducted to showcase the competitive advantages of the proposed method for simulating the time-dependent nonlinear thermo-electric coupling problems in composite structures, not only exceptional numerical accuracy, but also less computational cost.
\end{abstract}

% REQUIRED

% REQUIRED

\section{Introduction}
As the integrated circuit industry enters the post-Moore's Law era, the integrated circuits have shifted from the traditional approach of continuously reducing structural sizes to enhancing system performance through advanced packaging technologies. Three-dimensional (3D) packaging technology is regarded as a crucial scheme for extending Moore's Law. Electronic packaging structures involve a variety of materials with distinct properties, while their components differ by several orders of magnitude in spatial size. Hence, electronic packaging structure is a complex three-dimensional composite structure with multiscale spatial features \cite{R1}. Electronic packaging structures server under extremely electrical-thermal coupling environments, whose localized heat flux densities can reach up to $1kW/cm^2$. These composite structures will exhibit significant nonlinear effects, including temperature-dependent material properties and Joule heating effects \cite{R2,R3}, which leads to tremendous cost by use of classical numerical methods. Therefore, efficient and high-accuracy multiscale computation is of paramount importance for performance simulation and prediction for nonlinear thermo-electric coupling behaviors of composite structures in electronic packaging industry.

On the basis of the first law of thermodynamics, the Fourier's law and the conservation law of electric charge \cite{R3,R4} in the continuum mechanics (CM) framework, the governing equations for time-dependent nonlinear thermo-electric coupling problems of composite structures with periodically microscopic configurations can be established in domain $\Omega$, where $\Omega\in\mathbb{R}^n(n=2,3)$ is a bounded convex domain with the boundary $\partial\Omega=\partial\Omega_{u}\cup\partial\Omega_{\phi}\cup\partial\Omega_{q}\cup\partial\Omega_{d}$.
\begin{equation}
\left\{
\begin{aligned}
&\rho^{\varepsilon} ({\bm{x}},{u^{\varepsilon}})c^{\varepsilon}({\bm{x}},{u^{\varepsilon}})\frac{{\partial {u^{\varepsilon}(\bm{x},t)}}}{{\partial t}}-\frac{\partial }{{\partial {x_i}}}\Big( {{k_{ij}^{\varepsilon}}({\bm{x}},{u^{\varepsilon}})\frac{{\partial {u^{\varepsilon}(\bm{x},t)}}}{{\partial {x_j}}}}\Big)\\
&\quad\quad\quad\quad\quad\quad={{\sigma_{ij}^{\varepsilon}}({\bm{x}},{u^{\varepsilon}})\frac{\partial \phi^{\varepsilon}(\bm{x},t)}{\partial {x_i}}}\frac{\partial \phi^{\varepsilon}(\bm{x},t)}{\partial {x_j}}+f_u(\bm{x},t),\;\;\text{in}\;\;\Omega\times(0,T),\\
&- \frac{\partial }{{\partial {x_i}}}\Big({{\sigma_{ij}^{\varepsilon}}({\bm{x}},{u^{\varepsilon} })\frac{{\partial \phi^{\varepsilon}(\bm{x},t)}}{{\partial {x_j}}}} \Big) = f_\phi(\bm{x},t),\;\;\text{in}\;\;\Omega\times(0,T),\\
&u^{{\varepsilon}}(\bm{x},t) = \widehat u(\bm{x},t),\;\;\text{on}\;\;\partial {\Omega_u}\times(0,T),\\
&\phi^{{\varepsilon}}(\bm{x},t)=\widehat\phi(\bm{x},t),\;\;\text{on}\;\;\partial\Omega_{\phi}\times(0,T),\\
&{k_{ij}^{{\varepsilon}}({\bm{x}},{u^{\varepsilon}})\frac{\partial u^{{\varepsilon}}(\bm{x},t)}{\partial {x_j}}}{n_i} = \bar q(\bm{x},t),\;\;\text{on}\;\;\partial {\Omega_q}\times(0,T),\\
&{\sigma_{ij}^{{\varepsilon}}({\bm{x}},{u^{\varepsilon}})\frac{\partial \phi^{{\varepsilon}}(\bm{x},t)}{\partial {x_j}}}{n_i}=\bar{d}(\bm{x},t),\;\;\text{on}\;\;\partial\Omega_{d}\times(0,T),\\
&u^{{\varepsilon}}({\bm{x}},0)=\widetilde u(\bm{x}),\;\;\text{in}\;\;\Omega.
\end{aligned} \right.
\end{equation}
In mathematical model (1.1), the temperature field $u^{\varepsilon}(\bm{x},t)$ and electric potential field $\phi^{\varepsilon}(\bm{x},t)$ are targeted and undetermined, in which $\varepsilon$ represents the characteristic length of periodic unit cell. $\widehat u(\bm{x},t)$ and $\widehat\phi(\bm{x},t)$ are the prescribed temperature and electric potential on the domain boundaries $\partial\Omega_{u}$ and $\partial\Omega_{\phi}$ respectively. $\bar q(\bm{x},t)$ and $\bar{d}(\bm{x},t)$ are the prescribed heat flux and surface electric charge on the domain boundaries $\partial\Omega_{q}$ and $\partial\Omega_{d}$ respectively. $\widetilde u(\bm{x})$ represents the initial temperature of domain $\Omega$. $f_u(\bm{x},t)$ and $f_\phi(\bm{x},t)$ denote the internal heat source and electric charge density respectively. Here $\rho^{\varepsilon}({\bm{x}},{u^{\varepsilon}})$, $c^{\varepsilon}({\bm{x}},{u^{\varepsilon}})$ $\{{k_{ij}^{\varepsilon}}({\bm{x}},{u^{\varepsilon}})\}$ and $\{{\sigma_{ij}^{\varepsilon}}({\bm{x}},{u^{\varepsilon}})\}$ are, respectively, the mass density, the specific heat, the thermal conductivity tensor and the electric conductivity tensor, which all are temperature-dependent. Mathematical model (1.1) appears in the thermal and electric coupling simulation of composite structures in the field of electronic packaging, whose multiscale nature stems from the high-frequency oscillating coefficients owing to periodic heterogeneity at microscale when $0<\varepsilon\ll1$, necessitating an exceptionally refined discretization for high-resolution simulation and making the overall overhead prohibitive, particularly for time-dependent problems. Additionally, the Laplace transform technique \cite{R5,R38,R39} for time-dependent linear system can not work for theoretically analyzing and numerically simulating the nonlinear model (1.1) due to its temperature-dependent material properties (nonlinear equation coefficients) and Joule heating effects (nonlinear product terms of gradient).

In the last several decades, theoretical analysis and numerical computation for the time-dependent thermo-electric coupling problems has been done extensively. In 1992 and 2006, Allegretto et al. established the existence theory for solutions to nonlinear thermo-electric coupling problems \cite{R6} and developed a posteriori error analysis theories for finite element method solving these nonlinear coupling problems \cite{R7}. In 2005, Akrivis and Larsson proposed a linearly implicit finite element method for solving nonlinear thermo-electric coupling problems and demonstrated optimal error estimates under the assumption of sufficiently regular solutions \cite{R8}. In 2014 and 2021, Sun, Li and Gao systematically developed Crank-Nicolson finite element method \cite{R9} and mixed finite element method \cite{R10} for solving nonlinear thermo-electric coupling problems, and obtained optimal error estimates based on space-time error splitting technique. In 2019 and 2023, Shi et al. developed nonconforming finite element method \cite{R11} and linearized Galerkin finite element method \cite{R12} for nonlinear thermo-electric coupling problems, and also achieved unconditionally superconvergent error estimates. In 2010, Fern$\rm{\acute{a}}$ndez and Kuttler established a fully discrete numerical scheme for the nonlinear electro-thermo-mechanical coupling problem based on finite element method and a forward Euler scheme, and also obtained a linear convergence estimate under appropriate regularity assumptions \cite{R13}. In 2017, M$\rm{\mathring{a}}$lqvist and Stillfjord presented a fully discrete numerical scheme for the nonlinear electro-thermo-mechanical coupling problem based on finite element method in space and a semi-implicit Euler scheme in time, while proving optimal convergence orders, i.e. second-order in space and first-order in time \cite{R14}. In 2023, Jiang et al. established a Galerkin mixed finite element method for the nonlinear electro-thermo-mechanical coupling problem, and employed space-time error splitting technique to obtain optimal error estimate \cite{R15}. However, the above-mentioned computational methods are only applicable for the computation of homogeneous materials and don't consider the impact of microscopic material heterogeneities on nonlinear multiphysics behaviors of composite structures.

To overcome the challenging issues inherent in the multiscale feature of composite structures, a diverse array of multiscale methods have been established by researchers to balance the accuracy and efficiency, such as asymptotic homogenization method (AHM) \cite{R16,R17}, multiscale finite element method (MsFEM) \cite{R18,R19}, heterogeneous multiscale method (HMM) \cite{R20}, variational multiscale method (VMS) \cite{R21}, multiscale eigenelement method (MEM) \cite{R22} and localized orthogonal decomposition method (LOD) \cite{R23}, etc. In the past thirty years, Cui and his research team systematically developed a class of higher-order multiscale approaches for precisely and efficiently simulating the thermal, mechanical and multiphysics behaviors of composite structures, as shown in references \cite{R5,R24,R25,R26,R27} for the more details. Nevertheless, it is essential to note that the aforementioned multiscale methodologies are often restricted to the computation of linear multiscale problems in the practical scenarios. With the expansion of engineering applications and the demand for high-accuracy engineering simulations of composite structures, research attention of scientists and engineers inevitably focused on multiscale modeling and computation for nonlinear problems of inhomogeneous solids. In references \cite{R28,R29,R30,R31}, researchers employed the asymptotic homogenization method to analyze and simulate static and transient nonlinear heat conduction problems of inhomogeneous solids with temperature-dependent conductivity coefficients. In \cite{R32}, Pankov systematically studied the G-convergence and homogenized theory of nonlinear partial differential operators. In \cite{R19}, Efendiev et al. established multiscale finite element methods for nonlinear elliptic equations and nonlinear parabolic equations. Furthermore, Dong et al. first developed higher-order multiscale method and corresponding convergence analysis for nonlinear thermo-mechanical coupling problems of composite structures and shells \cite{R33,R34}. To summarize, few studies have been published in which the nonlinear multiphysics problems are taken into account in the multiscale modeling and computation of composite structures. However, widespread engineering demands strongly prompt continued research about this challenging issue, especially for nonlinear thermo-electric coupling problems.

The remainder of this study is structured as below. Section 2 establishes the higher-order multiscale computational model for nonlinear thermo-electric coupling simulation of composite structures with microscopic periodic configurations by virtue of the multiscale asymptotic approach and the Taylor series technique. In Section 3, both local and global error analyses are established for the proposed multiscale solutions. Section 4 provides a two-stage numerical algorithm with off-line microscale computation, and on-line macroscale and multiscale computation to effectively simulate time-dependent nonlinear thermo-electric coupling problems of composite structures at length. The corresponding convergent analysis of the two-stage numerical algorithm is obtained in Section 5. Numerical examples are designed in Section 6 to validate the computational performance of the presented multiscale computational model and corresponding numerical algorithm. Eventually, concluding remarks and potential directions are proclaimed in Section 7. Throughout this paper, Einstein summation convention is employed to simplify repetitious indices.

\section{Higher-order multiscale computational model of nonlinear thermo-electric coupling problems}
\subsection{The statement of multiscale nonlinear coupling system}
In accordance with the framework of the asymptotic homogenization method, we adeptly denote $\bm{y}={\bm{x}}/{\varepsilon}=({x_1}/{\varepsilon},\cdots,{x_n}/{\varepsilon})=(y_1,\cdots,y_n)$ as microscopic coordinates of periodic unit cell (PUC) $\Theta=[0,1]^n$. As a result, material parameters $\rho^{\varepsilon} ({\bm{x}},{u^{\varepsilon}})$, $c^{\varepsilon} ({\bm{x}},{u^{\varepsilon}})$, ${k_{ij}^{\varepsilon}}({\bm{x}},{u^{\varepsilon}})$ and ${\sigma_{ij}^{\varepsilon}}({\bm{x}},{u^{\varepsilon}})$ in nonlinear thermo-electric system (1.1) can be expressed with new forms $\rho({\bm{y}},{u^{\varepsilon}})$, $c({\bm{y}},{u^{\varepsilon}})$, ${k_{ij}}({\bm{y}},{u^{\varepsilon}})$ and ${\sigma_{ij}}({\bm{y}},{u^{\varepsilon}})$, which imply that these material parameters are 1-periodic functions in microvariable ${\bm{y}}$.

Following the previous works \cite{R9,R24,R33,R34}, some assumptions for multiscale nonlinear equations (1.1) are presented as follows.
\begin{enumerate}
\item[(A$_1$)]
${k_{ij}^{\varepsilon}}({\bm{x}},{u^{\varepsilon}})$ and ${\sigma_{ij}^{\varepsilon}}({\bm{x}},{u^{\varepsilon}})$ are symmetric, and there exist two positive constants $\gamma_0$ and $\gamma_1$ irrespective of $\varepsilon$ for the following uniform elliptic conditions
\begin{displaymath}
\begin{aligned}
&k_{ij}^{\varepsilon}=k_{ji}^{\varepsilon},\;\gamma_0|\bm{\zeta}|^2\leq {k_{ij}^{\varepsilon}}({\bm{x}},{u^{\varepsilon}})\zeta_i\zeta_j  \le\gamma_1|\bm{\zeta}|^2,\\
&\sigma_{ij}^{\varepsilon}=\sigma_{ji}^{\varepsilon},\;\gamma_0|\bm{\zeta}|^2\leq {\sigma_{ij}^{\varepsilon}}({\bm{x}},{u^{\varepsilon}})\zeta_i\zeta_j  \le\gamma_1|\bm{\zeta}|^2,
\end{aligned}
\end{displaymath}
where $\bm{\zeta}=(\zeta_1,\cdots,\zeta_n)$ is an arbitrary vector with
real elements in $\mathbb{R}^n$, and ${\bm{x}}$ is an arbitrary point in $\Omega$.
\item[(A$_2$)]
$\rho^{\varepsilon} ({\bm{x}},{u^{\varepsilon}})$, $c^{\varepsilon} ({\bm{x}},{u^{\varepsilon}})$, ${k_{ij}^{\varepsilon}}({\bm{x}},{u^{\varepsilon}})$ and ${\sigma_{ij}^{\varepsilon}}({\bm{x}},{u^{\varepsilon}})\in L^\infty (\Omega)$; $0<\rho^{*}\leq\rho^{\varepsilon} ({\bm{x}},{u^{\varepsilon}})$, $0<c^{*}\leq c^{\varepsilon} ({\bm{x}},{u^{\varepsilon}})$, where $\rho^{*}$ and $c^{*}$ are constants irrespective of $\varepsilon$.
\item[(A$_3$)]
$f_{u}\in L^2(\Omega\times(0,T))$, $f_{\phi}\in L^2(\Omega\times(0,T))$, ${\widehat u}({\bm{x}},t)\in L^2(0,T;H^{1}(\Omega))$, $\widehat{\phi}({\bm{\alpha }},t)\in L^2(0,T;H^{1}(\Omega))$, $\bar q\in L^2(\Omega\times(0,T))$, $\bar d\in L^2(\Omega\times(0,T))$, $\widetilde{u}(\bm{x})\in {L^2 (\Omega)}$.
\end{enumerate}
\subsection{Higher-order multiscale analysis for nonlinear coupling problem}
In order to establish higher-order multiscale computational model, one firstly suppose that the basic field variables ${u^{\varepsilon} }(\bm{x},t)$ and $\phi^{\varepsilon}(\bm{x},t)$ are represented by the succeeding two-scale asymptotic expansion forms inspired by \cite{R16,R17}.
\begin{equation}
\left\{\begin{array}{l}
{u^{\varepsilon} }(\bm{x},t)={u_{0}}(\bm{x},\bm{y},t) + \varepsilon {u_{1}}(\bm{x},\bm{y},t) + {\varepsilon ^2}{u_{2}}(\bm{x},\bm{y},t) + {\rm O}({\varepsilon ^3}),\\
\phi^{\varepsilon}(\bm{x},t)=\phi_{0}(\bm{x},\bm{y},t) + \varepsilon \phi_{1}(\bm{x},\bm{y},t) + {\varepsilon ^2}\phi_{2}(\bm{x},\bm{y},t) + {\rm O}({\varepsilon ^3}).
\end{array}\right.
\end{equation}
In the preceding formula, $u_0$ and $\phi_0$ are zeroth-order expansion terms, $u_1$ and $\phi_1$ are first-order asymptotic terms (lower-order asymptotic terms), $u_2$ and $\phi_2$ are second-order asymptotic terms (higher-order asymptotic terms).

After that, the key idea to handle the temperature-dependent material parameters is introduced. Drawing support from the Taylor's formula with multi-index notation $\displaystyle f({x_0},{y_0} + \delta) = f({x_0},{y_0}) + {f_y}({x_0},{y_0})\delta + \frac{1}{2}{f_{yy}}({x_0},{y_0}){\delta^2} + O({\delta^3})= f({x_0},{y_0}) + \mathbf{D}^{(0,1)}{f}({x_0},{y_0})\delta + \frac{1}{2}\mathbf{D}^{(0,2)}{f}({x_0},{y_0}){\delta^2} + O({\delta^3})$ in \cite{R35}, material parameter $k_{ij}^\varepsilon ({\bm{x}},u^{\varepsilon})$, which strongly depends on temperature $u^\varepsilon$, shall be expanded around leading term $u_0$ by displacing ${x_0} = \bm{y}$, ${y_0} = {u_0}$ and $\delta = \varepsilon {u_1} + {\varepsilon^2}{u_2} + O({\varepsilon^3})$ in preceding Taylor's formula as below
\begin{equation}
\begin{aligned}
&\quad\; k_{ij}^\varepsilon ({\bm{x}},u^{\varepsilon})= {k_{ij}}({\bm{y}},u^{\varepsilon}) = {k_{ij}}({\bm{y}},{u_{0}} + \varepsilon{u_{1}} + {\varepsilon^2}{u_{2}} + {\rm O}({\varepsilon^3}))\\
& = {k_{ij}}({\bm{y}},{u_{0}}) + \mathbf{D}^{(0,1)}{k_{ij}}({\bm{y}},{u_{0}})\big[ {\varepsilon{u_{1}} + {\varepsilon^2}{u_{2}} + {\rm O}({\varepsilon^3})} \big]\\
&+ \frac{1}{2}\mathbf{D}^{(0,2)}{k_{ij}}({\bm{y}},{u_{0}}){\big[ {\varepsilon{u_{1}} + {\varepsilon^2}{u_{2}} + {\rm O}({\varepsilon^3})} \big]^2} + {\rm O}\Big( \big[ {\varepsilon{u_{1}} + {\varepsilon^2}{u_{2}} + {\rm O}({\varepsilon^3})} \big]^3 \Big)\\
& = {k_{ij}}({\bm{y}},{u_{0}}) + \varepsilon{u_{1}}\mathbf{D}^{(0,1)}{k_{ij}}({\bm{y}},{u_{0}})\\
&+ {\varepsilon^2}\big[ {{u_{2}}\mathbf{D}^{(0,1)}{k_{ij}}({\bm{y}},{u_{0}}) + \frac{1}{2}{{( {{u_{1}}} )}^2}\mathbf{D}^{(0,2)}{k_{ij}}({\bm{y}},{u_{0}})} \big] + {\rm O}({\varepsilon^3})\\
&= k_{ij}^{(0)}({\bm{y}},{u_{0}})+ \varepsilon k_{ij}^{(1)}({\bm{x}},{\bm{y}},{u_{0}})+ {\varepsilon^2}k_{ij}^{(2)}({\bm{x}},{\bm{y}},{u_{0}})+ {\rm O}({\varepsilon^3}).
\end{aligned}
\end{equation}
Making use of the aforementioned expanding approach as (2.2), remaining material parameters $\rho^{\varepsilon} ({\bm{x}},{u^{\varepsilon}})$, $c^{\varepsilon} ({\bm{x}},{u^{\varepsilon}})$ and ${\sigma_{ij}^{\varepsilon}}({\bm{x}},{u^{\varepsilon}})$ can be continuously expanded as below
\begin{equation}
\begin{aligned}
&{\rho^{\varepsilon}}({\bm{x}},{u^{\varepsilon}})= {\rho^{(0)}}({\bm{y}},{u_{0}}) + \varepsilon{\rho^{(1)}}({\bm{x}},{\bm{y}},{u_{0}}) + {\varepsilon^2}{\rho^{(2)}}({\bm{x}},{\bm{y}},{u_{0}}) + {\rm O}({\varepsilon^3}),\\
&{c^{\varepsilon}}({\bm{x}},{u^{\varepsilon}})= {c^{(0)}}({\bm{y}},{u_{0}}) + \varepsilon{c^{(1)}}({\bm{x}},{\bm{y}},{u_{0}}) + {\varepsilon^2}{c^{(2)}}({\bm{x}},{\bm{y}},{u_{0}}) + {\rm O}({\varepsilon^3}),\\
&\sigma_{ij}^{\varepsilon}({\bm{x}},{u^{\varepsilon}})= \sigma_{ij}^{(0)}({\bm{y}},{u_{0}}) + \varepsilon\sigma_{ij}^{(1)}({\bm{x}},{\bm{y}},{u_{0}}) + {\varepsilon^2}\sigma_{ij}^{(2)}({\bm{x}},{\bm{y}},{u_{0}}) + {\rm O}({\varepsilon ^3}).
\end{aligned}
\end{equation}

Considering the assumption of microscopic periodicity of composite structures, the chain rule achieves for multiscale asymptotic analysis as below.
\begin{equation}
\displaystyle\frac{\partial \Phi^{\varepsilon} (\bm{x},t)}{\partial x_i}=\frac{\partial \Phi (\bm{x},\bm{y},t)}{\partial x_i}+\frac{1}{\varepsilon}\frac{\partial \Phi (\bm{x},\bm{y},t)}{\partial y_i},\;\;(i=1,\cdots,n).
\end{equation}
In preceding formula (2.4), $\Phi^{\varepsilon}(\bm{x},t)$ stands for a function, which has multiscale characteristic, namely $\Phi^{\varepsilon}(\bm{x},t)=\Phi^{\varepsilon}(\bm{x},\bm{y},t)$, referring to \cite{R16}.

Then substituting (2.1)-(2.3) into two-scale nonlinear initial-boundary problem (1.1) and utilizing the chain rule provided in (2.4), we hence have a series of equations by grouping the power-like terms of small periodic parameter $\varepsilon$ as below
\begin{equation}
{\rm O}({\varepsilon^{ - 2}}):\left\{ \begin{aligned}
&\frac{\partial }{{\partial {y_i}}}\Big( {{k_{ij}^{(0)}}\frac{\partial u_{0}}{\partial {y_j}}}\Big)=-{\sigma_{ij}^{(0)}}\frac{\partial \phi_{0}}{\partial {y_i}}\frac{\partial \phi_{0}}{\partial {y_j}},\\
&\frac{\partial }{{\partial {y_i}}}\Big({{\sigma_{ij}^{(0)}}\frac{\partial \phi_{0}}{\partial {y_j}}}\Big)=0.
\end{aligned} \right.
\end{equation}
\begin{equation}
{\rm O}({\varepsilon^{ - 1}}):\left\{ \begin{aligned}
&\frac{\partial }{{\partial {y_i}}}\Big({{k_{ij}^{(0)}}\frac{\partial u_{0}}{\partial {x_j}}}\Big)+\frac{\partial }{{\partial {x_i}}}\Big({{k_{ij}^{(0)}}\frac{\partial u_{0}}{\partial {y_j}}}\Big)+\frac{\partial }{{\partial {y_i}}}\Big({{k_{ij}^{(0)}}\frac{\partial u_{1}}{\partial {y_j}}}\Big)\\
&+\frac{\partial }{{\partial {y_i}}}\Big({{k_{ij}^{(1)}}\frac{\partial u_{0}}{\partial {y_j}}}\Big)=-{\sigma_{ij}^{(0)}}\frac{\partial \phi_{0}}{\partial {x_i}}\frac{\partial \phi_{0}}{\partial {y_j}}-{\sigma_{ij}^{(0)}}\frac{\partial \phi_{0}}{\partial {y_i}}\frac{\partial \phi_{0}}{\partial {x_j}}\\
&-{\sigma_{ij}^{(1)}}\frac{\partial \phi_{0}}{\partial {y_i}}\frac{\partial \phi_{0}}{\partial {y_j}}-{\sigma_{ij}^{(0)}}\frac{\partial \phi_{0}}{\partial {y_i}}\frac{\partial \phi_{1}}{\partial {y_j}}-{\sigma_{ij}^{(0)}}\frac{\partial \phi_{1}}{\partial {y_i}}\frac{\partial \phi_{0}}{\partial {y_j}},\\
&\frac{\partial }{{\partial {y_i}}}\Big( {{\sigma_{ij}^{(0)}}\frac{\partial \phi_{0}}{\partial {x_j}}}\Big)+\frac{\partial }{{\partial {x_i}}}\Big({{\sigma_{ij}^{(0)}}\frac{\partial \phi_{0}}{\partial {y_j}}}\Big)+\frac{\partial }{{\partial {y_i}}}\Big({{\sigma_{ij}^{(0)}}\frac{\partial \phi_{1}}{\partial {y_j}}}\Big)\\
&+\frac{\partial }{{\partial {y_i}}}\Big({{\sigma_{ij}^{(1)}}\frac{\partial \phi_{0}}{\partial {y_j}}}\Big)=0.
\end{aligned} \right.
\end{equation}
\begin{equation}
{\rm O}({\varepsilon^0}):\left\{ \begin{aligned}
&\rho^{(0)}c^{(0)}\frac{{\partial {u_{0}}}}{{\partial t}}-\frac{\partial }{{\partial {x_i}}}\Big( {{k_{ij}^{(0)}}\frac{{\partial u_{0} }}{{\partial {x_j}}}} \Big)-\frac{\partial }{{\partial {x_i}}}\Big( {{k_{ij}^{(0)}}\frac{{\partial u_{1} }}{{\partial {y_j}}}} \Big)-\frac{\partial }{{\partial {y_i}}}\Big({{k_{ij}^{(0)}}\frac{{\partial u_{1} }}{{\partial {x_j}}}} \Big)\\
&-\frac{\partial }{{\partial {y_i}}}\Big({{k_{ij}^{(0)}}\frac{{\partial u_{2} }}{{\partial {y_j}}}} \Big)-\frac{\partial }{{\partial {x_i}}}\Big({{k_{ij}^{(1)}}\frac{{\partial u_{0} }}{{\partial {y_j}}}} \Big)-\frac{\partial }{{\partial {y_i}}}\Big({{k_{ij}^{(1)}}\frac{{\partial u_{0} }}{{\partial {x_j}}}} \Big)\\
&-\frac{\partial }{{\partial {y_i}}}\Big({{k_{ij}^{(1)}}\frac{{\partial u_{1} }}{{\partial {y_j}}}} \Big)-\frac{\partial }{{\partial {y_i}}}\Big({{k_{ij}^{(2)}}\frac{{\partial u_{0} }}{{\partial {y_j}}}} \Big)={\sigma_{ij}^{(0)}}\frac{\partial \phi_{0}}{\partial {x_i}}\frac{\partial \phi_{0}}{\partial {x_j}}+{\sigma_{ij}^{(0)}}\frac{\partial \phi_{0}}{\partial {x_i}}\frac{\partial \phi_{1}}{\partial {y_j}}\\
&+{\sigma_{ij}^{(0)}}\frac{\partial \phi_{1}}{\partial {x_i}}\frac{\partial \phi_{0}}{\partial {y_j}}+{\sigma_{ij}^{(0)}}\frac{\partial \phi_{0}}{\partial {y_i}}\frac{\partial \phi_{1}}{\partial {x_j}}+{\sigma_{ij}^{(0)}}\frac{\partial \phi_{0}}{\partial {y_i}}\frac{\partial \phi_{2}}{\partial {y_j}}+{\sigma_{ij}^{(0)}}\frac{\partial \phi_{1}}{\partial {y_i}}\frac{\partial \phi_{0}}{\partial {x_j}}\\
&+{\sigma_{ij}^{(0)}}\frac{\partial \phi_{1}}{\partial {y_i}}\frac{\partial \phi_{1}}{\partial {y_j}}+{\sigma_{ij}^{(0)}}\frac{\partial \phi_{2}}{\partial {y_i}}\frac{\partial \phi_{0}}{\partial {y_j}}+{\sigma_{ij}^{(1)}}\frac{\partial \phi_{0}}{\partial {x_i}}\frac{\partial \phi_{0}}{\partial {y_j}}+{\sigma_{ij}^{(1)}}\frac{\partial \phi_{0}}{\partial {y_i}}\frac{\partial \phi_{0}}{\partial {x_j}}\\
&+{\sigma_{ij}^{(1)}}\frac{\partial \phi_{0}}{\partial {y_i}}\frac{\partial \phi_{1}}{\partial {y_j}}+{\sigma_{ij}^{(1)}}\frac{\partial \phi_{1}}{\partial {y_i}}\frac{\partial \phi_{0}}{\partial {y_j}}+{\sigma_{ij}^{(2)}}\frac{\partial \phi_{0}}{\partial {y_i}}\frac{\partial \phi_{0}}{\partial {y_j}}+f_u,\\
&\frac{\partial }{{\partial {x_i}}}\Big({{\sigma_{ij}^{(0)}}\frac{{\partial \phi_{0} }}{{\partial {x_j}}}} \Big)+\frac{\partial }{{\partial {x_i}}}\Big( {{\sigma_{ij}^{(0)}}\frac{{\partial \phi_{1} }}{{\partial {y_j}}}} \Big)+\frac{\partial }{{\partial {y_i}}}\Big({{\sigma_{ij}^{(0)}}\frac{{\partial \phi_{1} }}{{\partial {x_j}}}} \Big)\\
&+\frac{\partial }{{\partial {y_i}}}\Big({{\sigma_{ij}^{(0)}}\frac{{\partial \phi_{2} }}{{\partial {y_j}}}} \Big)+\frac{\partial }{{\partial {x_i}}}\Big({{\sigma_{ij}^{(1)}}\frac{{\partial \phi_{0} }}{{\partial {y_j}}}} \Big)+\frac{\partial }{{\partial {y_i}}}\Big({{\sigma_{ij}^{(1)}}\frac{{\partial \phi_{0} }}{{\partial {x_j}}}} \Big)\\
&+\frac{\partial }{{\partial {y_i}}}\Big( {{\sigma_{ij}^{(1)}}\frac{{\partial \phi_{1} }}{{\partial {y_j}}}} \Big)+\frac{\partial }{{\partial {y_i}}}\Big( {{\sigma_{ij}^{(2)}}\frac{{\partial \phi_{0} }}{{\partial {y_j}}}} \Big)=-f_\phi.
\end{aligned} \right.
\end{equation}
%In the subsequent part, each equation with order $O({\varepsilon^{-2}})$, $O({\varepsilon^{- 1}})$ and $O({\varepsilon^{0}})$ is analyzed respectively so that the higher-order multiscale analysis can be implemented for deriving the detailed expressions of multiscale asymptotic solutions.

From $O({\varepsilon^{-2}})$-order equations (2.5), we firstly can deduce that
\begin{equation}
\begin{aligned}
u_{0}(\bm{x},\bm{y},t) =u_{0}(\bm{x},t),\;\phi_{0}(\bm{x},\bm{y},t)= \phi_{0}(\bm{x},t).
\end{aligned}
\end{equation}
By noticing the fact (2.8), the terms $\displaystyle\frac{\partial u_{0}}{\partial y_j}$ and $\displaystyle\frac{\partial \phi_{0}}{\partial y_j}$ in $O({\varepsilon^{-1}})$-order equations (2.6) both equate to zero. Subsequently, equations (2.6) can be further simplified as the subsequent equations with $\alpha_1=1,\cdots,n$.
\begin{equation}
\left\{ \begin{aligned}
&\frac{\partial }{{\partial {y_i}}}\Big({{k_{ij}^{(0)}}\frac{\partial u_{1}}{\partial {y_j}}}\Big)=-\frac{\partial {k_{i\alpha_1}^{(0)}}}{{\partial {y_i}}}{\frac{\partial u_{0}}{\partial {x_{\alpha_1}}}},\\
&\frac{\partial }{{\partial {y_i}}}\Big({{\sigma_{ij}^{(0)}}\frac{\partial \phi_{1}}{\partial {y_j}}}\Big)=-\frac{\partial {\sigma_{i\alpha_1}^{(0)}}}{{\partial {y_i}}}{\frac{\partial \phi_{0}}{\partial {x_{\alpha_1}}}}.
\end{aligned} \right.
\end{equation}
Taking advantage of equations (2.9) and using the separation of variables, the first-order correctors $u_{1}$ and $\phi_{1}$ can be decomposed as the following separation forms
\begin{equation}
\left\{
\begin{aligned}
&{u_{1}}(\bm{x},\bm{y},t) = {M_{\alpha_1}}(\bm{y},u_{0})\frac{\partial u_{0}(\bm{x},t)}{\partial x_{\alpha_1}},\\
&\phi_{1}(\bm{x},\bm{y},t) = N_{\alpha_1}(\bm{y},u_{0})\frac{\partial \phi_{0}(\bm{x},t)}{\partial x_{\alpha_1}},
\end{aligned}\right.
\end{equation}
where ${M_{\alpha_1}}$ and $N_{\alpha_1}$ are termed as first-order cell functions defined in PUC $\Theta$, that all rely upon macroscopic temperature field $u_0$. Furthermore, substituting (2.10) into (2.9), the following equations with homogeneous Dirichlet boundary condition are obtained after simplified calculation, which are referred as the unit cell problems.
\begin{equation}
\left\{
\begin{aligned}
&\frac{\partial}{\partial y_i}\big[ { k_{ij}^{(0)}{\frac{\partial M_{\alpha_1}}{\partial y_j}}} \big]= -\frac{\partial k_{i{\alpha_1}}^{(0)}}{\partial y_i},\;\;\;\bm{y}\in \Theta, \\
&M_{\alpha_1}(\bm{y},u_{0})\;\mathrm{is}\;1-\mathrm{periodic}\;\mathrm{in}\;{\bm{y}}{\rm{, }}\;\;\;{\int_{\Theta}}M_{\alpha_1}d\Theta=0.
\end{aligned} \right.
\end{equation}
\begin{equation}
\left\{
\begin{aligned}
&\frac{\partial}{\partial y_i}\big[ { \sigma_{ij}^{(0)}{\frac{\partial N_{\alpha_1}}{\partial y_j}}} \big]= -\frac{\partial \sigma_{i{\alpha_1}}^{(0)}}{\partial y_i},\;\;\;\bm{y}\in \Theta, \\
&N_{\alpha_1}(\bm{y},u_{0})\;\mathrm{is}\;1-\mathrm{periodic}\;\mathrm{in}\;{\bm{y}}{\rm{, }}\;\;\;{\int_{\Theta}}N_{\alpha_1}d\Theta=0.
\end{aligned} \right.
\end{equation}

Subsequently, performing a volume integral on both sides of equations (2.7) on microscopic unit cell $\Theta$ with the substitution (2.10) into (2.7) and exploiting the Gauss theorem on equations (2.7), these procedures lead to derive the macroscopic homogenized equations associated with multi-scale problem (1.1) as presented below
\begin{equation}
\left\{ \begin{aligned}
&\widehat S(u_{0})\frac{{\partial {u_{0}}(\bm{x},t)}}{{\partial t}}- \frac{\partial }{{\partial {x_i}}}\Big( {{\widehat k_{ij}}(u_{0})\frac{{\partial {u_{0}}(\bm{x},t)}}{{\partial {x_j}}}}\Big)\\
&\quad\quad\quad\quad\quad\quad = {\widehat \sigma _{ij}^*}(u_{0})\frac{\partial \phi_{0}(\bm{x},t)}{\partial {x_i}}\frac{\partial \phi_{0}(\bm{x},t)}{\partial {x_j}}+f_u(\bm{x},t),\;\;\text{in}\;\;\Omega\times(0,T),\\
&- \frac{\partial }{{\partial {x_i}}}\Big( {{\widehat \sigma _{ij}}(u_{0})\frac{{\partial \phi_0(\bm{x},t)}}{{\partial {x_j}}}}\Big)=f_\phi(\bm{x},t),\;\;\text{in}\;\;\Omega\times(0,T),\\
&u_0(\bm{x},t) = \widehat u(\bm{x},t),\;\;\text{on}\;\;\partial {\Omega_u}\times(0,T),\\
&\phi_0(\bm{x},t)=\widehat\phi(\bm{x},t),\;\;\text{on}\;\;\partial\Omega_{\phi}\times(0,T),\\
&{{\widehat k_{ij}}(u_{0})\frac{\partial u_0(\bm{x},t)}{\partial {x_j}}}{n_i} = \bar q(\bm{x},t),\;\;\text{on}\;\;\partial {\Omega_q}\times(0,T),\\
&{{\widehat \sigma _{ij}}(u_{0})\frac{\partial \phi_0(\bm{x},t)}{\partial {x_j}}}{n_i}=\bar{d}(\bm{x},t),\;\;\text{on}\;\;\partial\Omega_{d}\times(0,T),\\
&u_0({\bm{x}},0)=\widetilde u,\;\;\text{in}\;\;\Omega.
\end{aligned} \right.
\end{equation}
Here, the macroscopic homogenized material parameters in (2.13) are evaluated using the following formulas, which correspond to microscopic unit cell $\Theta$.
\begin{equation}
\begin{aligned}
&\widehat S(u_{0}) = \frac{1}{|\Theta|}{\int_{\Theta}}{\rho^{(0)}{c}^{(0)}}d\Theta,\;{\widehat k_{ij}}(u_{0}) = \frac{1}{|\Theta|}{\int_{\Theta}}\big({k_{ij}^{(0)} + k_{i\alpha_1}^{(0)}{\frac{\partial M_j}{\partial y_{\alpha_1}}}}\big)d\Theta,\\
&{\widehat \sigma _{ij}}(u_{0}) = \frac{1}{|\Theta|}{\int_{\Theta}}\big({\sigma_{ij}^{(0)}+ \sigma_{i\alpha_1}^{(0)}{\frac{\partial N_j}{\partial y_{\alpha_1}}}}\big)d\Theta,\\
&{\widehat \sigma _{ij}^*}(u_{0}) = \frac{1}{|\Theta|}{\int_{\Theta}}\big({\sigma_{ij}^{(0)}+ \sigma_{i\alpha_1}^{(0)}{\frac{\partial N_j}{\partial y_{\alpha_1}}}}+ {\sigma_{\alpha_1j}^{(0)}{\frac{\partial N_i}{\partial y_{\alpha_1}}}}+ \sigma_{\alpha_1\alpha_2}^{(0)}{\frac{\partial N_i}{\partial y_{\alpha_1}}{\frac{\partial N_j}{\partial y_{\alpha_2}}}}\big)d\Theta.
\end{aligned}
\end{equation}
\begin{rmk}
It is to be noted that all homogenized material parameters vary with the macroscopic homogenized solution $u_{0}$ due to the quasi-periodic properties of first-order cell functions. This distinction is significant comparing with linear composites.
\end{rmk}
\begin{rmk}
Stemming from the approach as outlined in \cite{R16,R17}, it can be proved in Section 5 that $\bar\gamma_0|\bm{\zeta}|^2\leq \widehat k_{ij}(u_{0})\zeta_i\zeta_j\leq\bar\gamma_1|\bm{\zeta}|^2$ and $\bar\gamma_0|\bm{\zeta}|^2\leq \widehat \sigma_{ij}(u_{0})\zeta_i\zeta_j\leq\bar\gamma_1|\bm{\zeta}|^2$.
\end{rmk}
\begin{rmk}
Despite ${\widehat \sigma _{ij}^*}(u_{0})$ and ${\widehat \sigma _{ij}}(u_{0})$ have different computational formulas, we can prove that ${\widehat \sigma _{ij}^*}(u_{0})={\widehat \sigma _{ij}}(u_{0})$ for any fixed $u_{0}$. The detailed proof is provided in Appendix A.
\end{rmk}

Now, we proceed to establish the vital second-order correctors $u_{2}$ and $\phi_{2}$. By substituting (2.8) and (2.10) into (2.7), and then subtracting (2.7) from (2.13), we can formulate the following equations after simplification and computation
\begin{equation}
\left\{\begin{aligned}
&\!\!\!\frac{\partial}{\partial y_i}\Big( {k_{ij}^{(0)}{\frac{\partial u_{2}}{\partial y_j}}} \Big)=\Big[ { \rho^{(0)}{c}^{(0)}-\widehat S(u_{0}) } \Big]\frac{{\partial {u_{0}}}}{{\partial t}}\\
&\!\!\!+ \Big[ { {\widehat k}_{\alpha_1\alpha_2}-  {k_{\alpha_1\alpha_2}^{(0)}}}{- \frac{\partial}{\partial y_i}\big( {k_{i\alpha_1}^{(0)}{M_{\alpha_2}}} \big)-{k_{\alpha_1j}\frac{\partial M_{\alpha_2}}{\partial y_j}}} \Big]\frac{\partial^2 u_{0}}{\partial x_{\alpha_1}\partial x_{\alpha_2}}\\
&\!\!\!+ \Big[ { \frac{\partial{\widehat k}_{i\alpha_1}}{\partial x_i} -  {\frac{\partial{k}_{i\alpha_1}^{(0)}}{\partial x_i}}- \frac{\partial}{\partial y_i}\big( {k_{ij}^{(0)}\frac{\partial M_{\alpha_1}}{\partial x_{j}}} \big)}{-\frac{\partial}{\partial x_i}\big( {k_{ij}^{(0)}\frac{\partial M_{\alpha_1}}{\partial y_{j}}} \big)} \Big]\frac{\partial u_{0}}{\partial x_{\alpha_1}}\\
&\!\!\!- \frac{\partial}{\partial y_i}\Big({{M_{\alpha_1}}\mathbf{D}^{(0,1)}{k_{i\alpha_2}^{(0)}}}+ {{M_{\alpha_1}}\mathbf{D}^{(0,1)}{k_{ij}^{(0)}}\frac{\partial M_{\alpha_2}}{\partial y_j}}\Big)\frac{\partial u_{0}}{\partial x_{\alpha_1}}\frac{\partial u_{0}}{\partial x_{\alpha_2}}\\
&\!\!\!+\!\Big[{\widehat \sigma}_{\alpha_1\alpha_2}^* - \sigma_{\alpha_1\alpha_2}^{(0)}
-\sigma_{i\alpha_2}^{(0)}\frac{\partial N_{\alpha_1}}{\partial y_i}- \sigma_{\alpha_1j}^{(0)}\frac{\partial N_{\alpha_2}}{\partial y_{j}}-\sigma_{ij}^{(0)}{\frac{\partial N_{\alpha_1}}{\partial y_{i}}{\frac{\partial N_{\alpha_2}}{\partial y_j}}}\Big]\frac{\partial \phi_{0}}{\partial x_{\alpha_1}}\frac{\partial \phi_{0}}{\partial x_{\alpha_2}},\\
&\!\!\!\frac{\partial}{\partial y_i}\Big( {\sigma_{ij}^{(0)}{\frac{\partial \phi_{2}}{\partial y_j}}} \Big)=\Big[ { {\widehat \sigma}_{\alpha_1\alpha_2}-  {\sigma_{\alpha_1\alpha_2}^{(0)}}}{- \frac{\partial}{\partial y_i}\big( {\sigma_{i\alpha_1}^{(0)}{N_{\alpha_2}}} \big)-{\sigma_{\alpha_1j}\frac{\partial N_{\alpha_2}}{\partial y_j}}} \Big]\frac{\partial^2 \phi_{0}}{\partial x_{\alpha_1}\partial x_{\alpha_2}}\\
&\!\!\!+ \Big[ { \frac{\partial{\widehat \sigma}_{i\alpha_1}}{\partial x_i} -  {\frac{\partial{\sigma}_{i\alpha_1}^{(0)}}{\partial x_i}}- \frac{\partial}{\partial y_i}\big( {\sigma_{ij}^{(0)}\frac{\partial N_{\alpha_1}}{\partial x_{j}}} \big)}{-\frac{\partial}{\partial x_i}\big( {\sigma_{ij}^{(0)}\frac{\partial N_{\alpha_1}}{\partial y_{j}}} \big)} \Big]\frac{\partial \phi_{0}}{\partial x_{\alpha_1}}\\
&\!\!\!- \frac{\partial}{\partial y_i}\Big({{M_{\alpha_1}}\mathbf{D}^{(0,1)}{\sigma_{i\alpha_2}^{(0)}}}+ {{M_{\alpha_1}}\mathbf{D}^{(0,1)}{\sigma_{ij}^{(0)}}\frac{\partial N_{\alpha_2}}{\partial y_j}}\Big)\frac{\partial u_{0}}{\partial x_{\alpha_1}}\frac{\partial \phi_{0}}{\partial x_{\alpha_2}}.
\end{aligned}\right.
\end{equation}
Given equations (2.15), then we construct the concrete expressions with $\alpha_2=1,\cdots,n$ for ${u_{2}}$ and $\phi_{2}$ as follows
\begin{equation}
\left\{
\begin{aligned}
{u_{2}}(\bm{x},\bm{y},t)&= Q({\bm{y}},{u_{0}})\frac{{\partial {u_{0}(\bm{x},t)}}}{{\partial t}}+{M_{\alpha_1\alpha_2}}({\bm{y}},{u_{0}})\frac{\partial^2 u_{0}(\bm{x},t)}{\partial x_{\alpha_1}\partial x_{\alpha_2}}\\
&+R_{\alpha_1}({\bm{y}},{u_{0}})\frac{\partial u_{0}(\bm{x},t)}{\partial x_{\alpha_1}}+{H_{\alpha_1\alpha_2}}({\bm{y}},{u_{0}})\frac{\partial u_{0}(\bm{x},t)}{\partial x_{\alpha_1}}\frac{\partial u_{0}(\bm{x},t)}{\partial x_{\alpha_2}}\\
&+{G_{\alpha_1\alpha_2}}({\bm{y}},{u_{0}})\frac{\partial \phi_{0}(\bm{x},t)}{\partial x_{\alpha_1}}\frac{\partial \phi_{0}(\bm{x},t)}{\partial x_{\alpha_2}},\\
\phi_{2}(\bm{x},\bm{y},t)&= N_{\alpha_1\alpha_2}({\bm{y}},{u_{0}})\frac{\partial^2 \phi_{0}(\bm{x},t)}{\partial x_{\alpha_1}\partial x_{\alpha_2}}+ Z_{\alpha_1}({\bm{y}},{u_{0}})\frac{\partial \phi_{0}(\bm{x},t)}{\partial x_{\alpha_1}}\\
&+W_{\alpha_1\alpha_2}({\bm{y}},{u_{0}})\frac{\partial u_{0}(\bm{x},t)}{\partial x_{\alpha_1}}\frac{\partial \phi_{0}(\bm{x},t)}{\partial x_{\alpha_2}}.
\end{aligned} \right.
\end{equation}
In the above expressions, $Q$, $M_{\alpha_1\alpha_2}$, $R_{\alpha_1}$, ${H_{\alpha_1\alpha_2}}$, $G_{\alpha_1\alpha_2}$, $N_{\alpha_1\alpha_2}$, $Z_{\alpha_1}$ and $W_{\alpha_1\alpha_2}$ are referred as second-order cell functions, which all rely upon macroscopic temperature field $u_0$. By substituting (2.16) into (2.15), a series of equations, which are subject to the homogeneous Dirichlet boundary condition, are derived as follows
\begin{equation}
\left\{
\begin{aligned}
&\frac{\partial}{\partial y_i}\big[ { k_{ij}^{(0)}{\frac{\partial Q}{\partial y_j}}} \big] = { \rho^{(0)} {c}^{(0)}-\widehat S },\;\;\;\bm{y}\in \Theta,\\
&Q({\bm{y}},{u_{0}})\;\mathrm{is}\;1-\mathrm{periodic}\;\mathrm{in}\;{\bm{y}}{\rm{, }}\;\;\;{\int_{\Theta}}Qd\Theta=0.
\end{aligned} \right.
\end{equation}
\begin{equation}
\left\{
\begin{aligned}
&\frac{\partial}{\partial y_i}\big[ { k_{ij}^{(0)}{\frac{\partial {M_{\alpha_1\alpha_2}}}{\partial y_j}}} \big]= { {\widehat k}_{\alpha_1\alpha_2}-  {k_{\alpha_1\alpha_2}^{(0)}}\!\!- \frac{\partial}{\partial y_i}\big( {k_{i\alpha_1}^{(0)}{M_{\alpha_2}}} \big)\!\!-{k_{\alpha_1j}^{(0)}\frac{\partial M_{\alpha_2}}{\partial y_j}}},\;\;\;\bm{y}\in \Theta,\\
&{M_{\alpha_1\alpha_2}}({\bm{y}},{u_{0}})\;\mathrm{is}\;1-\mathrm{periodic}\;\mathrm{in}\;{\bm{y}}{\rm{, }}\;\;\;{\int_{\Theta}}M_{\alpha_1\alpha_2}d\Theta=0.
\end{aligned} \right.
\end{equation}
\begin{equation}
\left\{
\begin{aligned}
&\frac{\partial}{\partial y_i}\big[ { k_{ij}^{(0)}{\frac{\partial R_{\alpha_1}}{\partial y_j}}} \big]\!\!=\!\!{ \frac{\partial{\widehat k}_{i\alpha_1}}{\partial x_i} \!\!-\!\!  {\frac{\partial{k}_{i\alpha_1}^{(0)}}{\partial x_i}}\!\!-\!\!\frac{\partial}{\partial y_i}\big( {k_{ij}^{(0)}\frac{\partial M_{\alpha_1}}{\partial x_{j}}} \big)\!\!-\!\!\frac{\partial}{\partial x_i}\big( {k_{ij}^{(0)}\frac{\partial M_{\alpha_1}}{\partial y_{j}}} \big)},\;\;\;\bm{y}\!\in\! \Theta,\\
&{R_{\alpha_1}}({\bm{y}},{u_{0}})\;\mathrm{is}\;1-\mathrm{periodic}\;\mathrm{in}\;{\bm{y}}{\rm{, }}\;\;\;{\int_{\Theta}}R_{\alpha_1}d\Theta=0.
\end{aligned} \right.
\end{equation}
\begin{equation}
\left\{
\begin{aligned}
&\frac{\partial}{\partial y_i}\big[ { k_{ij}^{(0)}{\frac{\partial H_{\alpha_1\alpha_2}}{\partial y_j}}} \big]\!\!=\!\!-\frac{\partial}{\partial y_i}\Big({{M_{\alpha_1}}\mathbf{D}^{(0,1)}{k_{i\alpha_2}^{(0)}}}\!+ \!{{M_{\alpha_1}}\mathbf{D}^{(0,1)}{k_{ij}^{(0)}}\frac{\partial M_{\alpha_2}}{\partial y_j}}\Big),\;\;\;\bm{y}\!\in\!\Theta,\\
&{H_{\alpha_1\alpha_2}}({\bm{y}},{u_{0}})\;\mathrm{is}\;1-\mathrm{periodic}\;\mathrm{in}\;{\bm{y}}{\rm{, }}\;\;\;{\int_{\Theta}}H_{\alpha_1\alpha_2}d\Theta=0.
\end{aligned} \right.
\end{equation}
\begin{equation}
\left\{
\begin{aligned}
&\frac{\partial}{\partial y_i}\big[ { k_{ij}^{(0)}{\frac{\partial {G_{\alpha_1\alpha_2}}}{\partial y_j}}} \big] = {\widehat \sigma}_{\alpha_1\alpha_2}^* - \sigma_{\alpha_1\alpha_2}^{(0)}
-\sigma_{i\alpha_2}^{(0)}\frac{\partial N_{\alpha_1}}{\partial y_i}\\
&\quad\quad\quad\quad\quad\quad\quad\;\;- \sigma_{\alpha_1j}^{(0)}\frac{\partial N_{\alpha_2}}{\partial y_{j}}-\sigma_{ij}^{(0)}{\frac{\partial N_{\alpha_1}}{\partial y_{i}}{\frac{\partial N_{\alpha_2}}{\partial y_j}}},\;\;\;\bm{y}\in \Theta,\\
&{G_{\alpha_1\alpha_2}}({\bm{y}},{u_{0}})\;\mathrm{is}\;1-\mathrm{periodic}\;\mathrm{in}\;{\bm{y}}{\rm{, }}\;\;\;{\int_{\Theta}}G_{\alpha_1\alpha_2}d\Theta=0.
\end{aligned} \right.
\end{equation}
\begin{equation}
\left\{
\begin{aligned}
&\frac{\partial}{\partial y_i}\big[ { \sigma_{ij}^{(0)}{\frac{\partial {N_{\alpha_1\alpha_2}}}{\partial y_j}}} \big]\!\!=\!\!{ {\widehat \sigma}_{\alpha_1\alpha_2}\!-\!  {\sigma_{\alpha_1\alpha_2}^{(0)}}\!-\!\frac{\partial}{\partial y_i}\big( {\sigma_{i\alpha_1}^{(0)}{N_{\alpha_2}}} \big)\!-\!{\sigma_{\alpha_1j}^{(0)}\frac{\partial N_{\alpha_2}}{\partial y_j}}},\;\;\;\bm{y}\in \Theta,\\
&{N_{\alpha_1\alpha_2}}({\bm{y}},{u_{0}})\;\mathrm{is}\;1-\mathrm{periodic}\;\mathrm{in}\;{\bm{y}}{\rm{, }}\;\;\;{\int_{\Theta}}N_{\alpha_1\alpha_2}d\Theta=0.
\end{aligned} \right.
\end{equation}
\begin{equation}
\left\{
\begin{aligned}
&\frac{\partial}{\partial y_i}\big[ { \sigma_{ij}^{(0)}{\frac{\partial Z_{\alpha_1}}{\partial y_j}}} \big]\!\!=\!\!{ \frac{\partial{\widehat \sigma}_{i\alpha_1}}{\partial x_i}\!\!-\!\!  {\frac{\partial{\sigma}_{i\alpha_1}^{(0)}}{\partial x_i}}\!\!-\!\!\frac{\partial}{\partial y_i}\big( {\sigma_{ij}^{(0)}\frac{\partial N_{\alpha_1}}{\partial x_{j}}} \big)\!\!-\!\!\frac{\partial}{\partial x_i}\big( {\sigma_{ij}^{(0)}\frac{\partial N_{\alpha_1}}{\partial y_{j}}} \big)},\;\;\;\bm{y}\!\in\!\Theta,\\
&{Z_{\alpha_1}}({\bm{y}},{u_{0}})\;\mathrm{is}\;1-\mathrm{periodic}\;\mathrm{in}\;{\bm{y}}{\rm{, }}\;\;\;{\int_{\Theta}}Z_{\alpha_1}d\Theta=0.
\end{aligned} \right.
\end{equation}
\begin{equation}
\left\{
\begin{aligned}
&\frac{\partial}{\partial y_i}\big[ { \sigma_{ij}^{(0)}{\frac{\partial W_{\alpha_1\alpha_2}}{\partial y_j}}} \big]\!\!=\!\!-\frac{\partial}{\partial y_i}\Big({{M_{\alpha_1}}\mathbf{D}^{(0,1)}{\sigma_{i\alpha_2}^{(0)}}}\!+\! {{M_{\alpha_1}}\mathbf{D}^{(0,1)}{\sigma_{ij}^{(0)}}\frac{\partial N_{\alpha_2}}{\partial y_j}}\Big),\;\;\;\bm{y}\!\in\!\Theta,\\
&{H_{\alpha_1\alpha_2}}({\bm{y}},{u_{0}})\;\mathrm{is}\;1-\mathrm{periodic}\;\mathrm{in}\;{\bm{y}}{\rm{, }}\;\;\;{\int_{\Theta}}H_{\alpha_1\alpha_2}d\Theta=0.
\end{aligned} \right.
\end{equation}
\begin{rmk}
By the hypotheses (A$_1$)-(A$_2$) and Lax-Milgram theorem, the existence and uniqueness of solutions for equations (2.11)-(2.12) and (2.17)-(2.24) are established for any fixed macroscopic temperature field $u_0$.
\end{rmk}

By now, lower-order multiscale (LOMS) solutions for temperature $u^{\varepsilon}$ and electric potential $\phi^{\varepsilon}$ are given by
\begin{equation}
{\begin{aligned}
{u^{(1\varepsilon)} }(\bm{x},t)&=%{\approx}
u_{0}(\bm{x},t)+\varepsilon{M_{\alpha_1}}(\bm{y},u_{0})\frac{\partial u_{0}(\bm{x},t)}{\partial x_{\alpha_1}}.
\end{aligned}}
\end{equation}
\begin{equation}
{\begin{aligned}
\phi^{(1\varepsilon)}(\bm{x},t)&=%{\approx}
\phi_{0}(\bm{x},t)+\varepsilon N_{\alpha_1}(\bm{y},u_{0})\frac{\partial \phi_{0}(\bm{x},t)}{\partial x_{\alpha_1}}.
\end{aligned}}
\end{equation}
In addition, higher-order multiscale (HOMS) solutions for temperature $u^{\varepsilon}$ and electric potential $\phi^{\varepsilon}$ are given by
\begin{equation}
{\begin{aligned}
{u^{(2\varepsilon)} }(\bm{x},t)&\!\!=\!\!%{\approx}
u_{0}(\bm{x},\!t)\!\!+\!\!\varepsilon{M_{\alpha_1}}(\bm{y},\!u_{0})\frac{\partial u_{0}(\bm{x},t)}{\partial x_{\alpha_1}}\!\!+\!\!\varepsilon^2\big[ Q({\bm{y}},\!{u_{0}})\frac{{\partial {u_{0}}(\bm{x},t)}}{{\partial t}}\\
&\!\!+\!\!{M_{\alpha_1\alpha_2}}({\bm{y}},\!{u_{0}})\frac{\partial^2 u_{0}(\bm{x},t)}{\partial x_{\alpha_1}\partial x_{\alpha_2}}\!\!+\!\!R_{\alpha_1}({\bm{y}},\!{u_{0}})\frac{\partial u_{0}(\bm{x},t)}{\partial x_{\alpha_1}}\\
&\!\!+\!\!{H_{\alpha_1\alpha_2}}({\bm{y}},\!{u_{0}})\frac{\partial u_{0}(\bm{x},t)}{\partial x_{\alpha_1}}\frac{\partial u_{0}(\bm{x},t)}{\partial x_{\alpha_2}}\!\!+\!\!{G_{\alpha_1\alpha_2}}({\bm{y}},\!{u_{0}})\frac{\partial \phi_{0}(\bm{x},t)}{\partial x_{\alpha_1}}\frac{\partial \phi_{0}(\bm{x},t)}{\partial x_{\alpha_2}}\big].
\end{aligned}}
\end{equation}
\begin{equation}
{\begin{aligned}
\phi^{(2\varepsilon)}(\bm{x},t)&=%{\approx}
\phi_{0}(\bm{x},t)+\varepsilon N_{\alpha_1}(\bm{y},u_{0})\frac{\partial \phi_{0}(\bm{x},t)}{\partial x_{\alpha_1}}+ \varepsilon^2\big[N_{\alpha_1\alpha_2}({\bm{y}},{u_{0}})\frac{\partial^2 \phi_{0}(\bm{x},t)}{\partial x_{\alpha_1}\partial x_{\alpha_2}}\\
&+ Z_{\alpha_1}({\bm{y}},{u_{0}})\frac{\partial \phi_{0}(\bm{x},t)}{\partial x_{\alpha_1}}+W_{\alpha_1\alpha_2}({\bm{y}},{u_{0}})\frac{\partial u_{0}(\bm{x},t)}{\partial x_{\alpha_1}}\frac{\partial \phi_{0}(\bm{x},t)}{\partial x_{\alpha_2}}\big].
\end{aligned}}
\end{equation}
Moreover, also from (2.27) and (2.28), it can be concluded that only the HOMS solutions can characterize the mutual coupling impact of temperature field and electric potential field owing to the presence of correction terms $\displaystyle {G_{\alpha_1\alpha_2}}({\bm{y}},{u_{0}})\frac{\partial \phi_{0}(\bm{x},t)}{\partial x_{\alpha_1}}\frac{\partial \phi_{0}(\bm{x},t)}{\partial x_{\alpha_2}}$ in (2.27) and $\displaystyle W_{\alpha_1\alpha_2}({\bm{y}},{u_{0}})\frac{\partial u_{0}(\bm{x},t)}{\partial x_{\alpha_1}}\frac{\partial \phi_{0}(\bm{x},t)}{\partial x_{\alpha_2}}$ in (2.28), which is one essential motivation to develop higher-order multiscale method for high-accuracy nonlinear thermo-electric coupling simulation of composite structures.
\section{The error analyses of multiscale asymptotic solutions}
%\overset{*}
\subsection{The proof of local balance of heat flux and electric charge by local error analysis}
Firstly, the residual functions for LOMS and HOMS approximate solutions are defined as below.
\begin{equation}
\left\{
\begin{aligned}
&{u_{\Delta}^{(1\varepsilon)} }(\bm{x},t)={u^{\varepsilon} }(\bm{x},t)-{u^{(1\varepsilon)} }(\bm{x},t),\;{\phi_{\Delta}^{(1\varepsilon)} }(\bm{x},t)={\phi^{\varepsilon} }(\bm{x},t)-{ \phi^{(1\varepsilon)}}(\bm{x},t),\\
&{u_{\Delta}^{(2\varepsilon)} }(\bm{x},t)={u^{\varepsilon} }(\bm{x},t)-{u^{(2\varepsilon)} }(\bm{x},t),\;{\phi_{\Delta}^{(2\varepsilon)} }(\bm{x},t)={\phi^{\varepsilon} }(\bm{x},t)-{ \phi^{(2\varepsilon)}}(\bm{x},t).
\end{aligned}\right.
\end{equation}

To analyze the local heat flux and electric charge balance of multiscale asymptotic solutions, substituting the above residual functions $u^{(1\varepsilon)}_{\Delta}$ and $\phi^{(1\varepsilon)}_{\Delta}$ into (1.1), the residual equations for the LOMS solutions can be derived as below.
\begin{equation}
\left\{
\begin{aligned}
&\rho^{\varepsilon}c^{\varepsilon}\frac{{\partial {u_{\Delta}^{(1\varepsilon)}}}}{{\partial t}}- \frac{\partial }{{\partial {x_i}}}\Big( {{k_{ij}^{\varepsilon}}\frac{{\partial {u_{\Delta}^{(1\varepsilon)} }}}{{\partial {x_j}}}}\Big)=F_0(\bm{x},\bm{y},t)+\varepsilon F_1(\bm{x},\bm{y},t),\;\;\text{in}\;\;\Omega\times(0,T),\\
&- \frac{\partial }{{\partial {x_i}}}\Big( {{\sigma_{ij}^{\varepsilon}}\frac{{\partial {\phi_{\Delta}^{(1\varepsilon)} }}}{{\partial {x_j}}}} \Big) = E_0(\bm{x},\bm{y},t)+\varepsilon E_1(\bm{x},\bm{y},t),\;\;\text{in}\;\;\Omega\times(0,T),\\
&{u_{\Delta}^{(1\varepsilon)}(\bm{x},t) } = -\varepsilon{M_{\alpha_1}}\frac{\partial u_{0}}{\partial x_{\alpha_1}}=\varepsilon \widehat \psi_1(\bm{x},t),\;\;\text{on}\;\;\partial {\Omega_u}\times(0,T),\\
&{\phi_{\Delta}^{(1\varepsilon)} }(\bm{x},t)=-\varepsilon N_{\alpha_1}\frac{\partial \phi_{0}}{\partial x_{\alpha_1}}=\varepsilon \widehat \chi_1(\bm{x},t),\;\;\text{on}\;\;\partial\Omega_{\phi}\times(0,T),\\
&{k_{ij}^{{\varepsilon}}\frac{\partial u_{\Delta}^{(1\varepsilon)}(\bm{x},t) }{\partial {x_j}}}{n_i} = \bar \zeta_{1i}(\bm{x},t){n_i},\;\;\text{on}\;\;\partial {\Omega_q}\times(0,T),\\
&{\sigma_{ij}^{{\varepsilon}}\frac{\partial {\phi_{\Delta}^{(1\varepsilon)} }(\bm{x},t)}{\partial {x_j}}}{n_i}=\bar{\eta}_{1i}(\bm{x},t){n_i},\;\;\text{on}\;\;\partial\Omega_{d}\times(0,T),\\
&u_{\Delta}^{(1\varepsilon)}({\bm{x}},0)=-\varepsilon{M_{\alpha_1}}\frac{\partial \widetilde u}{\partial x_{\alpha_1}}=\varepsilon \widetilde \omega_1(\bm{x}),\;\;\text{in}\;\;\Omega,
\end{aligned} \right.
\end{equation}
where the specific expressions of functions ${F}_{0}(\bm{x},\bm{y},t)$, ${F}_{1}(\bm{x},\bm{y},t)$, ${E}_{0}(\bm{x},\bm{y},t)$ and ${E}_{1}(\bm{x},\bm{y},t)$ are uncomplicated to obtain and be exhibited in Appendix B of the present study because of their lengthy expressions.

Then putting the residual functions $u^{(2\varepsilon)}_{\Delta}$ and $\phi^{(2\varepsilon)}_{\Delta}$ into (1.1), we have the following residual equations for the HOMS solutions.
\begin{equation}
\left\{
\begin{aligned}
&\rho^{\varepsilon}c^{\varepsilon}\frac{{\partial {u_{\Delta}^{(2\varepsilon)}}}}{{\partial t}}- \frac{\partial }{{\partial {x_i}}}\Big( {{k_{ij}^{\varepsilon}}\frac{{\partial {u_{\Delta}^{(2\varepsilon)}}}}{{\partial {x_j}}}}\Big)=\varepsilon F_2(\bm{x},\bm{y},t),\;\;\text{in}\;\;\Omega\times(0,T),\\
&\!- \frac{\partial }{{\partial {x_i}}}\Big( {{\sigma_{ij}^{\varepsilon}}\frac{{\partial {\phi_{\Delta}^{(2\varepsilon)} }}}{{\partial {x_j}}}} \Big) = \varepsilon E_2(\bm{x},\bm{y},t),\;\;\text{in}\;\;\Omega\times(0,T),\\
&\!{u_{\Delta}^{(2\varepsilon)}(\bm{x},t) } =-\varepsilon{M_{\alpha_1}}\frac{\partial u_{0}}{\partial x_{\alpha_1}}- \varepsilon^2\big[ Q\frac{{\partial {u_{0}}}}{{\partial t}}+{M_{\alpha_1\alpha_2}}\frac{\partial^2 u_{0}}{\partial x_{\alpha_1}\partial x_{\alpha_2}}+R_{\alpha_1}\frac{\partial u_{0}}{\partial x_{\alpha_1}}\\
&\!+{H_{\alpha_1\alpha_2}}\frac{\partial u_{0}}{\partial x_{\alpha_1}}\frac{\partial u_{0}}{\partial x_{\alpha_2}}+{G_{\alpha_1\alpha_2}}\frac{\partial \phi_{0}}{\partial x_{\alpha_1}}\frac{\partial \phi_{0}}{\partial x_{\alpha_2}}\big]=\varepsilon \widehat \psi_2(\bm{x},t),\;\;\text{on}\;\;\partial {\Omega_u}\times(0,T),\\
&\!{\phi_{\Delta}^{(2\varepsilon)} }(\bm{x},t)=-\varepsilon N_{\alpha_1}\frac{\partial \phi_{0}}{\partial x_{\alpha_1}}-\varepsilon^2\big[N_{\alpha_1\alpha_2}\frac{\partial^2 \phi_{0}}{\partial x_{\alpha_1}\partial x_{\alpha_2}}+ Z_{\alpha_1}\frac{\partial \phi_{0}}{\partial x_{\alpha_1}}\\
&\!+W_{\alpha_1\alpha_2}\frac{\partial u_{0}}{\partial x_{\alpha_1}}\frac{\partial \phi_{0}}{\partial x_{\alpha_2}}\big]=\varepsilon \widehat\chi_2(\bm{x},t),\;\;\text{on}\;\;\partial\Omega_{\phi}\times(0,T),\\
&\!{k_{ij}^{{\varepsilon}}\frac{\partial u_{\Delta}^{(2\varepsilon)}(\bm{x},t) }{\partial {x_j}}}{n_i} = \bar \zeta_{2i}(\bm{x},t){n_i},\;\;\text{on}\;\;\partial {\Omega_q}\times(0,T),\\
&\!{\sigma_{ij}^{{\varepsilon}}\frac{\partial {\phi_{\Delta}^{(2\varepsilon)} }(\bm{x},t)}{\partial {x_j}}}{n_i}=\bar \eta_{2i}(\bm{x},t){n_i},\;\;\text{on}\;\;\partial\Omega_{d}\times(0,T),\\
&\!u_{\Delta}^{(2\varepsilon)}({\bm{x}},0) =-\varepsilon{M_{\alpha_1}}\frac{\partial \widetilde u}{\partial x_{\alpha_1}}- \varepsilon^2\big[ Q\frac{{\partial {u_{0}}}}{{\partial t}}\big|_{t=0}+{M_{\alpha_1\alpha_2}}\frac{\partial^2 \widetilde u}{\partial x_{\alpha_1}\partial x_{\alpha_2}}+R_{\alpha_1}\frac{\partial \widetilde u}{\partial x_{\alpha_1}}\\
&\!+{H_{\alpha_1\alpha_2}}\frac{\partial \widetilde u}{\partial x_{\alpha_1}}\frac{\partial \widetilde u}{\partial x_{\alpha_2}}+{G_{\alpha_1\alpha_2}}\frac{\partial \phi_{0}}{\partial x_{\alpha_1}}\big|_{t=0}\frac{\partial \phi_{0}}{\partial x_{\alpha_2}}\big|_{t=0}\big]=\varepsilon \widetilde\omega_2(\bm{x}),\;\;\text{in}\;\;\Omega,
\end{aligned} \right.
\end{equation}
where the detailed expressions of functions ${F}_{2}(\bm{x},\bm{y},t)$ and ${E}_{2}(\bm{x},\bm{y},t)$ are also uncomplicated to achieve and be displayed in Appendix B of the present study owing to their lengthy expressions.

From the preceding error results in the pointwise sense, it is apparent that LOMS solutions fail to maintain local balance of heat flux and electric charge since the $\varepsilon$-independent terms ${F}_{0}(\bm{x},\bm{y},t)$ and ${E}_0(\bm{x},\bm{y},t)$ in (3.2) can not converge toward zero as the microstructural parameter $\varepsilon$ approaches zero. Benefitting from the higher-order correction terms, the HOMS solutions can guarantee the local heat flux balance of thermal equation and local electric charge balance of electric equation in the original governing equations (1.1) due to their $O(\varepsilon)$-order pointwise errors. This serves as the principal motivation for this study to establish the HOMS solutions exhibiting high-accuracy computation performance for composite structures.
\subsection{The convergence proof by global error estimation}
To acquire global error estimation of higher-order multiscale asymptotic solutions in the integral sense, certain assumptions are further presented as below.
\begin{enumerate}
\item[(B$_1$)]
Assume that $\Omega$ is a bounded, convex and smooth domain, which can be decomposed into an internal integral periodic region $\Omega_0$ and a boundary layer region $\Omega_1$. As shown in Fig.\hspace{1mm}1, $\bar{\Omega}_0=\cup_{\mathbf{z}\in I_{\varepsilon}}\varepsilon(\mathbf{z}+\bar{\Theta})$, where the index set $I_{\varepsilon}=\{\mathbf{z}=(z_1,\cdots,z_n)\in Z^n,\varepsilon(\mathbf{z}+\bar{\Theta})\subset \bar{\Omega}_0\}$. Besides, let $E_\mathbf{z}=\varepsilon(\mathbf{z}+\Theta)$ be the translational unit cell and $\partial E_\mathbf{z}$ be the boundary of $E_\mathbf{z}$.
\item[(B$_2$)]
Assume that $\displaystyle\frac{\partial \rho^{\varepsilon}({\bm{x}},u^{\epsilon})}{\partial t}$, $\displaystyle\frac{\partial c^{\varepsilon}({\bm{x}},u^{\epsilon})}{\partial t}$ and $\displaystyle\frac{\partial k^{\varepsilon}_{ij}({\bm{x}},u^{\epsilon})}{\partial t}\in L^\infty(\Omega\times(0,T))$.
\end{enumerate}
\begin{figure}[!htb]
\centering
\begin{minipage}[c]{0.48\textwidth}
  \centering
  \includegraphics[width=0.95\linewidth,totalheight=1.3in]{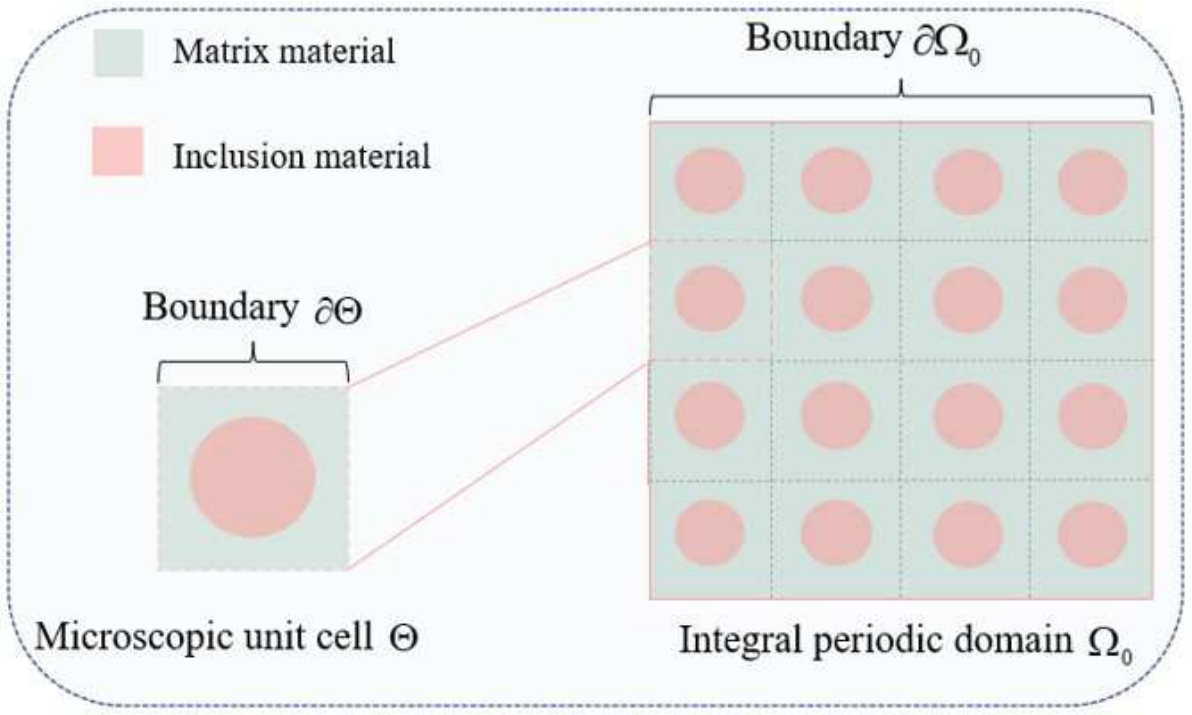} \\
  (a)
\end{minipage}
\begin{minipage}[c]{0.48\textwidth}
  \centering
  \includegraphics[width=0.95\linewidth,totalheight=1.3in]{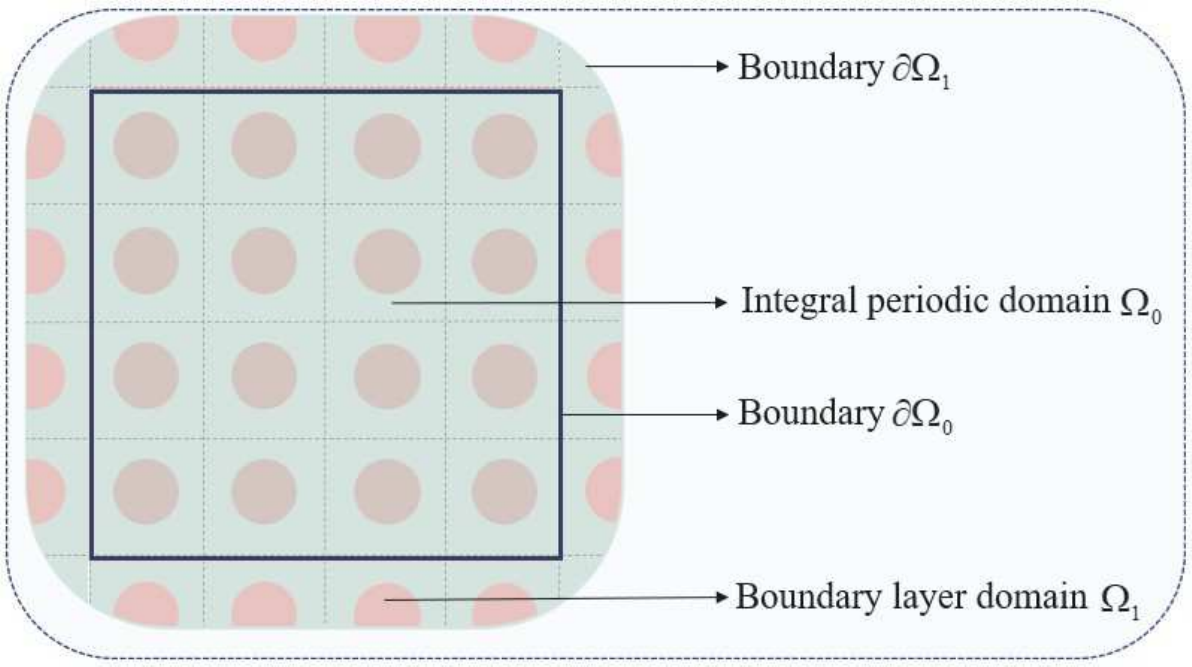} \\
  (b)
\end{minipage}
\caption{(a) Integral periodic domain $\Omega_0$; (b) non-integral periodic domain $\Omega$ with boundary layer $\Omega _1$ and ${\Omega_0}=\cup_{\mathbf{z}\in I_{\varepsilon}}\varepsilon(\mathbf{z}+\bar{\Theta})$, namely $\bar\Omega= {\bar \Omega _0} \cup {\bar \Omega _1}$.}
\end{figure}
Next, we give the ultimate result of global error estimation for the HOMS solutions of the time-dependent multiscale nonlinear equations (1.1) as the following theorem.
\begin{theorem}
Let ${u}^\varepsilon(\bm{x},t)$ and $\phi^\varepsilon(\bm{x},t)$ be the weak solutions of multiscale nonlinear equations (1.1), ${u}^{(0)}(\bm{x},t)$ and $\phi^{(0)}(\bm{x},t)$ be the weak solutions of corresponding homogenized equations (2.13), ${u^{(2\varepsilon)} }(\bm{x},t)$ and ${\phi^{(2\varepsilon)} }(\bm{x},t)$ be the HOMS solutions given by formulas (2.27) and (2.28). Under the above hypotheses (A$_1$)-(A$_3$) and (B$_1$)-(B$_2$), the following global error estimation are obtained.
\begin{equation}
{\begin{aligned}
\big\|{\phi^{\varepsilon} }-{\phi^{(2\varepsilon)} }\big\|_{L^\infty(0,T;H^1(\Omega))}\leq C(\Omega,T)\varepsilon^{1/2},
\end{aligned}}
\end{equation}
\begin{equation}
{\begin{aligned}
\big\|{u^{\varepsilon} }-{u^{(2\varepsilon)} }\big\|_{L^\infty(0,T;L^2(\Omega))}+
\big\|{u^{\varepsilon} }-{u^{(2\varepsilon)} }\big\|_{L^2(0,T;H^1(\Omega))}\leq C(\Omega,T)\varepsilon^{1/2},
\end{aligned}}
\end{equation}
\end{theorem}
where $C(\Omega,T)$ is a positive constant irrespective of $\varepsilon$, but dependent of $\Omega$ and $T$.
$\mathbf{Proof:}$
The residual equations (3.3) are employed to accomplish the global error estimation. Firstly, let us introduce the cut-off function $m_{\phi,\varepsilon}(\bm{x})\in C^{\infty}(\bar{\Omega})$ in \cite{R16} for electric potential field defined as follows
\begin{equation}
\left\{ {\begin{aligned}
&m_{\phi,\varepsilon}(\bm{x})=1,\;\mathrm{if}\;\mathrm{dist}(\bm{x},\partial\Omega_{\phi})\leq \varepsilon,\\
&m_{\phi,\varepsilon}(\bm{x})=0,\;\mathrm{if}\;\mathrm{dist}(\bm{x},\partial\Omega_{\phi})\geq 2\varepsilon,\\
&\|\bigtriangledown m_{\phi,\varepsilon}(\bm{x})\|_{L^\infty(\Omega)}\leq C\varepsilon^{-1},\;\mathrm{if}\; \varepsilon<\mathrm{dist}(\bm{x},\partial\Omega_{\phi})<2\varepsilon.
\end{aligned}} \right.
\end{equation}
Then, set new residual function for electric potential field as below
\begin{equation}
\wideparen\phi^{(2\varepsilon)}_{\Delta}(\bm{x},t)=\phi^{(2\varepsilon)}_{\Delta}(\bm{x},t)-\varepsilon m_{\phi,\varepsilon}(\bm{x})\widehat\chi_2(\bm{x},t).
\end{equation}
Then, multiplying on both sides of electric residual equation in (3.3) by $\displaystyle \wideparen\phi^{(2\varepsilon)}_{\Delta}(\bm{x},t)$ and integrating on $\Omega$, we derive the following equality.
\begin{equation}
\begin{aligned}
-\int_{\Omega} \frac{\partial}{\partial x_i}\Big( {\sigma_{ij}^{\varepsilon}({\bm{x}},{u^\varepsilon })\frac{{\partial {\phi_\Delta^{(2\varepsilon)} }}}{\partial {x_j}}}\Big)\wideparen\phi^{(2\varepsilon)}_{\Delta}(\bm{x},t)d\Omega
=\int_{\Omega}\varepsilon E_2(\bm{x},\bm{y},t)\wideparen\phi^{(2\varepsilon)}_{\Delta}(\bm{x},t)d\Omega.
\end{aligned}
\end{equation}
Subsequently, applying Green's formula on (3.8) and substituting boundary condition in it, the above identity can be rewritten as below
\begin{equation}
\begin{aligned}
&\int_{\Omega} {\sigma_{ij}^{\varepsilon}({\bm{x}},{u^\varepsilon })\frac{{\partial \wideparen\phi^{(2\varepsilon)}_{\Delta}}}{\partial {x_j}}}\frac{{\partial {\wideparen\phi_\Delta^{(2\varepsilon)} }}}{\partial {x_i}}d\Omega=\int_{\Omega}\varepsilon E_2(\bm{x},\bm{y},t)\wideparen\phi^{(2\varepsilon)}_{\Delta}d\Omega\\
&+{\int_{\partial\Omega_{d}}}\bar \eta_{2i}(\bm{x},t){n_i}\wideparen\phi^{(2\varepsilon)}_{\Delta}ds-\int_{\Omega} {\sigma_{ij}^{\varepsilon}({\bm{x}},{u^\varepsilon })\frac{{\partial \varepsilon m_{\phi,\varepsilon}(\bm{x})\widehat\chi_2(\bm{x},t)}}{\partial {x_j}}}\frac{{\partial {\wideparen\phi_\Delta^{(2\varepsilon)} }}}{\partial {x_i}}d\Omega.
\end{aligned}
\end{equation}
Furthermore, employing assumption (A$_1$) and Poincar$\mathrm{\acute{e}}$-Friedrichs inequality to the left side of (3.9), it follows that
\begin{equation}
\Big|\int_{\Omega} {\sigma_{ij}^{\varepsilon}({\bm{x}},{u^\varepsilon })\frac{{\partial \wideparen\phi^{(2\varepsilon)}_{\Delta}}}{\partial {x_j}}}\frac{{\partial {\wideparen\phi_\Delta^{(2\varepsilon)} }}}{\partial {x_i}}d\Omega\Big| \geq C\left \|\wideparen\phi^{(2\varepsilon)}_{\Delta}\right\|_{H_0^1(\Omega)}^2.
\end{equation}
After that, exploiting the Schwarz's inequality and lemma 2.2 in chapter 2 of reference \cite{R17}, the following inequality is obtained by transforming the right side of (3.9)
\begin{equation}
\begin{aligned}
&\Big|\int_{\Omega}\varepsilon E_2(\bm{x},\bm{y},t)\wideparen\phi^{(2\varepsilon)}_{\Delta}d\Omega+{\int_{\partial\Omega_{d}}}\bar \eta_{2i}(\bm{x},t){n_i}\wideparen\phi^{(2\varepsilon)}_{\Delta}ds\\
&-\int_{\Omega} {\sigma_{ij}^{\varepsilon}({\bm{x}},{u^\varepsilon })\frac{{\partial \big[\varepsilon m_{\phi,\varepsilon}(\bm{x})\widehat\chi_2(\bm{x},t)}\big]}{\partial {x_j}}}\frac{{\partial {\wideparen\phi_\Delta^{(2\varepsilon)} }}}{\partial {x_i}}d\Omega\Big|\\
&\leq \left\|\varepsilon E_2(\bm{x},\bm{y},t)\right\|_{L^2(\Omega)} \left\|\wideparen\phi^{(2\varepsilon)}_{\Delta}\right\|_{L^2(\Omega)}+C\varepsilon^{{1}/{2}}\left\|\wideparen\phi^{(2\varepsilon)}_{\Delta}\right\|_{H_0^1(\Omega)}\\
&+C\left\|\varepsilon m_{\phi,\varepsilon}(\bm{x})\widehat\chi_2(\bm{x},t)\right\|_{H^1(\Omega)} \left\|\wideparen\phi^{(2\varepsilon)}_{\Delta}\right\|_{H_0^1(\Omega)}\\
&\leq C\varepsilon\left\| E_2(\bm{x},\bm{y},t)\right\|_{L^2(\Omega)} \left\|\wideparen\phi^{(2\varepsilon)}_{\Delta}\right\|_{H_0^1(\Omega)}+C\varepsilon^{{1}/{2}}\left\|\wideparen\phi^{(2\varepsilon)}_{\Delta}\right\|_{H_0^1(\Omega)}\\
&+C\varepsilon \left\|m_{\phi,\varepsilon}(\bm{x})\widehat\chi_2(\bm{x},t)\right\|_{H^1(\Omega)} \left\|\wideparen\phi^{(2\varepsilon)}_{\Delta}\right\|_{H_0^1(\Omega)}.
\end{aligned}
\end{equation}
A combination of (3.10) and (3.11) leads to the following inequality
\begin{equation}
\left\|\wideparen\phi^{(2\varepsilon)}_{\Delta}\right\|_{H_0^1(\Omega)}\leq C\varepsilon\left\| E_2(\bm{x},\bm{y},t)\right\|_{L^2(\Omega)}+C\varepsilon^{{1}/{2}}
+C\varepsilon \left\|m_{\phi,\varepsilon}(\bm{x})\widehat\chi_2(\bm{x},t)\right\|_{H^1(\Omega)}.
\end{equation}
Then, in virtue of the triangle inequality, it is easy to obtain
\begin{equation}
\left\|\phi^{(2\varepsilon)}_{\Delta}\right\|_{H^1(\Omega)}\leq C\varepsilon\left\| E_2(\bm{x},\bm{y},t)\right\|_{L^2(\Omega)}+C\varepsilon^{{1}/{2}}
+2C\varepsilon \left\|m_{\phi,\varepsilon}(\bm{x})\widehat\chi_2(\bm{x},t)\right\|_{H^1(\Omega)}.
\end{equation}
Following the proof of the error estimate in chapter 7 of reference \cite{R16}, we can obtain
\begin{equation}
\left\|m_{\phi,\varepsilon}(\bm{x})\widehat\chi_2(\bm{x},t)\right\|_{H^1(\Omega)} =\left\|m_{\phi,\varepsilon}(\bm{x})\widehat\chi_2(\bm{x},t)\right\|_{H^1(K_{\varepsilon})}
\leq C\varepsilon^{-{1}/{2}},
\end{equation}
where $K_{\varepsilon}=\left\{{\bm{x}|\mathrm{dist}(\bm{x},\partial\Omega_{\phi})\leq 2\varepsilon}\right\}\cap\Omega$. Then, substituting (3.14) into (3.13), we have
\begin{equation}
{\begin{aligned}
\big\|{\phi^{\varepsilon} }(\bm{x},t)-{\phi^{(2\varepsilon)}(\bm{x},t) }\big\|_{H^1(\Omega)}\leq C(\Omega)\varepsilon^{1/2}.
\end{aligned}}
\end{equation}
With the arbitrariness of temporal variable $t$ in (3.15), we obtain the explicit convergence estimation (3.4) from (3.15).

Similarly, let us introduce the cut-off function $m_{u,\varepsilon}(\bm{x})\in C^{\infty}(\bar{\Omega})$ in \cite{R16} for temperature field defined as follows
\begin{equation}
\left\{ {\begin{aligned}
&m_{u,\varepsilon}(\bm{x})=1,\;\mathrm{if}\;\mathrm{dist}(\bm{x},\partial\Omega_{u})\leq \varepsilon,\\
&m_{u,\varepsilon}(\bm{x})=0,\;\mathrm{if}\;\mathrm{dist}(\bm{x},\partial\Omega_{u})\geq 2\varepsilon,\\
&\|\bigtriangledown m_{u,\varepsilon}(\bm{x})\|_{L^\infty(\Omega)}\leq C\varepsilon^{-1},\;\mathrm{if}\; \varepsilon<\mathrm{dist}(\bm{x},\partial\Omega_{u})<2\varepsilon.
\end{aligned}} \right.
\end{equation}
Then, define novel residual function for temperature field as below
\begin{equation}
\wideparen u^{(2\varepsilon)}_{\Delta}(\bm{x},t)=u^{(2\varepsilon)}_{\Delta}(\bm{x},t)-\varepsilon m_{u,\varepsilon}(\bm{x})\widehat\psi_2(\bm{x},t).
\end{equation}
Afterwards, multiplying on both sides of temperature residual equation in (3.3) by $\displaystyle \wideparen u^{(2\varepsilon)}_{\Delta}(\bm{x},t)$ and integrating on $\Omega$, the following equality achieves
\begin{equation}
\begin{aligned}
&\int_{\Omega}\rho^{\varepsilon}({\bm{x}},{u^\varepsilon })c^{\varepsilon}({\bm{x}},{u^\varepsilon })\frac{{\partial {u_\Delta^{(2\varepsilon)} }}}{{\partial t}}\wideparen u_\Delta^{(2\varepsilon)}d\Omega-\int_{\Omega} \frac{\partial}{\partial x_i}\Big( {k_{ij}^{\varepsilon}({\bm{x}},{u^\varepsilon })\frac{{\partial {u_\Delta^{(2\varepsilon)} }}}{\partial {x_j}}}\Big)\wideparen u_\Delta^{(2\varepsilon)}d\Omega \\
&=\int_{\Omega}\varepsilon F_2(\bm{x},\bm{y},t)\wideparen u_\Delta^{(2\varepsilon)}d\Omega.
\end{aligned}
\end{equation}
Furthermore, applying Green's formula on (3.18) and substituting boundary condition in it, the above identity can be rewritten as below
\begin{equation}
\begin{aligned}
&\int_{\Omega}\rho^{\varepsilon}({\bm{x}},{u^\varepsilon })c^{\varepsilon}({\bm{x}},{u^\varepsilon })\frac{{\partial {\wideparen u_\Delta^{(2\varepsilon)} }}}{{\partial t}}\wideparen u_\Delta^{(2\varepsilon)}d\Omega+\int_{\Omega} {k_{ij}^{\varepsilon}({\bm{x}},{u^\varepsilon })\frac{{\partial {\wideparen u_\Delta^{(2\varepsilon)} }}}{\partial {x_j}}}\frac{{\partial {\wideparen u_\Delta^{(2\varepsilon)} }}}{\partial {x_i}}d\Omega\\
&=\int_{\Omega}\varepsilon  F_2(\bm{x},\bm{y},t)\wideparen u_\Delta^{(2\varepsilon)}d\Omega+{\int_{\partial\Omega_{q}}}\bar \zeta_{2i}(\bm{x},t){n_i}\wideparen u^{(2\varepsilon)}_{\Delta}ds\\
&-\int_{\Omega}\rho^{\varepsilon}({\bm{x}},{u^\varepsilon })c^{\varepsilon}({\bm{x}},{u^\varepsilon })\frac{{\partial \big[{\varepsilon m_{u,\varepsilon}(\bm{x})\widehat\psi_2(\bm{x},t)}}\big]}{{\partial t}}\wideparen u_\Delta^{(2\varepsilon)}d\Omega\\
&-\int_{\Omega} {k_{ij}^{\varepsilon}({\bm{x}},{u^\varepsilon })\frac{{\partial \big[\varepsilon m_{u,\varepsilon}(\bm{x})\widehat\psi_2(\bm{x},t)}\big]}{\partial {x_j}}}\frac{{\partial {\wideparen u_\Delta^{(2\varepsilon)} }}}{\partial {x_i}}d\Omega.
\end{aligned}
\end{equation}
After that, it is easy to obtain the following equality for multi-scale thermal equation with temperature-dependent material parameters.
\begin{equation}
\begin{aligned}
&\frac{1}{2}\frac{\partial}{\partial t}\Big[\int_{\Omega}\rho^{\varepsilon}({\bm{x}},{u^\varepsilon })c^{\varepsilon}({\bm{x}},{u^\varepsilon })(\wideparen u_\Delta^{(2\varepsilon)})^2d\Omega\Big]+\int_{\Omega} {k_{ij}^{\varepsilon}({\bm{x}},{u^\varepsilon })\frac{{\partial {\wideparen u_\Delta^{(2\varepsilon)} }}}{\partial {x_j}}}\frac{{\partial {\wideparen u_\Delta^{(2\varepsilon)} }}}{\partial {x_i}}d\Omega\\
&=\frac{1}{2}\int_{\Omega}\frac{\partial \rho^{\varepsilon}({\bm{x}},{u^\varepsilon })}{\partial t}c^{\varepsilon}({\bm{x}},{u^\varepsilon })(\wideparen u_\Delta^{(2\varepsilon)})^2d\Omega+\frac{1}{2}\int_{\Omega}\rho^{\varepsilon}({\bm{x}},{u^\varepsilon })\frac{\partial c^{\varepsilon}({\bm{x}},{u^\varepsilon })}{\partial t}(\wideparen u_\Delta^{(2\varepsilon)})^2d\Omega\\
&+\int_{\Omega}\varepsilon  F_2(\bm{x},\bm{y},t)\wideparen u_\Delta^{(2\varepsilon)}d\Omega+{\int_{\partial\Omega_{q}}}\bar \zeta_{2i}(\bm{x},t){n_i}\wideparen u^{(2\varepsilon)}_{\Delta}ds\\
&-\int_{\Omega}\rho^{\varepsilon}({\bm{x}},{u^\varepsilon })c^{\varepsilon}({\bm{x}},{u^\varepsilon })\frac{{\partial \big[{\varepsilon m_{u,\varepsilon}(\bm{x})\widehat\psi_2(\bm{x},t)}}\big]}{{\partial t}}\wideparen u_\Delta^{(2\varepsilon)}d\Omega\\
&-\int_{\Omega} {k_{ij}^{\varepsilon}({\bm{x}},{u^\varepsilon })\frac{{\partial \big[\varepsilon m_{u,\varepsilon}(\bm{x})\widehat\psi_2(\bm{x},t)}\big]}{\partial {x_j}}}\frac{{\partial {\wideparen u_\Delta^{(2\varepsilon)} }}}{\partial {x_i}}d\Omega.
\end{aligned}
\end{equation}
Subsequently, we integrate both sides of (3.20) from $0$ to $t$ $(0<t\leq T)$ and substitute the initial condition of temperature residual equation in (3.3) after integration. Then it follows that
\begin{equation}
\begin{aligned}
&\int_{\Omega}\rho^{\varepsilon}c^{\varepsilon}(\wideparen u_\Delta^{(2\varepsilon)})^2d\Omega+\int_0^t\int_{\Omega}2 {k_{ij}^{\varepsilon}\frac{{\partial {\wideparen u_\Delta^{(2\varepsilon)}}}}{\partial {x_j}}}\frac{{\partial {\wideparen u_\Delta^{(2\varepsilon)}}}}{\partial {x_i}}d\Omega d\tau={\int_{\Omega}}\rho^{\varepsilon}c^{\varepsilon}\big(\varepsilon\widetilde\omega_2(\bm{x})\big)^2d\Omega\\
&+\int_0^t\int_{\Omega}\frac{\partial \rho^{\varepsilon}}{\partial \tau}c^{\varepsilon}(\wideparen u_\Delta^{(2\varepsilon)})^2d\Omega d\tau+\int_0^t\int_{\Omega}\rho^{\varepsilon}\frac{\partial c^{\varepsilon}}{\partial \tau}(\wideparen u_\Delta^{(2\varepsilon)})^2d\Omega d\tau\\
&+\int_0^t\int_{\Omega}2\varepsilon F_2(\bm{x},\bm{y},\tau)\wideparen u_\Delta^{(2\varepsilon)}d\Omega d\tau+\int_0^t{\int_{\partial\Omega_{q}}}2\bar \zeta_{2i}(\bm{x},\tau){n_i}\wideparen u^{(2\varepsilon)}_{\Delta}dsd\tau\\
&-\int_0^t\int_{\Omega}2\rho^{\varepsilon}c^{\varepsilon}\frac{{\partial \big[{\varepsilon m_{u,\varepsilon}(\bm{x})\widehat\psi_2(\bm{x},\tau)}}\big]}{{\partial \tau}}\wideparen u_\Delta^{(2\varepsilon)}d\Omega d\tau\\
&-\int_0^t\int_{\Omega} 2{k_{ij}^{\varepsilon}\frac{{\partial \big[\varepsilon m_{u,\varepsilon}(\bm{x})\widehat\psi_2(\bm{x},\tau)}\big]}{\partial {x_j}}}\frac{{\partial {\wideparen u_\Delta^{(2\varepsilon)} }}}{\partial {x_i}}d\Omega d\tau.
\end{aligned}
\end{equation}
Owing to assumptions (A$_1$) and (A$_2$), and employing Poincar$\rm{\acute{e}}$-Friedrichs inequality, the following inequality can be easily obtained from left side of equality (3.21).
\begin{equation}
\begin{aligned}
&\int_{\Omega}\rho^{\varepsilon}c^{\varepsilon}(\wideparen u_\Delta^{(2\varepsilon)})^2d\Omega+\int_0^t\int_{\Omega}2 {k_{ij}^{\varepsilon}\frac{{\partial {\wideparen u_\Delta^{(2\varepsilon)}}}}{\partial {x_j}}}\frac{{\partial {\wideparen u_\Delta^{(2\varepsilon)}}}}{\partial {x_i}}d\Omega d\tau\\
&\geq\rho^*c^*\left \|\wideparen u_{\Delta}^{(2\varepsilon)}(\bm{x},t)\right\|_{L^2(\Omega)}^2
+C_1\int_0^t\left \|\wideparen u_{\Delta}^{(2\varepsilon)}(\bm{x},\tau)\right\|_{H^1_0(\Omega)}^2d\tau.
\end{aligned}
\end{equation}
After that, employing Schwarz's inequality and Young's inequality $\displaystyle ab\leq\frac{1}{2}(\lambda a^2+\frac{1}{\lambda}b^2),\forall\lambda\in R^{+}$ with parameter $\lambda$, we derive the
following inequality by transforming the right side of equality (3.21).
\begin{equation}
\begin{aligned}
&{\int_{\Omega}}\rho^{\varepsilon}c^{\varepsilon}\big(\varepsilon\widetilde\omega_2(\bm{x})\big)^2d\Omega+\int_0^t\int_{\Omega}2\varepsilon F_2(\bm{x},\bm{y},\tau)\wideparen u_\Delta^{(2\varepsilon)}d\Omega d\tau\\
&+\int_0^t\int_{\Omega}\frac{\partial \rho^{\varepsilon}}{\partial \tau}c^{\varepsilon}(\wideparen u_\Delta^{(2\varepsilon)})^2d\Omega d\tau+\int_0^t\int_{\Omega}\rho^{\varepsilon}\frac{\partial c^{\varepsilon}}{\partial \tau}(\wideparen u_\Delta^{(2\varepsilon)})^2d\Omega d\tau\\
&+\int_0^t{\int_{\partial\Omega_{q}}}2\bar \zeta_{2i}(\bm{x},\tau){n_i}\wideparen u^{(2\varepsilon)}_{\Delta}dsd\tau\\
&-\int_0^t\int_{\Omega}2\rho^{\varepsilon}c^{\varepsilon}\frac{{\partial \big[{\varepsilon m_{u,\varepsilon}(\bm{x})\widehat\psi_2(\bm{x},\tau)}}\big]}{{\partial \tau}}\wideparen u_\Delta^{(2\varepsilon)}d\Omega d\tau\\
&-\int_0^t\int_{\Omega} 2{k_{ij}^{\varepsilon}\frac{{\partial \big[\varepsilon m_{u,\varepsilon}(\bm{x})\widehat\psi_2(\bm{x},\tau)}\big]}{\partial {x_j}}}\frac{{\partial {\wideparen u_\Delta^{(2\varepsilon)} }}}{\partial {x_i}}d\Omega d\tau\\
&\leq C_2\varepsilon^2+2\int_0^t\int_{\Omega}\frac{(\varepsilon F_2(\bm{x},\bm{y},\tau))^2+(\wideparen u_{\Delta}^{(2\varepsilon)}(\bm{x},\tau))^2}{2}d\Omega d\tau\\
&+C_2\int_0^t\left \|\wideparen u_{\Delta}^{(2\varepsilon)}(\bm{x},\tau)\right\|_{L^2(\Omega)}^2d\tau+C_2\int_0^t\left \| \wideparen u_{\Delta}^{(2\varepsilon)}(\bm{x},\tau)\right\|_{L^2(\Omega)}^2d\tau\\
&+\int_0^tC_2\varepsilon^{{1}/{2}}\left\|\wideparen u^{(2\varepsilon)}_{\Delta}\right\|_{H_0^1(\Omega)}d\tau+\int_0^tC_2\varepsilon \left\|\frac{{\partial {m_{u,\varepsilon}(\bm{x})\widehat\psi_2}}}{{\partial \tau}}\right\|_{L^2(\Omega)} \left\|\wideparen u^{(2\varepsilon)}_{\Delta}\right\|_{L^2(\Omega)}d\tau\\
&+\int_0^tC_2\varepsilon \left\|m_{u,\varepsilon}(\bm{x})\widehat\psi_2(\bm{x},\tau)\right\|_{H^1(\Omega)} \left\|\wideparen u^{(2\varepsilon)}_{\Delta}\right\|_{H_0^1(\Omega)}d\tau\\
&\leq C_2\varepsilon^2+C_2\int_0^t\left \|\wideparen u_{\Delta}^{(2\varepsilon)}(\bm{x},\tau)\right\|_{L^2(\Omega)}^2d\tau+C_2\int_0^t\varepsilon^{{1}/{2}}\left\|\wideparen u^{(2\varepsilon)}_{\Delta}\right\|_{H_0^1(\Omega)}d\tau\\
&+C_2\int_0^t\varepsilon\left\|\wideparen u_{\Delta}^{(2\varepsilon)}(\bm{x},\tau)\right\|_{L^2(\Omega)}d\tau\\
&\leq C_2\varepsilon^2+C_2\int_0^t\left \|\wideparen u_{\Delta}^{(2\varepsilon)}(\bm{x},\tau)\right\|_{L^2(\Omega)}^2d\tau+C_2\int_0^t\Big[\frac{1}{2\lambda}\varepsilon+\frac{\lambda}{2}\left\|\wideparen u^{(2\varepsilon)}_{\Delta}\right\|_{H_0^1(\Omega)}^2\Big]d\tau\\
&\leq C_2\varepsilon+C_2\int_0^t\left \|\wideparen u_{\Delta}^{(2\varepsilon)}(\bm{x},\tau)\right\|_{L^2(\Omega)}^2d\tau+\frac{C_2\lambda}{2}\int_0^t\left\|\wideparen u^{(2\varepsilon)}_{\Delta}\right\|_{H_0^1(\Omega)}^2d\tau.
\end{aligned}
\end{equation}
Afterwards, choosing a sufficiently small $\lambda$ fulfilled $C_1-C_2\lambda/2>0$, and combining (3.22) and (3.23), it is apparent that
\begin{equation}
\begin{aligned}
&\rho^*c^*\left \|\wideparen u_{\Delta}^{(2\varepsilon)}(\bm{x},t)\right\|_{L^2(\Omega)}^2
+(C_1-\frac{C_2\lambda}{2})\int_0^t\left \|\wideparen u_{\Delta}^{(2\varepsilon)}(\bm{x},\tau)\right\|_{H^1_0(\Omega)}^2d\tau\\
&\leq C_2\varepsilon+C_2\int_0^t\left \|\wideparen u_{\Delta}^{(2\varepsilon)}(\bm{x},\tau)\right\|_{L^2(\Omega)}^2d\tau+C_2\int_0^t\int_0^\tau\left \| \wideparen u_{\Delta}^{(2\varepsilon)}(\bm{x},s)\right\|_{H^1_0(\Omega)}^2dsd\tau.
\end{aligned}
\end{equation}
Without loss of generality, we define $C=C_2/\mathrm{min}(\rho^*c^*,C_1-C_2\lambda/2)$ and $\displaystyle\Re(t)=\big \| \wideparen u_{\Delta}^{(2\varepsilon)}\big\|_{L^2(\Omega)}^2
+\mathlarger{\int}_0^t\big \|\wideparen u_{\Delta}^{(2\varepsilon)}\big\|_{H^1_0(\Omega)}^2d\tau$, then we have $\Re(t)\leq C(\Omega)\varepsilon+C(\Omega)\mathlarger{\int}_0^t\Re(\tau)d\tau$ from (3.24). It follows from Gronwall's inequality in chapter 12 of reference \cite{R16} that $\Re(t)\leq C(\Omega,T)\varepsilon$. Consequently, there holds the following inequality
\begin{equation}
\begin{aligned}
\big \|\wideparen u_{\Delta}^{(2\varepsilon)}\big\|_{L^2(\Omega)}^2
+\mathlarger{\int}_0^t\big \|\wideparen u_{\Delta}^{(2\varepsilon)}\big\|_{H^1_0(\Omega)}^2d\tau\leq C(\Omega,T)\varepsilon.
\end{aligned}
\end{equation}
Furthermore, in virtue of the triangle inequality and by adopting a similar error estimate as (3.14), it is easy to obtain
\begin{equation}
\begin{aligned}
&\big \|u_{\Delta}^{(2\varepsilon)}\big\|_{L^2(\Omega)}^2
+\mathlarger{\int}_0^t\big \|u_{\Delta}^{(2\varepsilon)}\big\|_{H^1(\Omega)}^2d\tau\\
&\leq C(\Omega,T)\varepsilon+\big \|\varepsilon m_{u,\varepsilon}(\bm{x})\widehat\psi_2(\bm{x},t)\big\|_{L^2(\Omega)}^2
+\mathlarger{\int}_0^t\big \|\varepsilon m_{u,\varepsilon}(\bm{x})\widehat\psi_2(\bm{x},t)\big\|_{H^1(\Omega)}^2d\tau\\
&\leq C(\Omega,T)\varepsilon+(\varepsilon C\varepsilon^{-1/2})^2+(\varepsilon C\varepsilon^{-1/2})^2\leq C(\Omega,T)\varepsilon.
\end{aligned}
\end{equation}
Then using the arithmetic and geometric means inequality $\displaystyle{(a+b)}/{2}\leq\sqrt{{(a^2+b^2)}/{2}}$ to the left side of the inequality (3.26) and squaring root on both sides of the inequality (3.26), the following inequality is obtained
\begin{equation}
\begin{aligned}
\big \|u_{\Delta}^{(2\varepsilon)}\big\|_{L^2(\Omega)}
+\big\|u_{\Delta}^{(2\varepsilon)}\big\|_{L^2(0,t;H^1(\Omega))}\leq C(\Omega,T)\varepsilon^{{1}/{2}}.
\end{aligned}
\end{equation}
To account for the arbitrariness of time variable $t$, we obtain the explicit convergence estimation (3.5) from (3.27).

Since the HOMS solutions of multiscale nonlinear coupling problem (1.1) do not satisfy the boundary conditions on $\Omega$ in a general domain, obtaining the optimal convergence order is impeded by the resulting boundary error. To gain the optimal error estimation, certain hypotheses are further presented as follows.
\begin{enumerate}
\item[(I)]
Suppose that $\Omega$ is a bounded and integral periodic region, i.e. $\bar{\Omega}=\cup_{\mathbf{z}\in I_{\varepsilon}}\varepsilon(\mathbf{z}+\bar{\Theta})$, where the index set $I_{\varepsilon}=\{\mathbf{z}=(z_1,\cdots,z_n)\in Z^n,\varepsilon(\mathbf{z}+\bar{\Theta})\subset \bar{\Omega}\}$. Besides, let $E_\mathbf{z}=\varepsilon(\mathbf{z}+\Theta)$ be the translational unit cell and $\partial E_\mathbf{z}$ be its boundary.
\item[(II)]
Apply the homogeneous Dirichlet boundary condition to replace the periodic boundary condition for whole auxiliary cell functions \cite{R5,R24,R36}.
\item[(III)]
Let $k_{ij}(\bm{y},{u^\varepsilon})=k_{ii}(\bm{y},{u^\varepsilon})\delta_{ij}$ and $\sigma_{ij}(\bm{y},{u^\varepsilon})=\sigma_{ii}(\bm{y},{u^\varepsilon})\delta_{ij}$, and $\delta_{ij}$ is a Kronecker symbol. Moreover, let $\Delta_1,\cdots,\Delta_n$ be the middle hyperplanes of PUC $\Theta=(0,1)^n$. Then suppose that material parameters $\rho(\bm{y},{u^\varepsilon})$, $c(\bm{y},{u^\varepsilon})$, $k_{ii}(\bm{y},{u^\varepsilon})$ and $\sigma_{ii}(\bm{y},{u^\varepsilon})$ are symmetric with respect to $\Delta_1,\cdots,\Delta_n$ for stationary ${u^\varepsilon }\in [u_{min},u_{max}+C_*]$.
\item[(IV)]
The multiscale nonlinear problem (1.1) is imposed with pure Dirichlet boundary conditions.
\end{enumerate}
On the basis of the above assumptions, one important lemma is obtained, which will be utilized to conduct optimal error estimation of the HOMS solutions.
\begin{lemma}
Defining two derivative operators $\displaystyle\sigma_{u \Theta}(\chi)=n_i k_{ij}(\bm{y},{u^\varepsilon })\frac{\partial \chi}{\partial y_j}$ and $\displaystyle\sigma_{\phi\Theta}(\chi)=n_i \sigma_{ij}(\bm{y},{u^\varepsilon })\frac{\partial \chi}{\partial y_j}$, then on the basis of foregoing assumptions (A$_1$)-(A$_2$) and (II)-(III), the normal derivatives $\sigma_{u\Theta}(M_{\alpha_1})$, $\sigma_{u\Theta}(Q)$, $\sigma_{u\Theta}(M_{\alpha_1\alpha_2})$, $\sigma_{u\Theta}(R_{\alpha_1})$, $\sigma_{u\Theta}(H_{\alpha_1\alpha_2})$, $\sigma_{u\Theta}(G_{\alpha_1\alpha_2})$, $\sigma_{\phi\Theta}(N_{\alpha_1})$, $\sigma_{\phi\Theta}(N_{\alpha_1\alpha_2})$, $\sigma_{\phi\Theta}(Z_{\alpha_1})$ and $\sigma_{\phi\Theta}(W_{\alpha_1\alpha_2})$ are continuous on the boundary of PUC $\Theta$ via the identical approach in Refs. \cite{R5,R24,R36}.
\end{lemma}
\begin{corollary}
Assume that $\Omega$ is the integral periodic region. Let ${u}^\varepsilon(\bm{x},t)$ and $\phi^\varepsilon(\bm{x},t)$ be the weak solutions of multiscale nonlinear equations (1.1), ${u}^{(0)}(\bm{x},t)$ and $\phi^{(0)}(\bm{x},t)$ be the weak solutions of corresponding homogenized equations (2.13), ${u^{(2\varepsilon)} }(\bm{x},t)$ and ${\phi^{(2\varepsilon)} }(\bm{x},t)$ be the HOMS solutions given by formulas (2.27) and (2.28). Under the above hypotheses (A$_1$)-(A$_3$), (B$_1$)-(B$_2$) and (I)-(IV), the following global error estimation are obtained.
\begin{equation}
{\begin{aligned}
\big\|{\phi^{\varepsilon} }-{\phi^{(2\varepsilon)} }\big\|_{L^\infty(0,T;H^1(\Omega))}\leq C(\Omega,T)\varepsilon,
\end{aligned}}
\end{equation}
\begin{equation}
{\begin{aligned}
\big\|{u^{\varepsilon} }-{u^{(2\varepsilon)} }\big\|_{L^\infty(0,T;L^2(\Omega))}+
\big\|{u^{\varepsilon} }-{u^{(2\varepsilon)} }\big\|_{L^2(0,T;H^1_0(\Omega))}\leq C(\Omega,T)\varepsilon,
\end{aligned}}
\end{equation}
\end{corollary}
where $C(\Omega,T)$ is a positive constant irrespective of $\varepsilon$, but dependent of $\Omega$ and $T$.
$\mathbf{Proof:}$ Recalling the above proof again, based on assumptions (I)-(IV), the error order $\varepsilon^{{1}/{2}}$ generating from boundary $\partial \Omega$ will not appear in the proof. However, it should be noted that the new auxiliary cell functions with homogeneous Dirichlet boundary condition don't have enough regularity on the outer boundary of unit cell $\Theta$ in general case.

At this time, $\sigma_{\phi\Theta}(\phi_{\Delta}^{(2\varepsilon)})$ and $\sigma_{u\Theta}(u_{\Delta }^{(2\varepsilon)})$ arise from employing the Green's formula on interface $\partial E_\mathbf{z}$. Recalling lemma 3.2, we shall hereby get the following results for the integral terms on the boundary $\partial E_\mathbf{z}$.
\begin{equation}
\left\{ %\begin{array}{l}
\begin{aligned}
&\sum\limits _{\mathbf{z}\in I_\varepsilon}\int_{\partial E_\mathbf{z}}\sigma_{\phi\Theta}(\phi_{\Delta}^{(2\varepsilon)})\phi_{\Delta}^{(2\varepsilon)}d\Gamma_{\bm{y}}=\sum\limits _{\mathbf{z}\in I_\varepsilon}\int_{\partial E_\mathbf{z}}\sigma_{\phi\Theta}(\phi^{\varepsilon}-\phi^{(2\varepsilon)})\phi_{\Delta}^{(2\varepsilon)}d\Gamma_{\bm{y}}
\\
&=-\sum\limits _{\mathbf{z}\in I_\varepsilon}\int_{\partial E_\mathbf{z}}\sigma_{\phi\Theta}(\phi^{(2\varepsilon)})\phi_{\Delta}^{(2\varepsilon)}d\Gamma_{\bm{y}}=0,\\
&\sum\limits _{\mathbf{z}\in I_\varepsilon}\int_{\partial E_\mathbf{z}}\sigma_{u\Theta}(u_{\Delta}^{(2\varepsilon)})u_{\Delta}^{(2\varepsilon)}d\Gamma_{\bm{y}}=\sum\limits _{\mathbf{z}\in I_\varepsilon}\int_{\partial E_\mathbf{z}}\sigma_{u\Theta}(u^{\varepsilon}-u^{(2\varepsilon)})u_{\Delta}^{(2\varepsilon)}d\Gamma_{\bm{y}}
\\
&=-\sum\limits _{\mathbf{z}\in I_\varepsilon}\int_{\partial E_\mathbf{z}}\sigma_{u\Theta}(u^{(2\varepsilon)})u_{\Delta}^{(2\varepsilon)}d\Gamma_{\bm{y}}=0.
\end{aligned} \right.
\end{equation}
Finally, following along the lines of the proof of theorem 3.1, we obtain the proof of corollary 3.3.
\section{Two-stage numerical algorithm}
The proposed multiscale computational framework consists of microscopic cell models, macroscopic homogenized model and higher-order multiscale solutions, which comprise a closed solving system. Noting that all microscopic cell functions defined by (2.11)-(2.12) and (2.17)-(2.24) are dependent on macroscopic temperature $u_0$, the continuous property of microscopic cell functions thereby can be proved in Appendix C by utilizing the similar idea of \cite{R37}. According to this continuous property, microscopic cell functions only need to be evaluated corresponding to few representative macroscopic temperatures rather than all appeared temperature points, and then exploit the interpolation technique to solve the auxiliary cell functions involving in simulation process \cite{R33,R34,R35}. In the following, we present the following two-stage numerical algorithm comprising of off-line and on-line stages for efficiently simulating the time-dependent nonlinear thermo-electric coupling problem (1.1) of composite structures, as elaborated in Fig.\hspace{1mm}2.
\begin{figure}[!htb]
\centering
\begin{minipage}[c]{0.9\textwidth}
  \centering
  \includegraphics[width=1.0\linewidth,totalheight=2.3in]{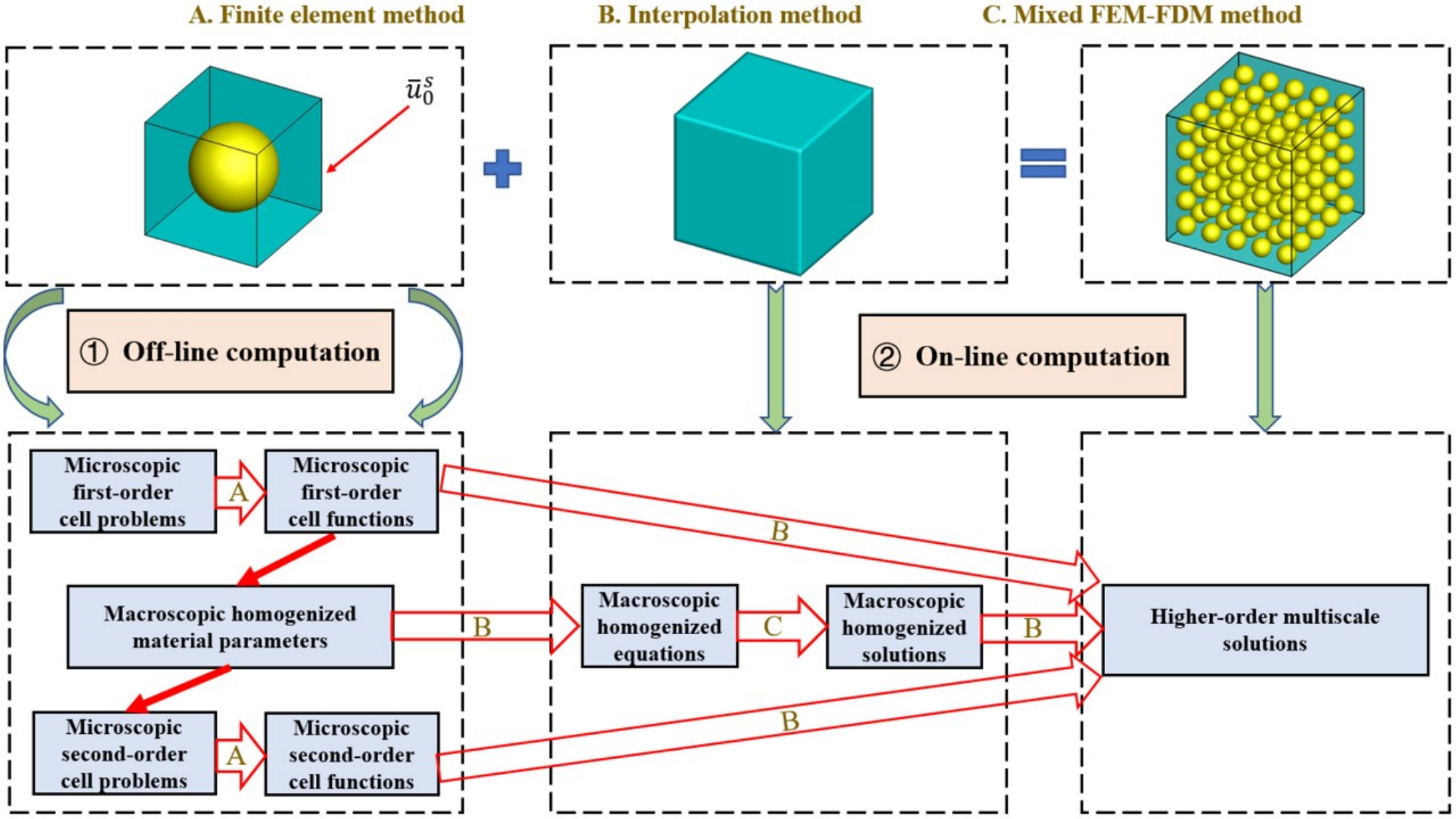}
\end{minipage}
\caption{The schematic diagram of two-stage multiscale numerical algorithm.}
\end{figure}
\subsection{Off-line stage: computation for microscopic cell problems}
\begin{enumerate}
\item[(1)]
Determine the geometric configuration of PUC $\Theta=(0,1)^n$ and create tetrahedra finite element mesh family $T_{h_1}=\{K\}$ for PUC $\Theta$, where $h_1=$max$_K\{h_K\}$. Whereupon denote the linear conforming finite element space $S_{h_1}(\Theta)=\{\nu\in C^0(\bar{\Theta}):\nu\mid_{\partial \Theta}=0,\nu\mid_{K}\in P_1(K)\}\subset H^1(\Theta)$ for auxiliary cell problems.
\item[(2)]
Define computational temperature range $[u_0^{min},u_0^{max}]$ and choose a certain number of representative macroscopic temperature $\bar u_0^{s}$ in concerned temperature range. Next, employ FEM to solve the first-order cell functions defined by (2.11)-(2.12) on $S_{h_1}(\Theta)$ corresponding to distinct representative macroscopic temperature $\bar u_0^{s}$. Note that classical periodic boundary condition of auxiliary cell problems is replaced by homogeneous Dirichlet boundary condition for practical numerical implementation \cite{R5,R24,R36}. The specific finite element scheme for first-order unit cell problem (2.11) is established as follows
\begin{equation}
\begin{aligned}
&\int_{\Theta}{ k_{ij}^{(0)}(\bm{y},\bar u_0^{s}){\frac{\partial M_{\alpha_1}(\bm{y},\bar u_0^{s})}{\partial y_j}}}\frac{\partial \upsilon^{h_1}}{\partial y_i}d\Theta\\
&=-\int_{\Theta}k_{i{\alpha_1}}^{(0)}(\bm{y},\bar u_0^{s})\frac{\partial \upsilon^{h_1}}{\partial y_i}d\Theta,\;\forall\upsilon^{h_1}\in S_{h_1}(\Theta).
\end{aligned}
\end{equation}
\item[(3)]
The macroscopic material parameters $\widehat{S}(u_0)$, $\widehat{k}_{ij}(u_0)$ and $\widehat\sigma_{ij}(u_0)$ are evaluated by formula (2.14) associated to distinct macroscopic temperature $\bar u_0^{s}$.
\item[(4)]
Utilizing the same mesh as first-order cell functions, second-order auxiliary cell functions defined by (2.17)-(2.24), which correspond to distinct representative macroscopic temperature $\bar u_0^{s}$ at macroscale, are evaluated on $S_{h_1}(\Theta)$ by employing FEM respectively.
\end{enumerate}
\subsection{On-line stage: computation for macroscopic homogenized problem}
\begin{enumerate}
\item[(1)]
Let $T_{h_0}=\{e\}$ be a tetrahedra finite element mesh family of the macroscopic region $\Omega$, where $h_0=$max$_e\{h_e\}$. Then define the linear conforming finite element spaces $S_{h_0}(\Omega)=\{\nu\in C^0(\bar{\Omega}):\nu\mid_{\partial\Omega_{u}}=0,\nu\mid_{e}\in P_1(e)\}\subset H^1(\Omega)$ and $S_{h_0}^{*}(\Omega)=\{\nu\in C^0(\bar{\Omega}):\nu\mid_{\partial\Omega_{\phi}}=0,\nu\mid_{e}\in P_1(e)\}\subset H^1(\Omega)$ for the macroscopic homogenized equations (2.13), the homogenized material parameters can be calculated by interpolation approach on each node $\bm{x}$ of $S_{h_0}(\Omega)$ and $S_{h_0}^{*}(\Omega)$.
\item[(2)]
Solve the macroscopic homogenized equations (2.13) without oscillatory coefficients by mixed FDM-FEM proposed in reference \cite{R9} on a coarse mesh and with a larger time step on the computational domain $\Omega\times(0,T)$, which means FEM is employed in spatial discretization and FDM is used to discretize time-domain. Using the equidistant time step $\displaystyle\Delta t={T}/{N}$ to discretize time-domain $(0,T)$ as $0=t_0<t_1<\cdots<t_N=T$ and $t_n=n\Delta t(n=0,\cdots,N)$, then we define $u_0^{n}=u_0(\bm{x},t_n)$ and $\phi_0^{n}=\phi_0(\bm{x},t_n)$. Moreover, define $\widehat u_{0}^{n+1/2}=(3u_{0}^{n}-u_{0}^{n-1})/2$ and $\bar u_{0}^{n+1/2}=(u_{0}^{n}+u_{0}^{n+1})/2$ for $n=0,\cdots,N-1$. Then, the computational scheme is introduced in detail. Firstly, the following computational scheme is employed to precompute $\phi_0^0$ and $\widehat u_{0}^{1/2}$.
\begin{equation}
\left\{\begin{aligned}
&\int_{\Omega}{{\widehat \sigma _{ij}}(u_{0}^0)\frac{{\partial \phi_0^0}}{{\partial {x_j}}}} \frac{\partial \varphi}{{\partial {x_i}}}d\Omega=\int_{\Omega}f_\phi(\bm{x},t_0)\varphi d\Omega\\
&+{\int_{\partial\Omega_{q}}}\bar{d}(\bm{x},t_0)\varphi ds,\;\;\forall\varphi\in S_{h_0}^{*}(\Omega),\\
&\phi_0^0=\widehat\phi(\bm{x},t_0),\;\;\text{on}\;\;\partial\Omega_{\phi}.
%&{\widehat \sigma_{ij}({\bm{x}},u_{0}^0)\frac{\partial \phi_0^0}{\partial {x_j}}}{n_i}=\bar{d}(\bm{x},t_0),\;\;\text{on}\;\;\partial\Omega_{d}.
\end{aligned}\right.
\end{equation}
\begin{equation}
\left\{ \begin{aligned}
&\int_{\Omega}\widehat S(u_{0}^0)\frac{{\widehat u_{0}^{1/2}-u_{0}^{0}}}{{\Delta t/2}}\upsilon d\Omega+ \int_{\Omega}{{\widehat k_{ij}}(u_{0}^0)\frac{{\partial {\widehat u_{0}^{1/2}}}}{{\partial {x_j}}}}\frac{\partial \upsilon}{{\partial {x_i}}}d\Omega\\
&=\int_{\Omega} {\widehat \sigma _{ij}^*}(u_{0}^0)\frac{\partial \phi_{0}^0}{\partial {x_i}}\frac{\partial \phi_{0}^0}{\partial {x_j}}\upsilon d\Omega+\int_{\Omega} f_u(\bm{x},t_{1/2})\upsilon d\Omega\\
&+{\int_{\partial\Omega_{u}}}\bar{q}(\bm{x},t_{1/2})\upsilon ds,\;\;\forall\upsilon\in S_{h_0}(\Omega),\\
&\widehat u_0^{1/2} = \widehat u(\bm{x},t_{1/2}),\;\;\text{on}\;\;\partial {\Omega_u}.
%&{{\widehat k_{ij}}(u_{0}^0)\frac{\partial \widehat u_0^{1/2}}{\partial {x_j}}}{n_i} = \bar q(\bm{x},t_{1/2}),\;\;\text{on}\;\;\partial {\Omega_q}.
\end{aligned} \right.
\end{equation}
Next, the following computational scheme is employed to compute $\phi_0^{n+1/2}$ and $u_{0}^{n+1}$ for $n=0,\cdots,N-1$.
\begin{equation}
\left\{\begin{aligned}
&\int_{\Omega}{{\widehat \sigma _{ij}}(\widehat u_{0}^{n+1/2})\frac{{\partial \phi_0^{n+1/2}}}{{\partial {x_j}}}}\frac{\partial \varphi}{{\partial {x_i}}}d\Omega=\int_{\Omega}f_\phi(\bm{x},t_{n+1/2})\varphi d\Omega\\
&+{\int_{\partial\Omega_{q}}}\bar{d}(\bm{x},t_{n+1/2})\varphi ds,\;\;\forall\varphi\in S_{h_0}^{*}(\Omega),\\
&\phi_0^{n+1/2}=\widehat\phi(\bm{x},t_{n+1/2}),\;\;\text{on}\;\;\partial\Omega_{\phi}.
%&{\widehat \sigma_{ij}({\bm{x}},\widehat u_{0}^{n+1/2})\frac{\partial \phi_0^{n+1/2}}{\partial {x_j}}}{n_i}=\bar{d}(\bm{x},t_{n+1/2}),\;\;\text{on}\;\;\partial\Omega_{d}.
\end{aligned}\right.
\end{equation}
\begin{equation}
\left\{ \begin{aligned}
&\int_{\Omega}\widehat S(\widehat u_{0}^{n+1/2})\frac{{u_{0}^{n+1}-u_{0}^{n}}}{{\Delta t}}\upsilon d\Omega+\int_{\Omega}{{\widehat k_{ij}}(\widehat u_{0}^{n+1/2})\frac{{\partial {\bar u_{0}^{n+1/2}}}}{{\partial {x_j}}}}\frac{\partial \upsilon}{{\partial {x_i}}}d\Omega\\
&= \int_{\Omega}{\widehat \sigma _{ij}^*}(\widehat u_{0}^{n+1/2})\frac{\partial \phi_{0}^{n+1/2}}{\partial {x_i}}\frac{\partial \phi_{0}^{n+1/2}}{\partial {x_j}}\upsilon d\Omega+\int_{\Omega}f_u(\bm{x},t_{n+1/2})\upsilon d\Omega\\
&+{\int_{\partial\Omega_{u}}}\bar{q}(\bm{x},t_{n+1})\upsilon ds,\;\;\forall\upsilon\in S_{h_0}(\Omega),\\
&u_0^{n+1} = \widehat u(\bm{x},t_{n+1}),\;\;\text{on}\;\;\partial {\Omega_u}.
%&{{\widehat k_{ij}}(\widehat u_{0}^{n+1/2})\frac{\partial u_0^{n+1}}{\partial {x_j}}}{n_i} = \bar q(\bm{x},t_{n+1}),\;\;\text{on}\;\;\partial {\Omega_q}.
\end{aligned} \right.
\end{equation}
\item[(3)]
The preceding fully discrete scheme results in two sub-problems. As a result, we can solve the macroscopic electric potential and temperature fields at each temporal step via two sub-problems by turn.
\end{enumerate}
\subsection{On-line stage: computation for higher-order multiscale solutions}
\begin{enumerate}
\item[(1)]
For arbitrary point $(\bm{x},t)\in \Omega\times(0,T)$, we employ the interpolation technique to solve the corresponding values of first-order cell functions, second-order cell functions and homogenized solutions.
\item[(2)]
The average technique on relative elements \cite{R36} is employed to evaluate the spatial derivatives, and the difference scheme is utilized to evaluate the temporal derivative at every time steps involving in formulas (2.27) and (2.28).
\item[(3)]
Ultimately, the temperature field ${u}^{(2\varepsilon)}(\bm{x},t)$ and the electric potential field $\phi^{(2\varepsilon)}(\bm{x},t)$ are computed by the formulas (2.27) and (2.28) separately. Besides, we can further utilize higher-order interpolation and post-processing techniques to gain high-accuracy HOMS solutions \cite{R36,R40}.
\end{enumerate}
\begin{rmk}
In this study, the equidistant macroscopic temperature values are employed as representative macroscopic temperatures.
\end{rmk}
\section{Error estimation for two-stage multiscale numerical algorithm}
\begin{lemma}
Define $M_{\alpha_1}^{h_1}$, $Q^{h_1}$, $M_{\alpha_1\alpha_2}^{h_1}$, $R_{\alpha_1}^{h_1}$, $H_{\alpha_1\alpha_2}^{h_1}$, $G_{\alpha_1\alpha_2}^{h_1}$, $N_{\alpha_1}^{h_1}$, $N_{\alpha_1\alpha_2}^{h_1}$, $Z_{\alpha_1}^{h_1}$ and $W_{\alpha_1\alpha_2}^{h_1}$ are the corresponding finite element solutions for microscopic cell functions, respectively. If all microscopic cell functions belong to $H^2(\Theta)$ for any fixed $u_0$, then the following inequalities hold
\begin{equation}
\begin{aligned}
&\left\|M_{\alpha_1}^{h_1}(\bm{y},u_0)-M_{\alpha_1}(\bm{y},u_0)\right\|_{H^m(\Theta)}\leq Ch_1^{2-m}\left\| M_{\alpha_1}(\bm{y},u_0)\right\|_{H^2(\Theta)},\\
&\left\|Q^{h_1}(\bm{y},u_0)-Q(\bm{y},u_0)\right\|_{H^m(\Theta)}\leq Ch_1^{2-m}\left\| Q(\bm{y},u_0)\right\|_{H^2(\Theta)},\\
&\left\|M_{\alpha_1\alpha_2}^{h_1}(\bm{y},u_0)-M_{\alpha_1\alpha_2}(\bm{y},u_0)\right\|_{H^m(\Theta)}\leq Ch_1^{2-m}\left\|M_{\alpha_1\alpha_2}(\bm{y},u_0)\right\|_{H^2(\Theta)},\\
&\left\|R_{\alpha_1}^{h_1}(\bm{y},u_0)-R_{\alpha_1}(\bm{y},u_0)\right\|_{H^m(\Theta)}\leq Ch_1^{2-m}\left\|R_{\alpha_1}(\bm{y},u_0)\right\|_{H^2(\Theta)},\\
&\left\|H_{\alpha_1\alpha_2}^{h_1}(\bm{y},u_0)-H_{\alpha_1\alpha_2}(\bm{y},u_0)\right\|_{H^m(\Theta)}\leq Ch_1^{2-m}\left\|H_{\alpha_1\alpha_2}(\bm{y},u_0)\right\|_{H^2(\Theta)},\\
&\left\|G_{\alpha_1\alpha_2}^{h_1}(\bm{y},u_0)-G_{\alpha_1\alpha_2}(\bm{y},u_0)\right\|_{H^m(\Theta)}\leq Ch_1^{2-m}\left\|G_{\alpha_1\alpha_2}(\bm{y},u_0)\right\|_{H^2(\Theta)},\\
&\left\|N_{\alpha_1}^{h_1}(\bm{y},u_0)-N_{\alpha_1}(\bm{y},u_0)\right\|_{H^m(\Theta)}\leq Ch_1^{2-m}\left\|N_{\alpha_1}(\bm{y},u_0)\right\|_{H^2(\Theta)},\\
&\left\|N_{\alpha_1\alpha_2}^{h_1}(\bm{y},u_0)-N_{\alpha_1\alpha_2}(\bm{y},u_0)\right\|_{H^m(\Theta)}\leq Ch_1^{2-m}\left\|N_{\alpha_1\alpha_2}(\bm{y},u_0)\right\|_{H^2(\Theta)},\\
&\left\|Z_{\alpha_1}^{h_1}(\bm{y},u_0)-Z_{\alpha_1}(\bm{y},u_0)\right\|_{H^m(\Theta)}\leq Ch_1^{2-m}\left\|Z_{\alpha_1}(\bm{y},u_0)\right\|_{H^2(\Theta)},\\
&\left\|W_{\alpha_1\alpha_2}^{h_1}(\bm{y},u_0)-W_{\alpha_1\alpha_2}(\bm{y},u_0)\right\|_{H^m(\Theta)}\leq Ch_1^{2-m}\left\|W_{\alpha_1\alpha_2}(\bm{y},u_0)\right\|_{H^2(\Theta)},
\end{aligned}
\end{equation}
where $m=0,1$ and $C$ is the finite element estimate constant independent of $h_1$ and dependent on $\Theta$.
\end{lemma}
$\mathbf{Proof:}$ By employing the classical finite element theory, the above inequalities is easily obtained.
\begin{lemma}
Denote $\widehat k_{ij}^{h_1}(u_{0})$ and $\widehat \sigma_{ij}^{h_1}(u_{0})$ be the finite element approximation of the corresponding homogenized parameters, the following results hold
\begin{equation}
\begin{aligned}
\left|\widehat k_{ij}^{h_1}(u_{0})\!-\!\widehat k_{ij}(u_{0})\right|\leq Ch_1^2\left\| M_j(\bm{y},u_0)\right\|_{H^2(\Theta)}^2,\;
\widetilde\gamma_0|\bm{\zeta}|^2\leq \widehat k_{ij}^{h_1}(u_{0})\zeta_i\zeta_j\leq\widetilde\gamma_1|\bm{\zeta}|^2,
\end{aligned}
\end{equation}
\begin{equation}
\begin{aligned}
\left|\widehat \sigma_{ij}^{h_1}(u_{0})\!-\!\widehat \sigma_{ij}(u_{0})\right|\leq Ch_1^2\left\| N_j(\bm{y},u_0)\right\|_{H^2(\Theta)}^2,\;
\widetilde\gamma_0|\bm{\zeta}|^2\leq \widehat \sigma_{ij}^{h_1}(u_{0})\zeta_i\zeta_j\leq\widetilde\gamma_1|\bm{\zeta}|^2,
\end{aligned}
\end{equation}
where $C$ is a constant independent of $h_1$.
\end{lemma}
$\mathbf{Proof:}$ By employing the definitions of macroscopic homogenized material parameters in (2.13), assumption $A_1$ and lemma 5.1, it follows that
\begin{equation}
\begin{aligned}
&\left|\widehat k_{ij}^{h_1}(u_{0})-\widehat k_{ij}(u_{0})\right|\\
&=\left|\frac{1}{|\Theta|}{\int_{\Theta}}\big({k_{ij}^{(0)} + k_{i\alpha_1}^{(0)}{\frac{\partial M_j^{h_1}}{\partial y_{\alpha_1}}}}\big)d\Theta-\frac{1}{|\Theta|}{\int_{\Theta}}\big({k_{ij}^{(0)} + k_{i\alpha_1}^{(0)}{\frac{\partial M_j}{\partial y_{\alpha_1}}}}\big)d\Theta\right|\\
&=\frac{1}{|\Theta|}\left|{\int_{\Theta}}\frac{\partial}{\partial y_{k}}\big(M_i^{h_1}-M_i\big)k_{kl}^{(0)}{\frac{\partial }{\partial y_{l}}\big(M_j^{h_1}-M_j\big)}d\Theta\right|\\
&\leq C\left\|M_j^{h_1}-M_j\right\|_{H^1(\Theta)}^2\leq Ch_1^2\left\|M_j\right\|_{H^2(\Theta)}^2.
\end{aligned}
\end{equation}
Furthermore, choosing a sufficiently small $h_1>0$ satisfies
\begin{equation}
Ch_1^2\left\| M_j(\bm{y},u_0)\right\|_{H^2(\Theta)}^2\leq\bar\gamma_0/2.
\end{equation}
Hence, we can verify that the lower bound in (5.2) holds
\begin{equation}
\widehat k_{ij}^{h_1}(u_{0})\zeta_i\zeta_j=\widehat k_{ij}(u_{0})\zeta_i\zeta_j+\big[\widehat k_{ij}^{h_1}(u_{0})-\widehat k_{ij}(u_{0})\big]\zeta_i\zeta_j\geq(\bar\gamma_0-\bar\gamma_0/2)\zeta_i\zeta_i=\widetilde\gamma_0|\bm{\zeta}|^2,
\end{equation}
where $\widetilde\gamma_0=\bar\gamma_0/2$ is a constant independent of $h_1$. Moreover, the upper bound in (5.2) is easily derived when setting $\widetilde\gamma_1=\bar\gamma_1+\bar\gamma_0/2$. Finally, following the similar way, we can obtain the result (5.3).

As shown in lemmas 5.1 and 5.2, the values of macroscopic homogenized material parameters $\widehat k_{ij}(u_{0})$ and $\widehat \sigma_{ij}(u_{0})$ depend on the finite element computations of
the auxiliary cell functions $M_j(\bm{y},u_0)$ and $N_j(\bm{y},u_0)$. Therefore, in practice, we need to numerically solve the modified homogenized equations as below
\begin{equation}
\left\{ \begin{aligned}
&\widehat S(u_{0}^{h_1})\frac{{\partial {u_{0}^{h_1}}(\bm{x},t)}}{{\partial t}}- \frac{\partial }{{\partial {x_i}}}\Big( {{\widehat k_{ij}^{h_1}}(u_{0}^{h_1})\frac{{\partial {u_{0}^{h_1}}(\bm{x},t)}}{{\partial {x_j}}}}\Big)\\
&\quad\quad\quad\quad\quad\quad = {\widehat \sigma _{ij}^{h_1}}(u_{0}^{h_1})\frac{\partial \phi_{0}^{h_1}(\bm{x},t)}{\partial {x_i}}\frac{\partial \phi_{0}^{h_1}(\bm{x},t)}{\partial {x_j}}+f_u(\bm{x},t),\;\;\text{in}\;\;\Omega\times(0,T),\\
&- \frac{\partial }{{\partial {x_i}}}\Big( {{\widehat \sigma _{ij}^{h_1}}(u_{0}^{h_1})\frac{{\partial \phi_0^{h_1}(\bm{x},t)}}{{\partial {x_j}}}}\Big)=f_\phi(\bm{x},t),\;\;\text{in}\;\;\Omega\times(0,T),\\
&u_0^{h_1}(\bm{x},t) = \widehat u(\bm{x},t),\;\;\text{on}\;\;\partial {\Omega_u}\times(0,T),\\
&\phi_0^{h_1}(\bm{x},t)=\widehat\phi(\bm{x},t),\;\;\text{on}\;\;\partial\Omega_{\phi}\times(0,T),\\
&{{\widehat k_{ij}^{h_1}}(u_{0}^{h_1})\frac{\partial u_0^{h_1}(\bm{x},t)}{\partial {x_j}}}{n_i} = \bar q(\bm{x},t),\;\;\text{on}\;\;\partial {\Omega_q}\times(0,T),\\
&{{\widehat \sigma _{ij}^{h_1}}(u_{0}^{h_1})\frac{\partial \phi_0^{h_1}(\bm{x},t)}{\partial {x_j}}}{n_i}=\bar{d}(\bm{x},t),\;\;\text{on}\;\;\partial\Omega_{d}\times(0,T),\\
&u_0^{h_1}({\bm{x}},0)=\widetilde u,\;\;\text{in}\;\;\Omega.
\end{aligned} \right.
\end{equation}
\begin{lemma}
Denote $u_{0}^{h_1}$ and $\phi_{0}^{h_1}$ be the exact solution of the modified homogenized equations (5.5), the following estimates hold
\begin{equation}
{\begin{aligned}
\big\|{\phi_0^{h_1}}-{\phi_0}\big\|_{L^\infty(0,T;H^1(\Omega))}\leq Ch_1^2,
\end{aligned}}
\end{equation}
\begin{equation}
{\begin{aligned}
\big\|{u_0^{h_1}}-{u_0}\big\|_{L^\infty(0,T;L^2(\Omega))}\leq Ch_1^2,
\end{aligned}}
\end{equation}
where $C$ is a constant independent of $h_1$.
\end{lemma}
$\mathbf{Proof:}$ Through subtracting the electric equation in (2.13) from corresponding electric equation in (5.7), one can directly check that
\begin{equation}
\begin{aligned}
&- \frac{\partial }{{\partial {x_i}}}\Big[ {{\widehat \sigma _{ij}^{h_1}}(u_{0}^{h_1})\frac{{\partial \big(\phi_0^{h_1}-\phi_0\big)}}{{\partial {x_j}}}}\Big]=- \frac{\partial }{{\partial {x_i}}}\Big[ {{\big(\widehat \sigma _{ij}(u_{0})-\widehat \sigma _{ij}^{h_1}(u_{0}^{h_1})\big)}\frac{{\partial \phi_0}}{{\partial {x_j}}}}\Big]\\
&=- \frac{\partial }{{\partial {x_i}}}\Big[ {{\big(\widehat \sigma _{ij}(u_{0})-\widehat \sigma _{ij}^{h_1}(u_{0})\big)}\frac{{\partial \phi_0}}{{\partial {x_j}}}}\Big]- \frac{\partial }{{\partial {x_i}}}\Big[ {{\big(\widehat \sigma _{ij}^{h_1}(u_{0})-\widehat \sigma _{ij}^{h_1}(u_{0}^{h_1})\big)}\frac{{\partial \phi_0}}{{\partial {x_j}}}}\Big]
\end{aligned}
\end{equation}
Next, multiplying on both sides of equality (5.10) by $\phi_0^{h_1}-\phi_0$ and integrating on $\Omega$, we derive the following equality.
\begin{equation}
\begin{aligned}
&\int_{\Omega}{{\widehat \sigma _{ij}^{h_1}}(u_{0}^{h_1})\frac{{\partial \big(\phi_0^{h_1}-\phi_0\big)}}{{\partial {x_j}}}}\frac{\partial \big(\phi_0^{h_1}-\phi_0\big)}{{\partial {x_i}}}d\Omega\\
&=\int_{\Omega}{{\big(\widehat \sigma _{ij}(u_{0})-\widehat \sigma _{ij}^{h_1}(u_{0})\big)}\frac{{\partial \phi_0}}{{\partial {x_j}}}}\frac{\partial \big(\phi_0^{h_1}-\phi_0\big)}{{\partial {x_i}}}d\Omega\\
&+\int_{\Omega} {{\big(\widehat \sigma _{ij}^{h_1}(u_{0})-\widehat \sigma _{ij}^{h_1}(u_{0}^{h_1})\big)}\frac{{\partial \phi_0}}{{\partial {x_j}}}}\frac{\partial \big(\phi_0^{h_1}-\phi_0\big)}{{\partial {x_i}}}d\Omega.
\end{aligned}
\end{equation}
Recalling the inequality in (5.3) and employing Cauchy-Schwarz inequality, we can naturally obtain the following inequality from equality (5.11) if $|\widehat \sigma _{ij}^{h_1}(u_{0})-\widehat \sigma _{ij}^{h_1}(u_{0}^{h_1})|\leq C|u_{0}-u_{0}^{h_1}|$.
\begin{equation}
{\begin{aligned}
\big\|{\phi_0^{h_1}}-{\phi_0}\big\|_{H^1(\Omega)}\leq Ch_1^2+C\big\|{u_0^{h_1}}-{u_0}\big\|_{L^2(\Omega)}.
\end{aligned}}
\end{equation}

Afterwards, subtracting the thermal equation in (2.13) from corresponding thermal equation in (5.7), we can obtain
\begin{equation}
\begin{aligned}
&\widehat S(u_{0}^{h_1})\frac{{\partial \big(u_0^{h_1}-u_0\big)}}{{\partial t}}- \frac{\partial }{{\partial {x_i}}}\Big[ {{\widehat k_{ij}^{h_1}}(u_{0}^{h_1})\frac{{\partial \big(u_0^{h_1}-u_0\big)}}{{\partial {x_j}}}}\Big]=\big[\widehat S(u_{0})-\widehat S(u_{0}^{h_1})\big]\frac{{\partial u_0}}{{\partial t}}\\
&-\frac{\partial }{{\partial {x_i}}}\Big[ {{\big(\widehat k_{ij}(u_{0})-\widehat k_{ij}^{h_1}(u_{0})\big)}\frac{{\partial u_0}}{{\partial {x_j}}}}\Big]- \frac{\partial }{{\partial {x_i}}}\Big[ {{\big(\widehat k_{ij}^{h_1}(u_{0})-\widehat k_{ij}^{h_1}(u_{0}^{h_1})\big)}\frac{{\partial u_0}}{{\partial {x_j}}}}\Big]\\
&+ {\widehat \sigma _{ij}^{h_1}}(u_{0}^{h_1})\frac{\partial \big(\phi_0^{h_1}-\phi_0\big)}{\partial {x_i}}\frac{\partial \phi_{0}^{h_1}}{\partial {x_j}}+ \big[{\widehat \sigma _{ij}^{h_1}}(u_{0}^{h_1})-{\widehat \sigma _{ij}^{h_1}}(u_{0})\big]\frac{\partial \phi_0}{\partial {x_i}}\frac{\partial \phi_{0}^{h_1}}{\partial {x_j}}\\
&+ \big[{\widehat \sigma _{ij}^{h_1}}(u_{0})-{\widehat \sigma _{ij}}(u_{0})\big]\frac{\partial \phi_0}{\partial {x_i}}\frac{\partial \phi_{0}^{h_1}}{\partial {x_j}}+ {\widehat \sigma _{ij}}(u_{0})\frac{\partial \phi_0}{\partial {x_i}}\frac{\partial \big(\phi_0^{h_1}-\phi_0\big)}{\partial {x_j}}.
\end{aligned}
\end{equation}
Furthermore, multiplying on both sides of equality (5.13) by $u_0^{h_1}-u_0$ and integrating on $\Omega$, it follows that
\begin{equation}
\begin{aligned}
&\frac{1}{2}\frac{\partial}{\partial t}\Big[\int_{\Omega}\widehat S(u_{0}^{h_1})\big(u_0^{h_1}-u_0\big)^2d\Omega\Big]+\int_{\Omega}{{\widehat k_{ij}^{h_1}}(u_{0}^{h_1})\frac{{\partial \big(u_0^{h_1}-u_0\big)}}{{\partial {x_j}}}}\frac{\partial \big(u_0^{h_1}-u_0\big) }{{\partial {x_i}}}d\Omega\\
&=\frac{1}{2}\int_{\Omega}\frac{\partial \widehat S(u_{0}^{h_1})}{\partial t}\big(u_0^{h_1}-u_0\big)^2d\Omega+\int_{\Omega}\big[\widehat S(u_{0})-\widehat S(u_{0}^{h_1})\big]\frac{{\partial u_0}}{{\partial t}}\big(u_0^{h_1}-u_0\big)d\Omega\\
&+\int_{\Omega}{{\big(\widehat k_{ij}(u_{0})-\widehat k_{ij}^{h_1}(u_{0})\big)}\frac{{\partial u_0}}{{\partial {x_j}}}}\frac{\partial \big(u_0^{h_1}-u_0\big)}{{\partial {x_i}}}d\Omega\\
&+\int_{\Omega}{{\big(\widehat k_{ij}^{h_1}(u_{0})-\widehat k_{ij}^{h_1}(u_{0}^{h_1})\big)}\frac{{\partial u_0}}{{\partial {x_j}}}}\frac{\partial \big(u_0^{h_1}-u_0\big)}{{\partial {x_i}}}d\Omega\\
&+\int_{\Omega} {\widehat \sigma _{ij}^{h_1}}(u_{0}^{h_1})\frac{\partial \big(\phi_0^{h_1}-\phi_0\big)}{\partial {x_i}}\frac{\partial \phi_{0}^{h_1}}{\partial {x_j}}\big(u_0^{h_1}-u_0\big)d\Omega\\
&+ \int_{\Omega}\big[{\widehat \sigma _{ij}^{h_1}}(u_{0}^{h_1})-{\widehat \sigma _{ij}^{h_1}}(u_{0})\big]\frac{\partial \phi_0}{\partial {x_i}}\frac{\partial \phi_{0}^{h_1}}{\partial {x_j}}\big(u_0^{h_1}-u_0\big)d\Omega\\
&+ \int_{\Omega}\big[{\widehat \sigma _{ij}^{h_1}}(u_{0})-{\widehat \sigma _{ij}}(u_{0})\big]\frac{\partial \phi_0}{\partial {x_i}}\frac{\partial \phi_{0}^{h_1}}{\partial {x_j}}\big(u_0^{h_1}-u_0\big)d\Omega\\
&+ \int_{\Omega}{\widehat \sigma _{ij}}(u_{0})\frac{\partial \phi_0}{\partial {x_i}}\frac{\partial \big(\phi_0^{h_1}-\phi_0\big)}{\partial {x_j}}\big(u_0^{h_1}-u_0\big)d\Omega.
\end{aligned}
\end{equation}
Recalling inequalities (5.2), (5.3) and (5.12), and employing Cauchy-Schwarz inequality and Young inequality, we can obtain the following inequality from equality (5.14) if $|\widehat S(u_{0})-\widehat S(u_{0}^{h_1})|\leq C|u_{0}-u_{0}^{h_1}|$ and $|\widehat k_{ij}^{h_1}(u_{0})-\widehat k_{ij}^{h_1}(u_{0}^{h_1})|\leq C|u_{0}-u_{0}^{h_1}|$.
\begin{equation}
\begin{aligned}
&\frac{\partial}{\partial t}\Big[C\big\|u_0^{h_1}-u_0\big\|_{L^2(\Omega)}^2\Big]+C\big\|u_0^{h_1}-u_0\big\|_{H^1(\Omega)}^2\\
&\leq C\big\|u_0^{h_1}-u_0\big\|_{L^2(\Omega)}^2+C\big\|u_0^{h_1}-u_0\big\|_{L^2(\Omega)}^2+\frac{C}{2\lambda}h_1^4+\frac{\lambda}{2}\big\|u_0^{h_1}-u_0\big\|_{H^1(\Omega)}^2\\
&+\frac{C}{2\lambda}\big\|u_0^{h_1}-u_0\big\|_{L^2(\Omega)}^2+\frac{\lambda}{2}\big\|u_0^{h_1}-u_0\big\|_{H^1(\Omega)}^2+C\big\|\phi_0^{h_1}-\phi_0\big\|_{H^1(\Omega)}^2\\
&+C\big\|u_0^{h_1}-u_0\big\|_{L^2(\Omega)}^2+C\big\|u_0^{h_1}-u_0\big\|_{L^2(\Omega)}^2+C\big\|u_0^{h_1}-u_0\big\|_{L^2(\Omega)}^2+Ch_1^4\\
&+C\big\|u_0^{h_1}-u_0\big\|_{L^2(\Omega)}^2+C\big\|\phi_0^{h_1}-\phi_0\big\|_{H^1(\Omega)}^2+C\big\|u_0^{h_1}-u_0\big\|_{L^2(\Omega)}^2\\
&\leq Ch_1^4+C\big\|u_0^{h_1}-u_0\big\|_{L^2(\Omega)}^2+{\lambda}\big\|u_0^{h_1}-u_0\big\|_{H^1(\Omega)}^2
\end{aligned}
\end{equation}
When choosing a sufficiently small $\lambda$ and setting $\Upsilon(t)=\big\|u_0^{h_1}-u_0\big\|_{L^2(\Omega)}^2$, then we can derive $\displaystyle\frac{d\Upsilon(t)}{dt}=Ch_1^4+C\Upsilon(t)$ from (5.15). Consequently, by taking advantage of Gronwall inequality and the arbitrariness of time variable $t$, there holds the inequality (5.9). Furthermore, combining the inequalities (5.9) and (5.12), we can easily obtain the inequality (5.8).
\begin{lemma}
Denote $u_{0}^{h_1,h_0}$ and $\phi_{0}^{h_1,h_0}$ be the finite element solutions of the modified homogenized equations (5.7) by mixed FDM-FEM scheme proposed in reference \cite{R9}, the following estimate holds
\begin{equation}
{\begin{aligned}
&\mathop{\max}_{1\leq n\leq N}\big\|{u_0^{h_1,h_0}}(\bm{x},t_n)-{u_0^{h_1}(\bm{x},t_n)}\big\|_{L^2(\Omega)}\\
&+\mathop{\max}_{1\leq n\leq N}\big\|{\phi_0^{h_1,h_0}}(\bm{x},t_{n-1/2})-{\phi_0^{h_1}(\bm{x},t_{n-1/2})}\big\|_{L^{12/5}(\Omega)}\leq C(\Delta t)^2+Ch_0^2,
\end{aligned}}
\end{equation}
%\begin{equation}
%{\begin{aligned}
%&\mathop{\max}_{1\leq n\leq N}\big\|{u_0^{h_1,h_0}}(\bm{x},t_n)-{u_0^{h_1}(\bm{x},t_n)}\big\|_{H^1(\Omega)}\\
%&+\mathop{\max}_{1\leq n\leq N}\big\|{\phi_0^{h_1,h_0}}(\bm{x},t_{n-1/2})-{\phi_0^{h_1}(\bm{x},t_{n-1/2})}\big\|_{W^{1,12/5}(\Omega)}\leq C(\Delta t)^2+Ch_0,
%\end{aligned}}
%\end{equation}
where $C$ is a constant irrespective of $h_1$ and $\Delta t$.
\end{lemma}
$\mathbf{Proof:}$ As shown in the modified homogenized equations (5.7), it satisfy the conditions of error estimates in reference \cite{R9}. Hence, employing the same proof technique as reference \cite{R9}, the estimate (5.16) can be derived.
\begin{theorem}
Let $u_{0}^{h_1,h_0}$ and $\phi_{0}^{h_1,h_0}$ be the finite element solutions of the modified homogenized equations (5.7), and $u_{0}$ and $\phi_{0}$ be the exact solutions of the homogenized equations (2.13), then the following estimate holds
\begin{equation}
{\begin{aligned}
&\mathop{\max}_{1\leq n\leq N}\big\|{u_0^{h_1,h_0}}(\bm{x},t_n)-{u_0(\bm{x},t_n)}\big\|_{L^2(\Omega)}\\
&+\mathop{\max}_{1\leq n\leq N}\big\|{\phi_0^{h_1,h_0}}(\bm{x},t_{n-1/2})-{\phi_0(\bm{x},t_{n-1/2})}\big\|_{L^{12/5}(\Omega)}\leq C(\Delta t)^2\!+\!Ch_0^2\!+\!Ch_1^2,
\end{aligned}}
\end{equation}
where $C$ is a positive constant independent of $h_1$, $h_0$ and $\Delta t$.
\end{theorem}
$\mathbf{Proof:}$ Firstly, employing the triangle inequality, there exists a inequality such that
\begin{equation}
{\begin{aligned}
&\mathop{\max}_{1\leq n\leq N}\big\|{u_0^{h_1,h_0}}(\bm{x},t_n)-{u_0(\bm{x},t_n)}\big\|_{L^2(\Omega)}\\
&+\mathop{\max}_{1\leq n\leq N}\big\|{\phi_0^{h_1,h_0}}(\bm{x},t_{n-1/2})-{\phi_0(\bm{x},t_{n-1/2})}\big\|_{L^{12/5}(\Omega)}\\
&\leq \mathop{\max}_{1\leq n\leq N}\big\|{u_0^{h_1,h_0}}(\bm{x},t_n)-{u_0^{h_1}(\bm{x},t_n)}\big\|_{L^2(\Omega)}\\
&+\mathop{\max}_{1\leq n\leq N}\big\|{\phi_0^{h_1,h_0}}(\bm{x},t_{n-1/2})-{\phi_0^{h_1}(\bm{x},t_{n-1/2})}\big\|_{L^{12/5}(\Omega)}\\
&+\mathop{\max}_{1\leq n\leq N}\big\|{u_0^{h_1}}(\bm{x},t_n)-{u_0(\bm{x},t_n)}\big\|_{L^2(\Omega)}\\
&+\mathop{\max}_{1\leq n\leq N}\big\|{\phi_0^{h_1}}(\bm{x},t_{n-1/2})-{\phi_0(\bm{x},t_{n-1/2})}\big\|_{L^{12/5}(\Omega)}.
\end{aligned}}
\end{equation}
Moreover, with the help of the embedding theorem, suffice it to have the following inequality
\begin{equation}
\begin{aligned}
&\mathop{\max}_{1\leq n\leq N}\big\|{\phi_0^{h_1}}(\bm{x},t_{n-1/2})-{\phi_0(\bm{x},t_{n-1/2})}\big\|_{L^{12/5}(\Omega)}\\
&\leq
\big\|{\phi_0^{h_1}}-{\phi_0}\big\|_{L^\infty(0,T;H^1(\Omega))}\leq Ch_1^2.
\end{aligned}
\end{equation}
Finally, substituting (5.9), (5.16) and (5.19) into (5.18), one can directly obtain the error estimate (5.17).

In a summary, the above-mentioned theoretical analysis rigorously ensures the convergence of the proposed two-stage numerical algorithm in microscopic and macroscopic computation.
\section{Numerical examples and results}
In this section, numerical examples are presented to validate the computational performance of the proposed HOMS computational model and corresponding two-stage algorithm. Moreover, all numerical experiments are conducted on a HPS desktop workstation equipped with an Intel(R) Xeon(R) Gold 6146 CPU (3.20 GHz) and internal memory (96.0 GB), and implemented based on Freefem++ software.

In addition, since it is impossible to obtain the exact solutions for nonlinear multiscale problems, we replace $u^\varepsilon(\bm{x},t)$ and $\phi^\varepsilon(\bm{x},t)$ by direct numerical simulation (DNS) solutions $u_{\text{DNS}}^\varepsilon(\bm{x},t)$ and $\phi_{{\text{DNS}}}^\varepsilon(\bm{x},t)$ for evaluating the proposed HOMS method. Furthermore, we define the following notations: $\text{Terr0}$, $\text{Terr1}$ and $\text{Terr2}$ represent the relative errors for homogenized solutions, LOMS solutions and HOMS solutions of temperature field in $L^2$ norm, $\text{TErr0}$, $\text{TErr1}$ and $\text{TErr2}$ represent the relative errors for homogenized solutions, LOMS solutions and HOMS solutions of temperature field in $H^1$ semi-norm, $\text{Perr0}$, $\text{Perr1}$ and $\text{Perr2}$ represent the relative errors for homogenized solutions, LOMS solutions and HOMS solutions of electric potential field in $L^2$ norm, $\text{PErr0}$, $\text{PErr1}$ and $\text{PErr2}$ represent the relative errors for homogenized solutions, LOMS solutions and HOMS solutions of electric potential field in $H^1$ semi-norm.
\subsection{Validation of computational accuracy and efficiency}
In this example, the nonlinear thermo-electric coupling behavior of 2D composite structure is simulated, which is modeled as a periodic array of microscopic cells, each comprising of matrix and inclusion constituents. Here, the investigated composite structure $\Omega$ is defined as $\Omega=(x_1,x_2)=[0,1]\times[0,1]$ and $\varepsilon=1/10$ (see Fig.\hspace{1mm}3).
\begin{figure}[!htb]
\centering
\begin{minipage}[c]{0.32\textwidth}
  \centering
  \includegraphics[width=0.75\linewidth,totalheight=1.2in]{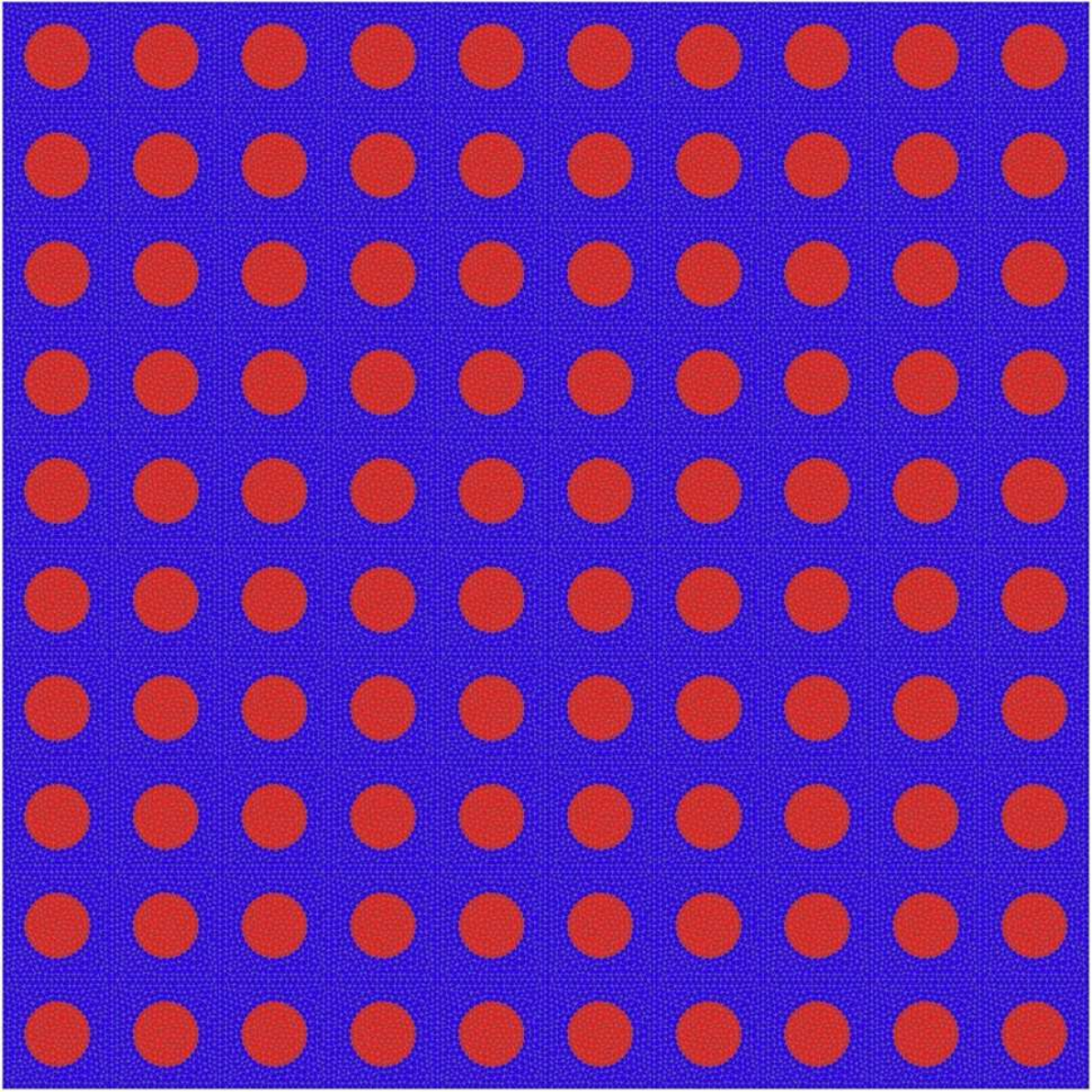} \\
  (a)
\end{minipage}
\begin{minipage}[c]{0.32\textwidth}
  \centering
  \includegraphics[width=0.75\linewidth,totalheight=1.2in]{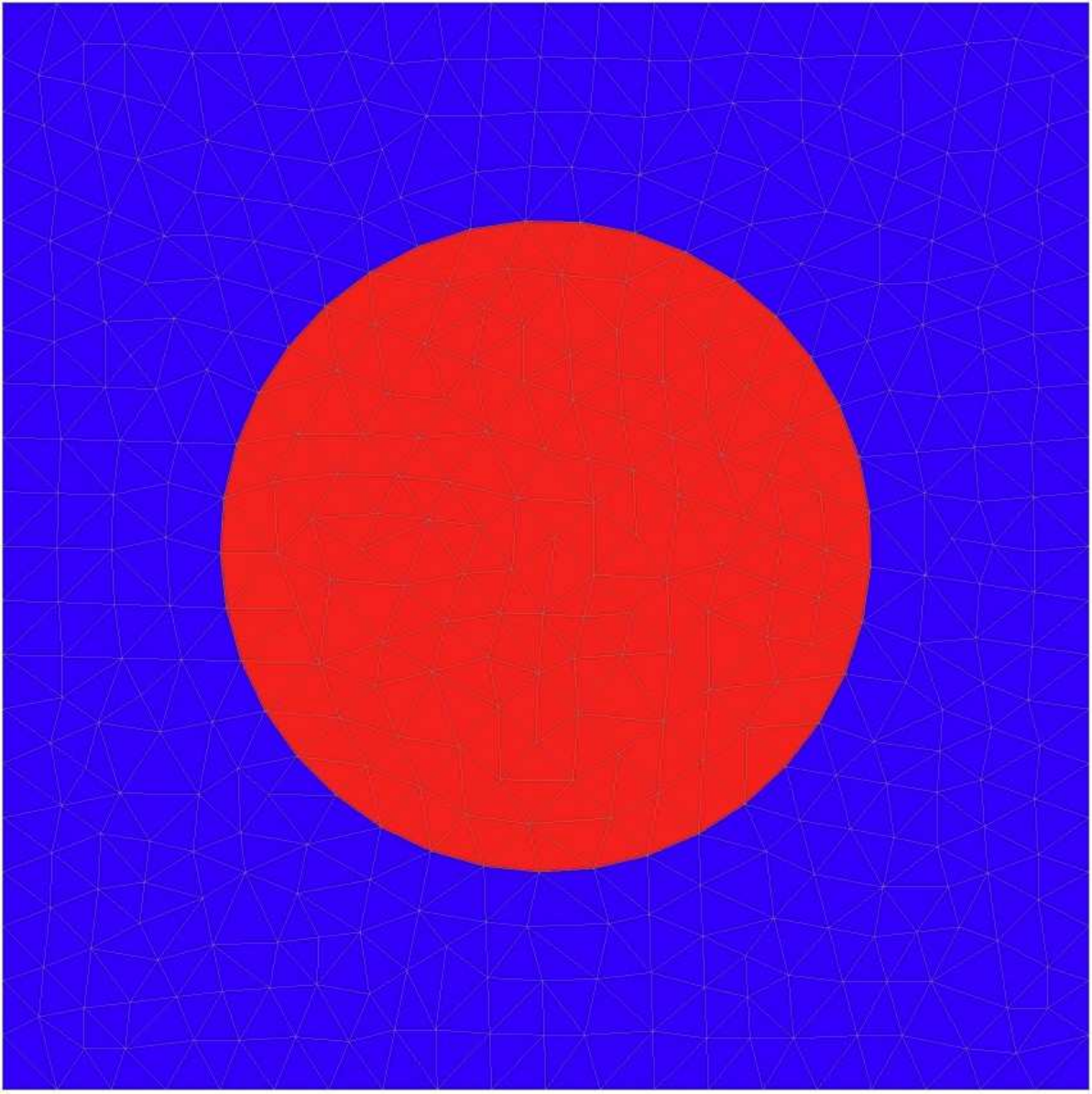} \\
  (b)
\end{minipage}
\begin{minipage}[c]{0.32\textwidth}
  \centering
  \includegraphics[width=0.75\linewidth,totalheight=1.2in]{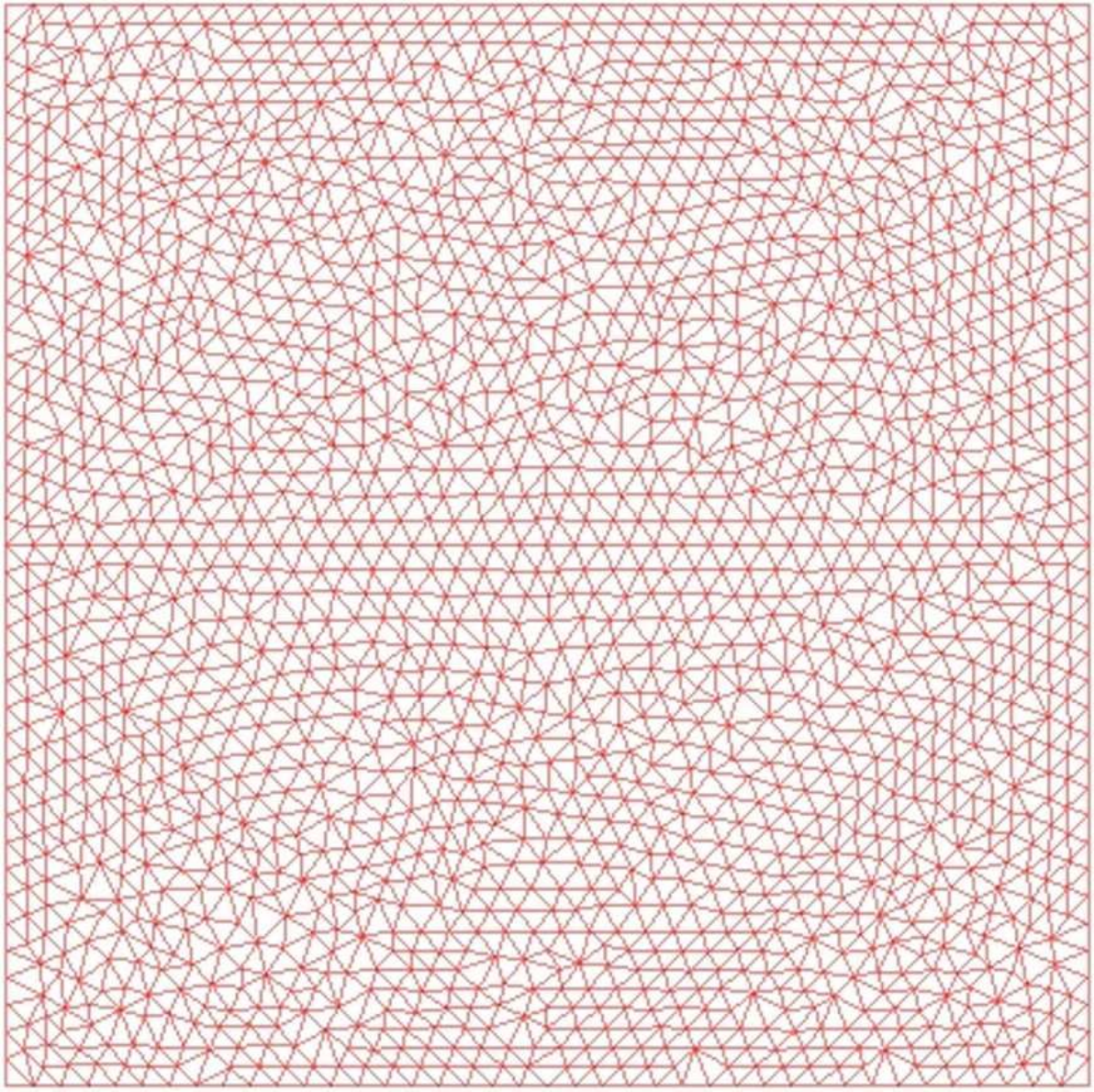} \\
  (c)
\end{minipage}
\caption{(a) The 2D composite structure $\Omega$; (b) PUC $\Theta$; (c) homogenized structure $\Omega$.}
\end{figure}

Moreover, the material parameters of composite structure are given in Table 1.
\begin{table}[h]{\caption{Material parameters of composite structure ($u$ stands for the temperature value).}}
\centering
\begin{tabular}{cc}
\hline
Material parameters & Matrix/Inclusion \\
\hline
$\rho^\varepsilon$ & 0.008/0.002  \\
$c^\varepsilon$ & 562.5/750.0  \\
$k_{ij}^\varepsilon$ & 4.0+0.0004$u$/0.04+0.000004$u$  \\
$\sigma_{ij}^{\varepsilon}$ & 300.0-0.015$u$/0.075-0.00001$u$ \\
\hline
\end{tabular}
\end{table}

In addition, the source items, initial conditions and boundary conditions in multi-scale nonlinear problem (1.1) of this example are given as follows.
\begin{equation}
\begin{aligned}
&f_u(\bm{x},t)=20000.0,\;f_\phi(\bm{x},t)=200.0,\;\;\text{in}\;\;\Omega,\\
&\widehat{{u}}(\bm{x},t)=300.0,\;\widehat{\phi}(\bm{x},t)=0.0,\;\;\text{in}\;\;\partial\Omega,\\
&\widetilde{u}(\bm{x})=300.0,\;\;\text{in}\;\;\Omega.
\end{aligned}
\end{equation}

Next, the tetrahedra finite element meshes are used for multiscale nonlinear problem (1.1), auxiliary cell problems and corresponding homogenized problem (2.13). In this example, the total auxiliary cell problems need to be solved off-line 580 times, in which the 4 first-order cell functions and 25 second-order cell functions are solved on 20 macroscopic temperature interpolation points. It is noteworthy that that the computation of auxiliary cell problems performed off-line prior to the on-line multiscale computation. The off-line computation results allows for their utilization in various composite structures made of the same material constituents. After that, macroscopic homogenized equations (2.13) and multiscale nonlinear equations (1.1) are on-line solved separately, where the temporal step is set as $\Delta t = 0.001$. The nonlinear thermo-electric coupling problem of 2D composite structure is simulated in the time interval $t\in[0,1]$. Subsequently, Table 2 gives a comparison of the numbers of FEM elements and nodes, and the computational times spent for direct finite element and multiscale simulations.
\begin{table}[!htb]{\caption{Summary of computational cost.}}
\centering
\begin{tabular}{cccc}
\hline
 & Multiscale eqs. & Cell eqs. & Homogenized eqs. \\
\hline
%\\begin{tabular}{rrrr}
FEM elements & 70800 & 856 & 3800\\
FEM nodes    & 35761 & 469 & 1981\\
\hline
& DNS & off-line stage & on-line stage \\
\hline
Computational time & 3050.549s  & 11.370s & 1967.181s\\
\hline
\end{tabular}
\end{table}

As demonstrated in Table 2, we can conclude that the computational cost of HOMS approach is far less than direct finite element simulation. The superiority of the proposed HOMS method over the full-scale DNS is obvious since a highly fine mesh is demanded to catch the microscopic oscillatory behaviors in this heterogeneous structure. Comparatively, the proposed HOMS approach achieves a significantly accelerated simulation for multiscale nonlinear problem, which can economize about $35.14\%$ computational time.

Figs.\hspace{1mm}4 and 5 illustrate the simulative results for solutions $u_{0}$, $u^{(1\varepsilon)}$, $u^{(2\varepsilon)}$, $u_{\rm{DNS}}^\varepsilon$, and $\phi_{0}$, $\phi^{(1\varepsilon)}$, $\phi^{(2\varepsilon)}$, $\phi_{{\rm{DNS}}}^\varepsilon$ at the final moment $t$=1.0, respectively. Furthermore, Fig.\hspace{1mm}6 displays the evolutive relative errors of temperature and electric potential fields in the $L^2$ norm and $H^1$ semi-norm senses.
\begin{figure}[!htb]
\centering
\begin{minipage}[c]{0.4\textwidth}
  \centering
  \includegraphics[width=45mm]{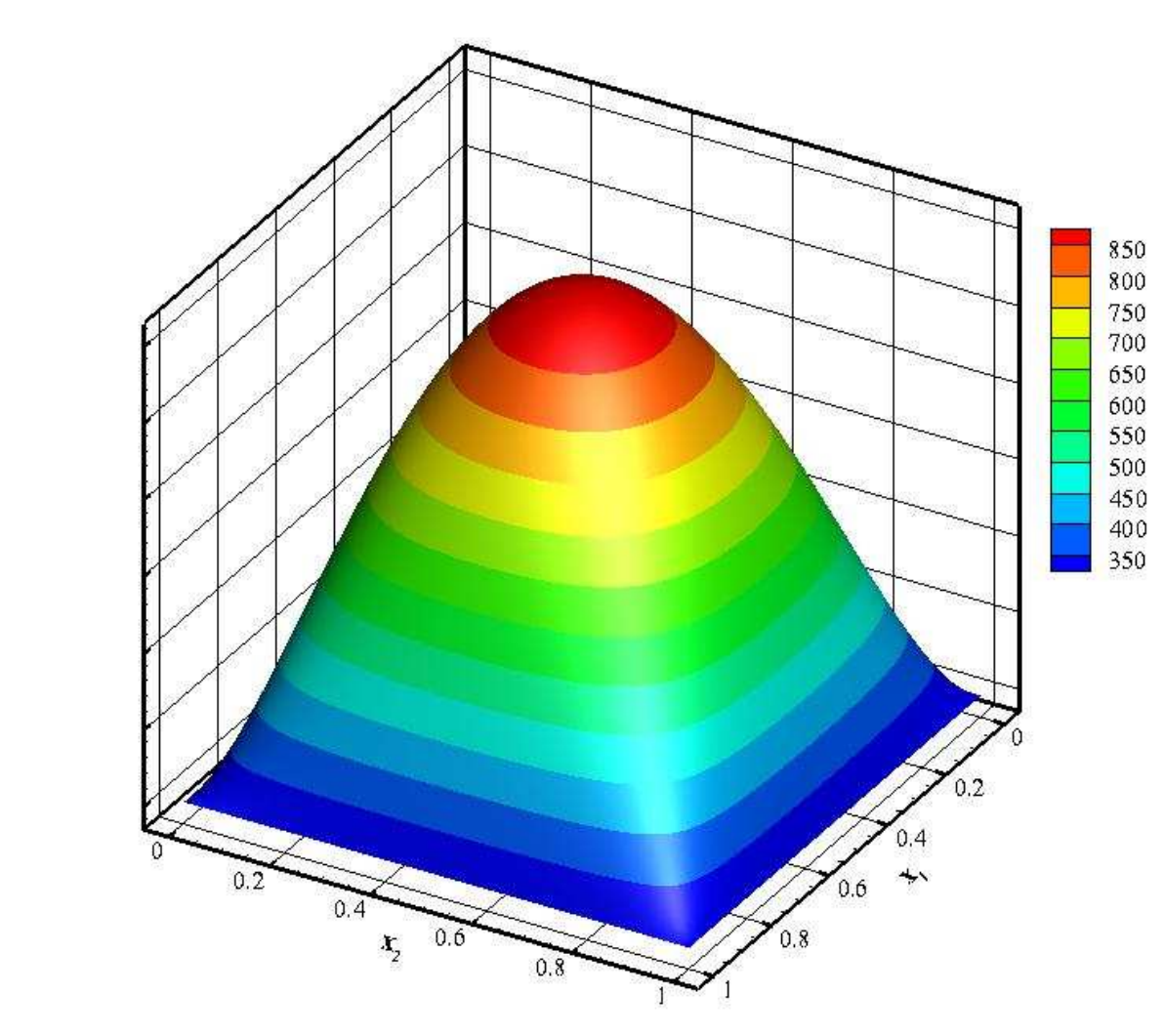}\\
  (a)
\end{minipage}
\begin{minipage}[c]{0.4\textwidth}
  \centering
  \includegraphics[width=45mm]{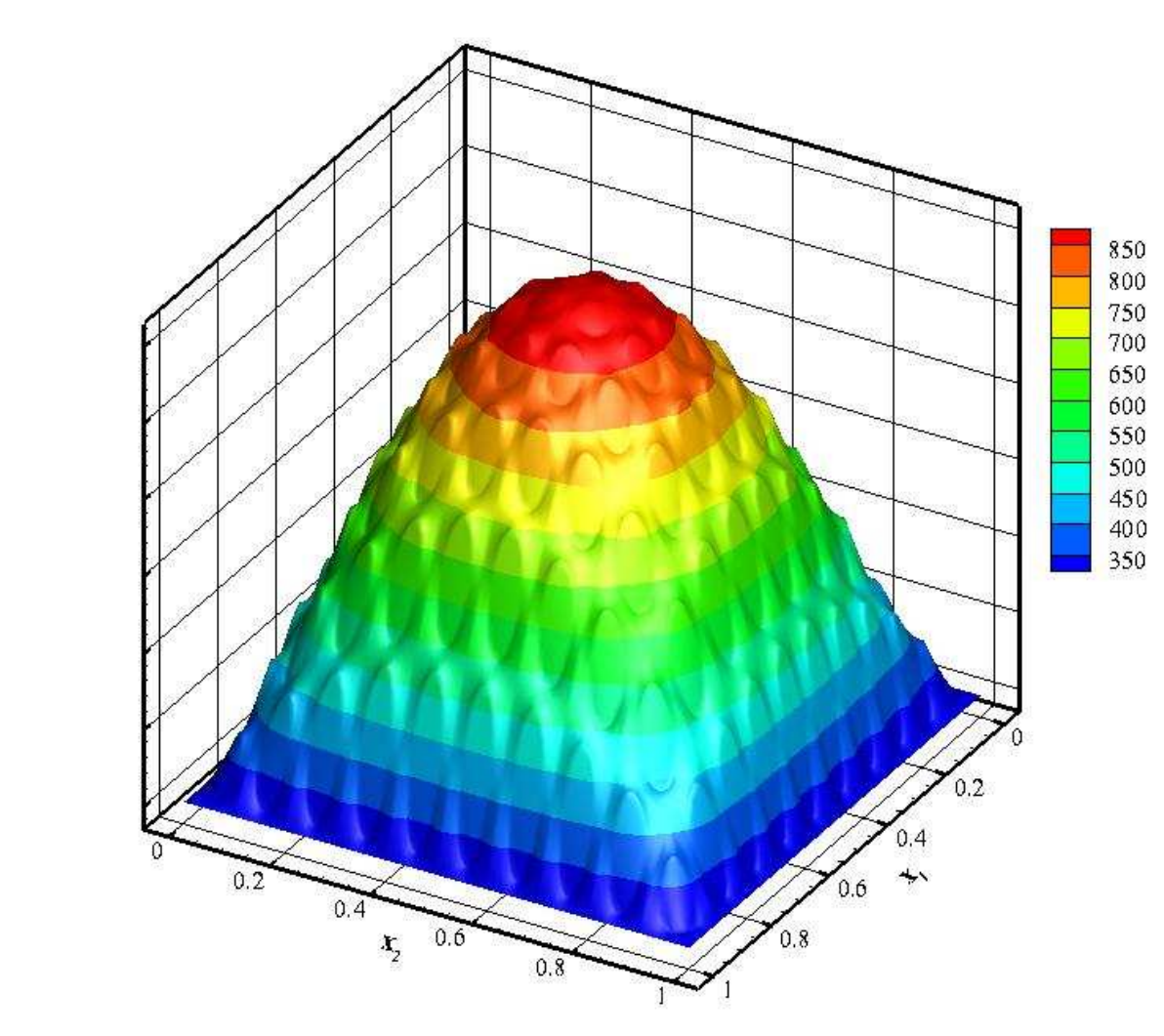}\\
  (b)
\end{minipage}
\begin{minipage}[c]{0.4\textwidth}
  \centering
  \includegraphics[width=45mm]{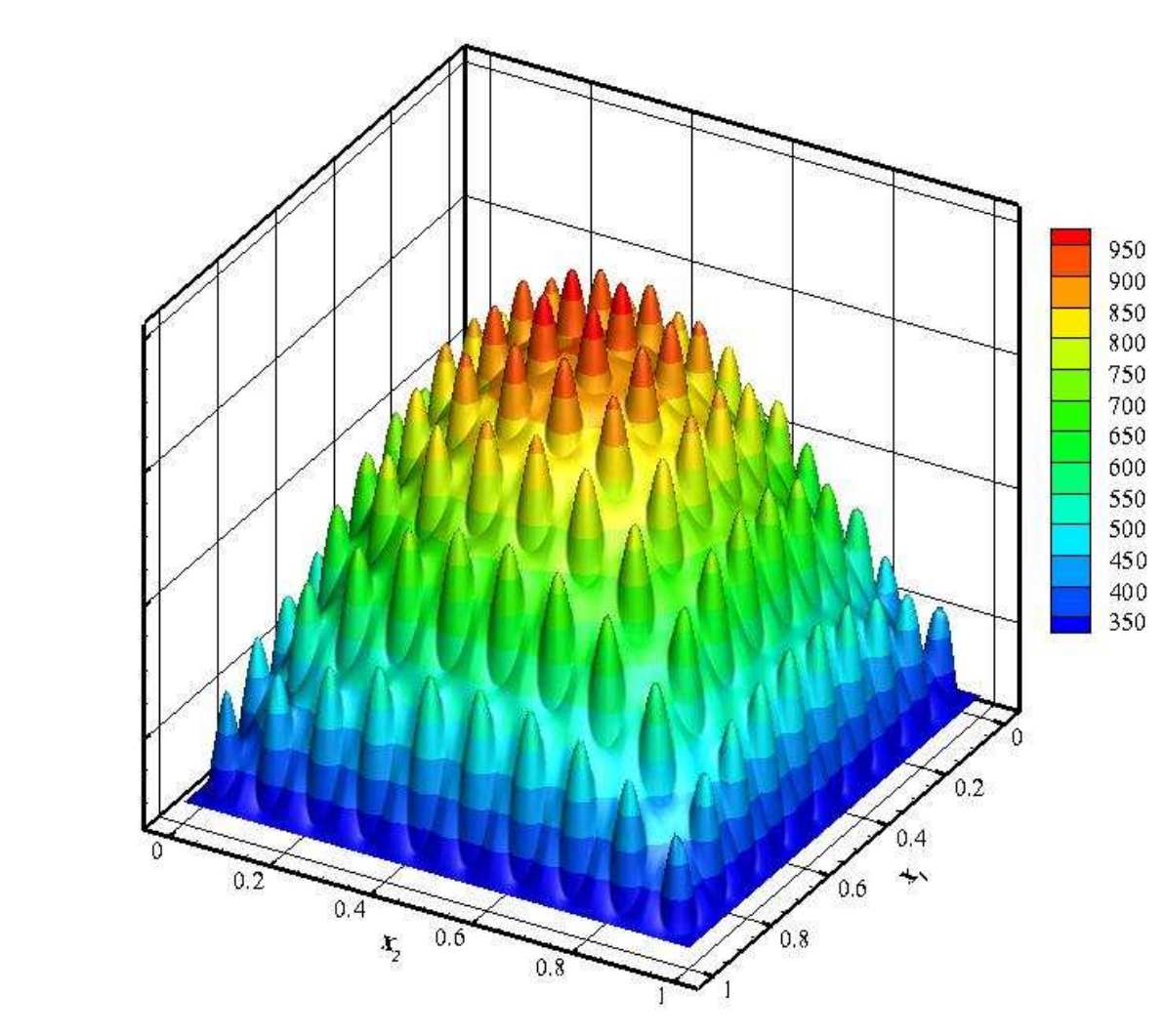}\\
  (c)
\end{minipage}
\begin{minipage}[c]{0.4\textwidth}
  \centering
  \includegraphics[width=45mm]{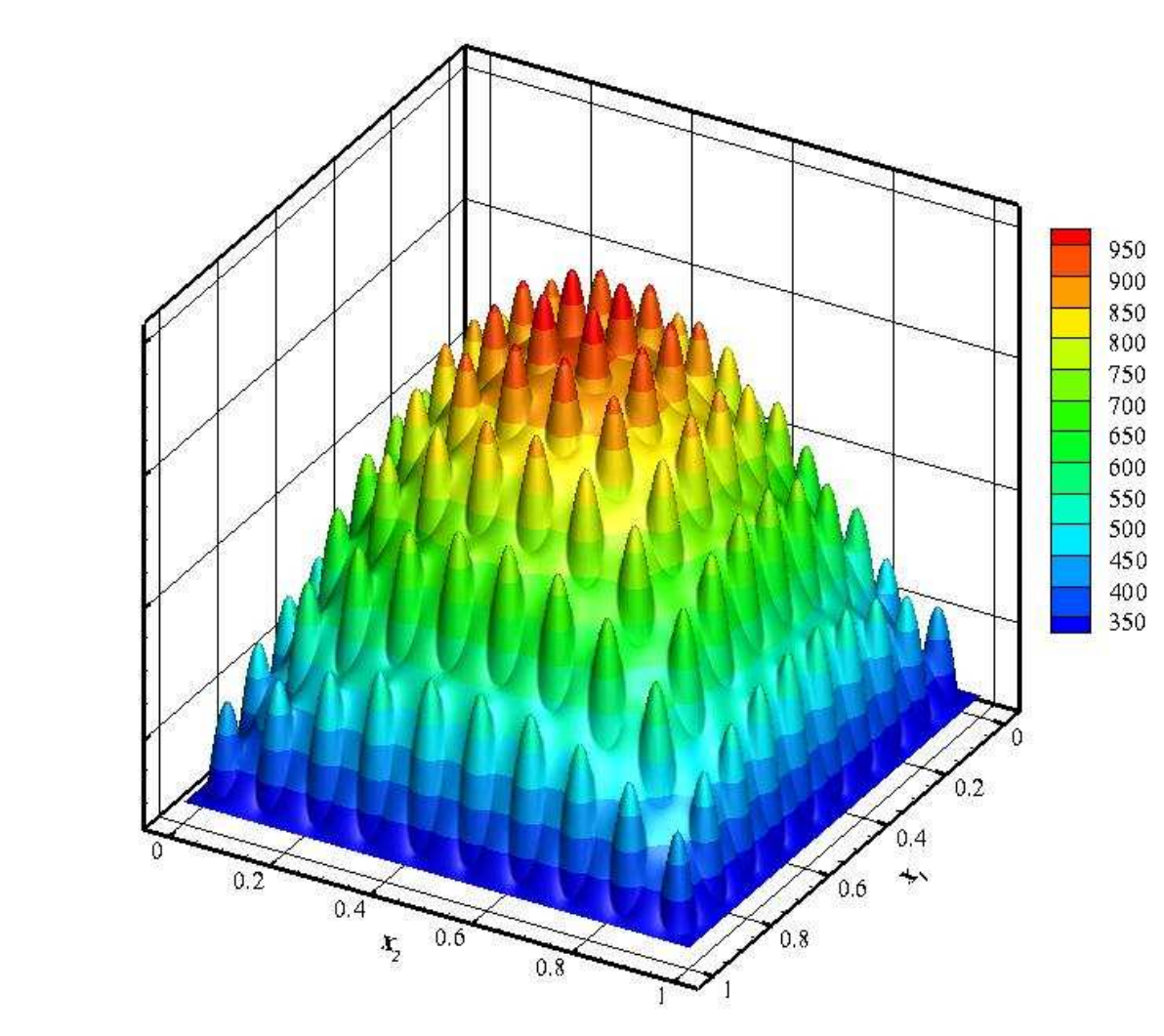}\\
  (d)
\end{minipage}
\caption{The temperature field at $t=1.0$: (a) $u_{0}$; (b) $u^{(1\varepsilon)}$; (c) $u^{(2\varepsilon)}$; (d) $u_{\rm{DNS}}^\varepsilon$.}
\end{figure}
\begin{figure}[!htb]
\centering
\begin{minipage}[c]{0.4\textwidth}
  \centering
  \includegraphics[width=45mm]{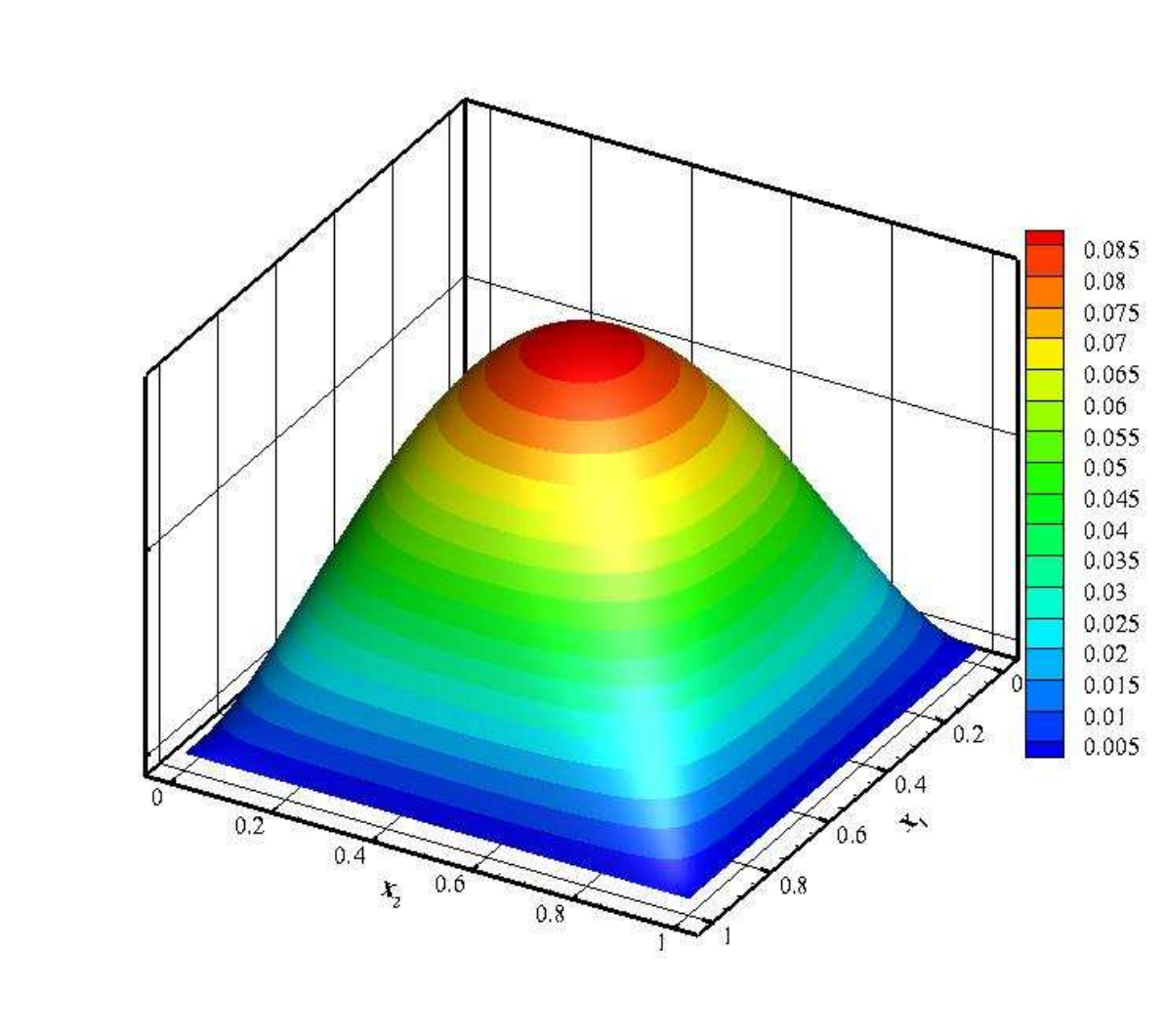}\\
  (a)
\end{minipage}
\begin{minipage}[c]{0.4\textwidth}
  \centering
  \includegraphics[width=45mm]{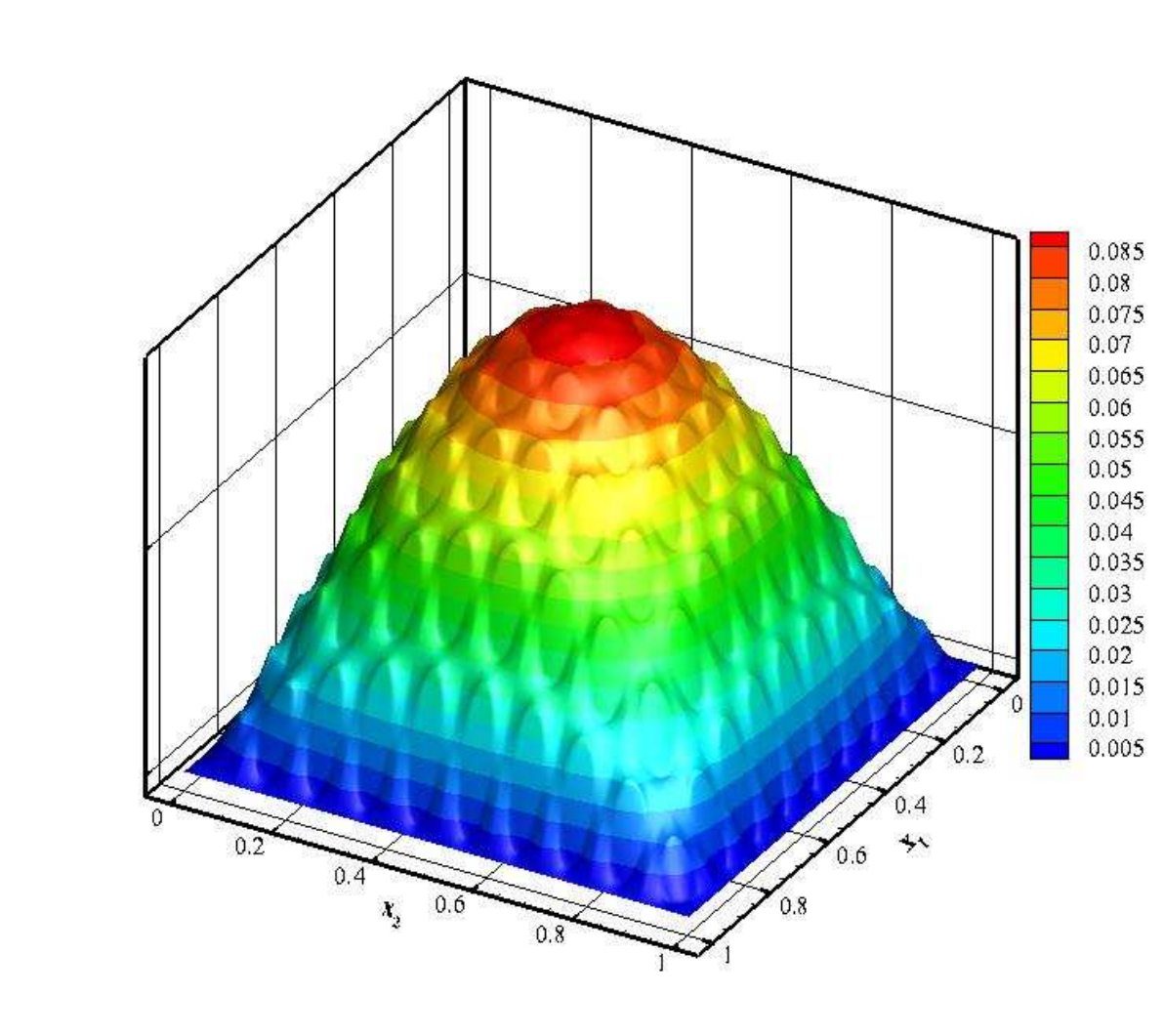}\\
  (b)
\end{minipage}
\begin{minipage}[c]{0.4\textwidth}
  \centering
  \includegraphics[width=45mm]{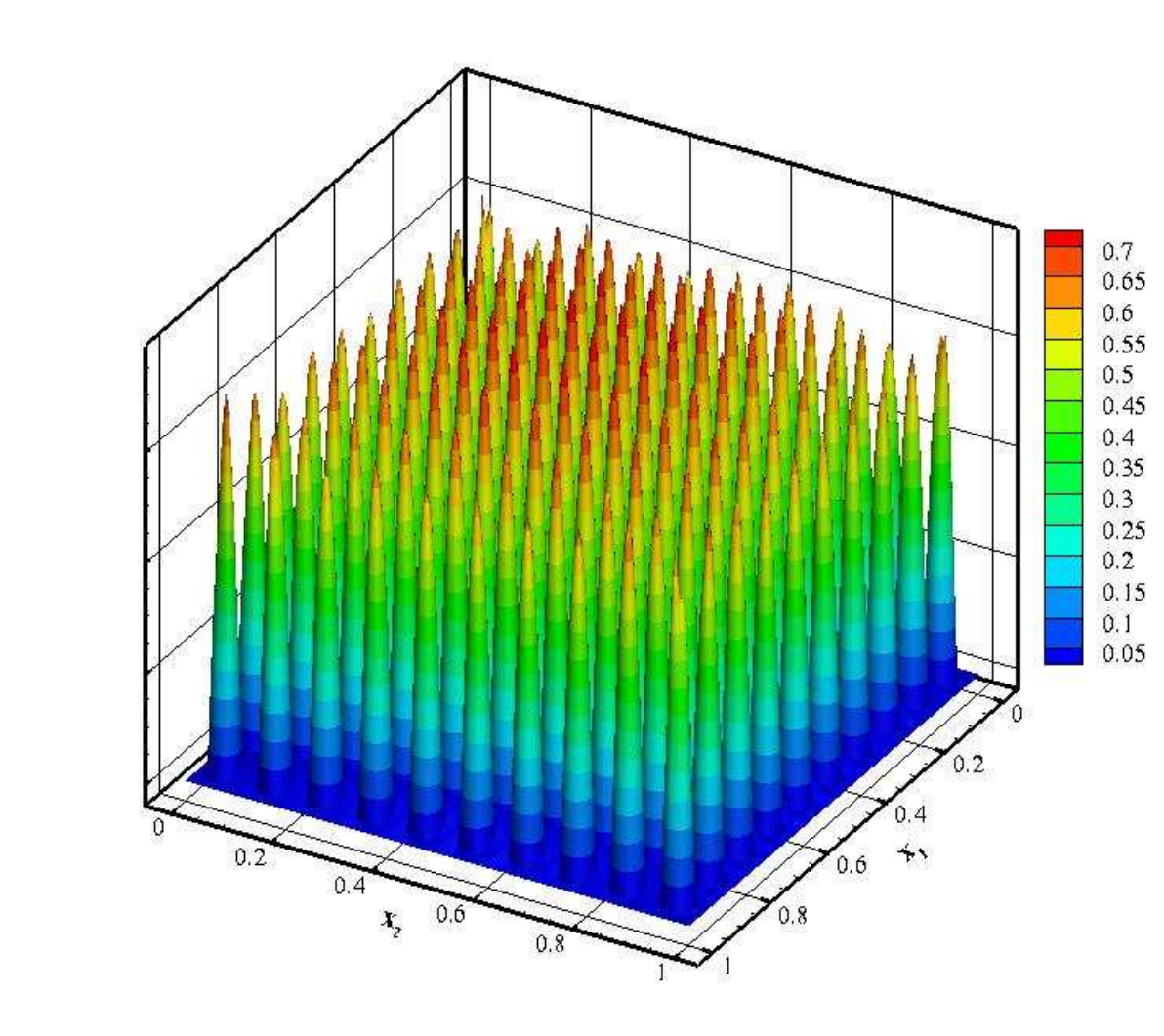}\\
  (c)
\end{minipage}
\begin{minipage}[c]{0.4\textwidth}
  \centering
  \includegraphics[width=45mm]{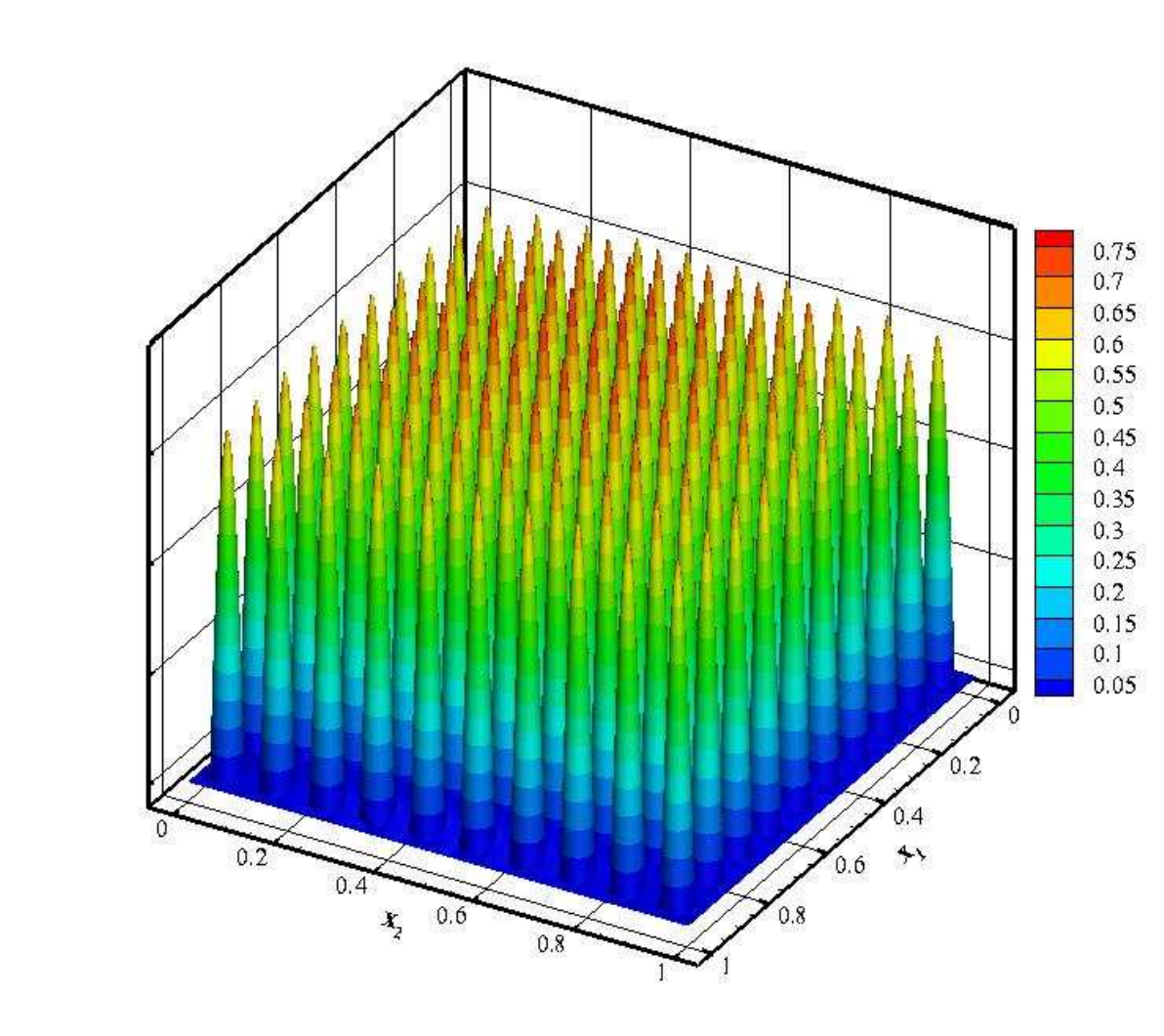}\\
  (d)
\end{minipage}
\caption{The electric potential field at $t=1.0$: (a) $\phi_{0}$; (b) $\phi^{(1\varepsilon)}$; (c) $\phi^{(2\varepsilon)}$; (d) $\phi_{{\rm{DNS}}}^\varepsilon$.}
\end{figure}
\begin{figure}[!htb]
\centering
\begin{minipage}[c]{0.4\textwidth}
  \centering
  \includegraphics[width=45mm]{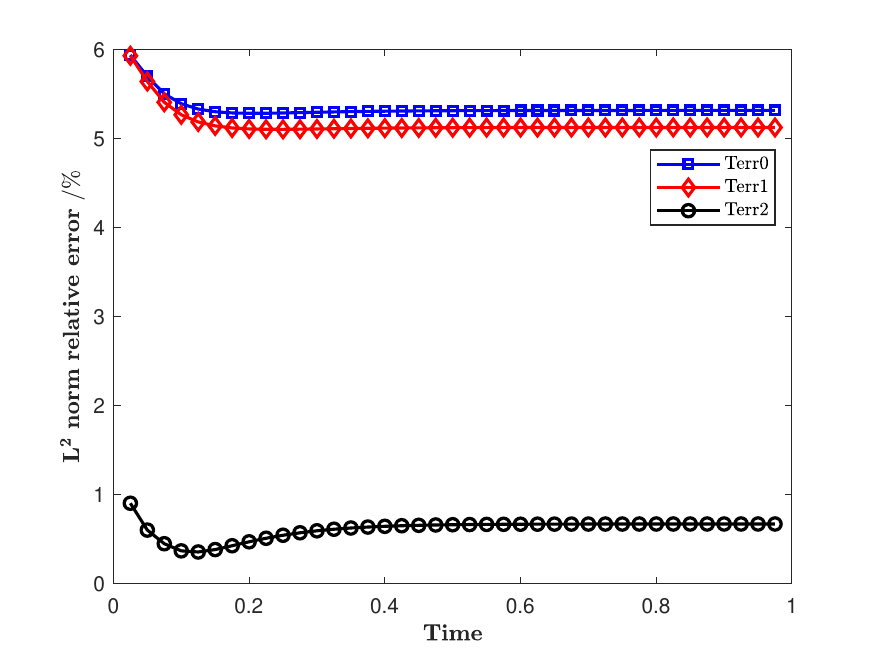}\\
  (a)
\end{minipage}
\begin{minipage}[c]{0.4\textwidth}
  \centering
  \includegraphics[width=45mm]{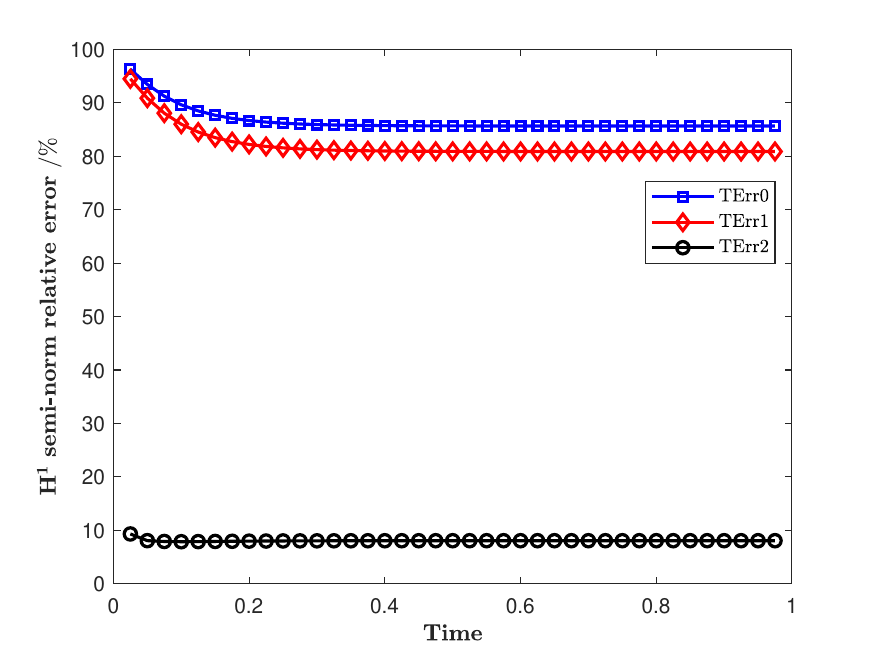}\\
  (b)
\end{minipage}
\begin{minipage}[c]{0.4\textwidth}
  \centering
  \includegraphics[width=45mm]{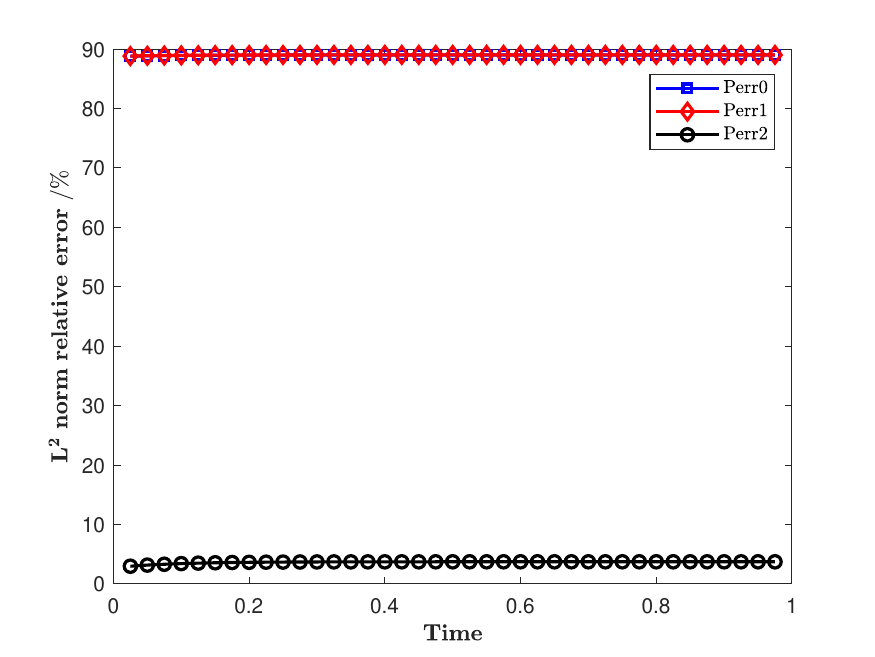}\\
  (c)
\end{minipage}
\begin{minipage}[c]{0.4\textwidth}
  \centering
  \includegraphics[width=45mm]{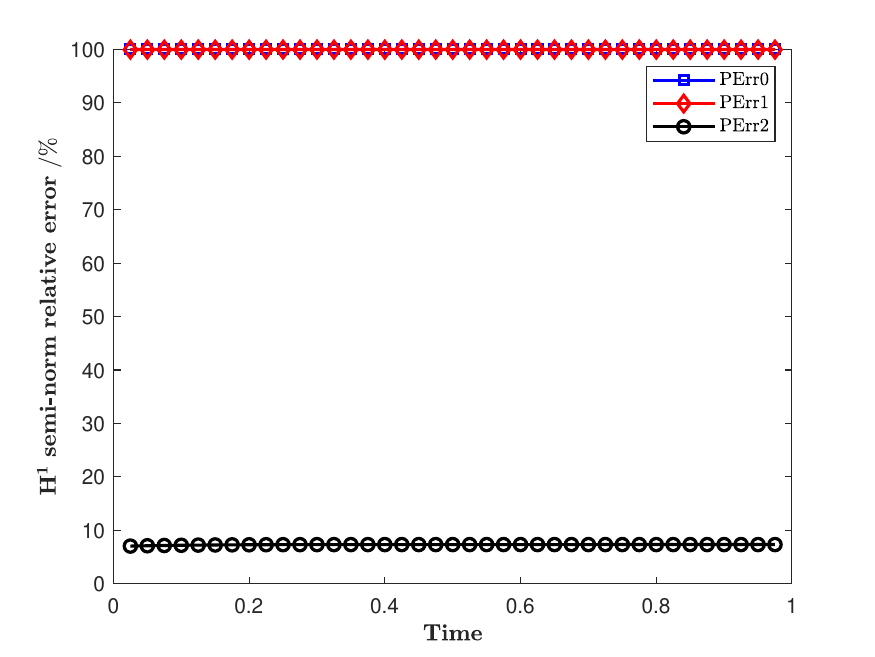}\\
  (d)
\end{minipage}
\caption{The evolutive relative errors of physical fields: (a) \rm{Terr}; (b) \rm{TErr}; (c) \rm{Perr}; (d) \rm{PErr}.}
\end{figure}

From the results in Figs.\hspace{1mm}4 and 5, HOMS solutions and full-scale simulation solutions demonstrate strong agreement, reflecting their ability to precisely capture the steeply oscillatory information within inhomogeneous structure. Conversely, homogenized and LOMS solutions can not provide enough numerical accuracy for simulating the dynamic nonlinear thermo-electric coupling problem of inhomogeneous structure. Besides, it is evident from Fig.\hspace{1mm}6 that the proposed two-stage multiscale algorithm is stable and efficient after long-time simulation. To conclude, the proposed HOMS method is ideal alternative for dynamic thermo-electric coupling simulation of composite structures.
\subsection{Application to 3D composite structure}
This example investigates the nonlinear thermo-electric coupling behaviors of 3D composite structure, which is composed of repeating unit cells with the scale separation parameter $\varepsilon=1/8$. The whole domain $\displaystyle\Omega=(x_1,x_2,x_3)=[0,1]\times[0,1]\times[0,1]$ and PUC $\Theta$ are displayed in Fig.\hspace{1mm}7.
\begin{figure}[!htb]
\centering
\begin{minipage}[c]{0.32\textwidth}
  \centering
  \includegraphics[width=0.75\linewidth,totalheight=1.2in]{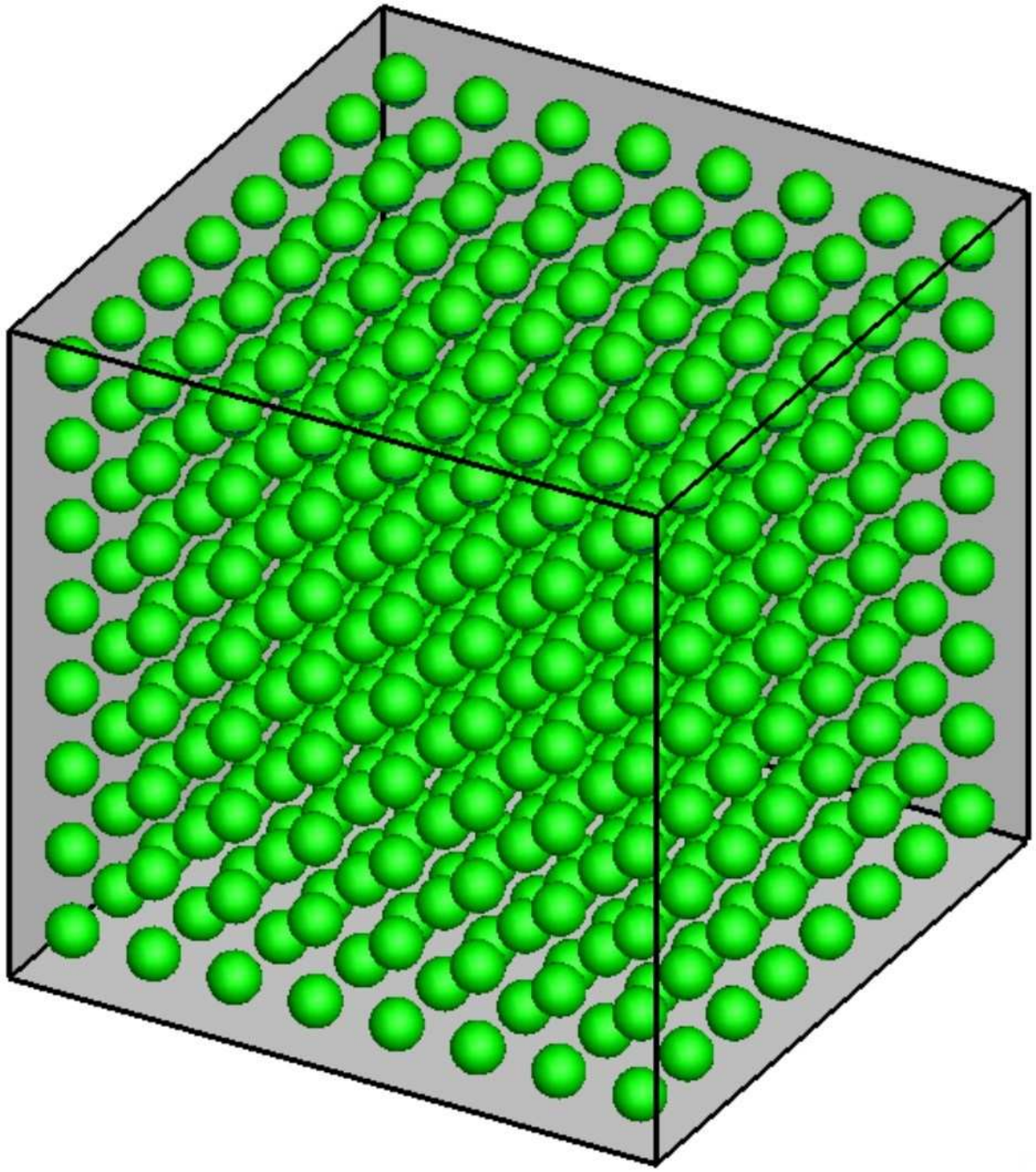} \\
  (a)
\end{minipage}
\begin{minipage}[c]{0.32\textwidth}
  \centering
  \includegraphics[width=0.75\linewidth,totalheight=1.2in]{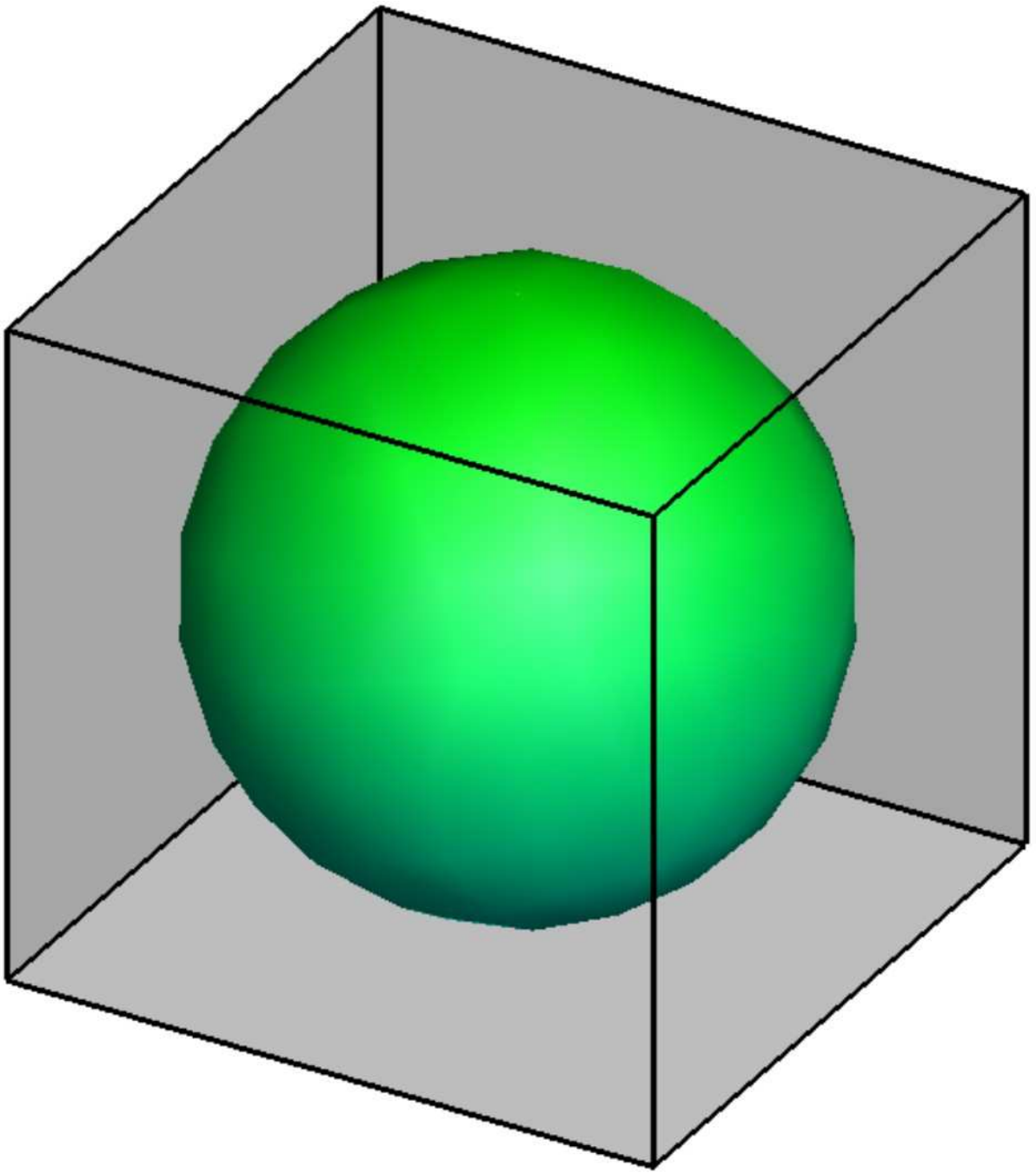} \\
  (b)
\end{minipage}
\begin{minipage}[c]{0.32\textwidth}
  \centering
  \includegraphics[width=0.75\linewidth,totalheight=1.2in]{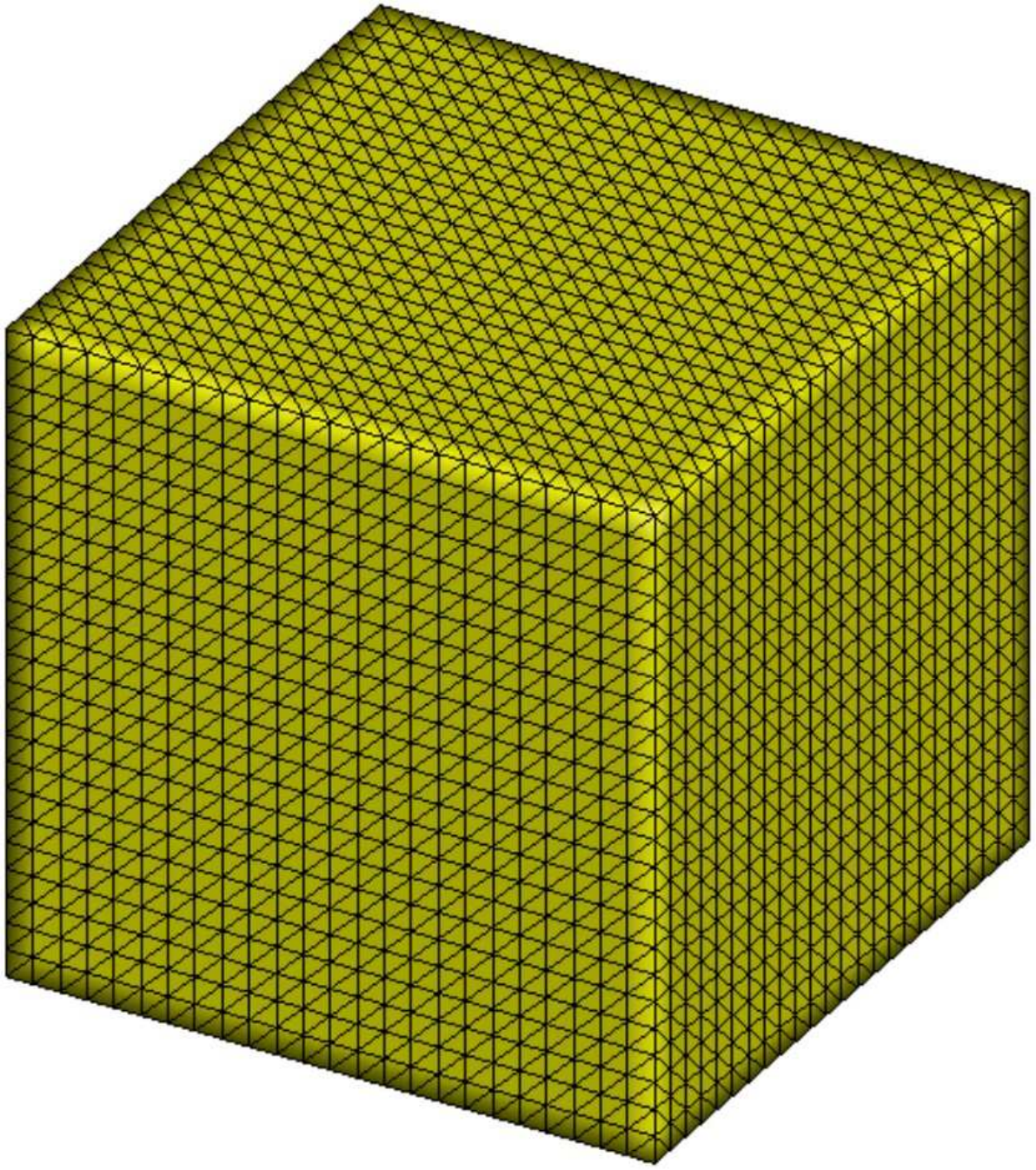} \\
  (c)
\end{minipage}
\caption{(a) The 3D composite structure $\Omega$; (b) PUC $\Theta$; (c) homogenized structure $\Omega$.}
\end{figure}

Additionally, its material parameters are defined with the same values as first example except that $\sigma_{ij}^{\varepsilon}=300.0-0.015u$ for matrix phase and $\sigma_{ij}^{\varepsilon}=0.75-0.0001u$ for inclusion phase. Moreover, the source items, initial conditions and boundary conditions in multi-scale nonlinear problem (1.1) of this example are presented as below.
\begin{equation}
\begin{aligned}
&f_u(\bm{x},t)=160000x_1(1-x_1)x_2(1-x_2),\;\;\text{in}\;\;\Omega,\\
&f_\phi(\bm{x},t)=6400x_1(1-x_1)x_2(1-x_2),\;\;\text{in}\;\;\Omega,\\
&\widehat{{u}}(\bm{x},t)=300.0,\;\widehat{\phi}(\bm{x},t)=0.0,\;\;\text{in}\;\;\partial\Omega,\\
&\widetilde{u}(\bm{x})=300.0,\;\;\text{in}\;\;\Omega.
\end{aligned}
\end{equation}

The investigated 3D heterogeneous structure in this example comprises a large number of microscopic unit cells. Direct finite element simulation needs very fine meshes and consume a significant amount of CPU time. And then, we create the finite element mesh for original multiscale equations, auxiliary cell equations and corresponding homogenized equations. The specific mesh information and computation time are presented in Table 3.
\begin{table}[!htb]{\caption{Summary of computational cost.}}
\centering
\begin{tabular}{cccc}
\hline
 & Multiscale eqs. & Cell eqs. & Homogenized eqs. \\
\hline
%\\begin{tabular}{rrrr}
FEM elements & 2924942 & 69866 & 82944\\
FEM nodes    & 468107 & 12203 & 15625\\
\hline
& DNS & off-line stage & on-line stage \\
\hline
Computational time & 179668.061s  & 1448.002s & 105521.112s\\
\hline
\end{tabular}
\end{table}

In this example, off-line 580 times computation is required for auxiliary cell problems totally, where the quantity of first-order and second-order auxiliary cell functions is set as 6 and 52 respectively. Moreover, we distribute 10 macroscopic interpolation temperature in one unit cell. The time-dependent nonlinear thermo-electric coupling responses of the 3D heterogeneous structure are simulated in the time interval $t\in[0,1]$. Setting the temporal step as $\Delta t = 0.001$, the macroscopic homogenized equations (2.13) and multiscale nonlinear equations (1.1) are on-line simulated respectively. Next, the final simulation results of temperature and electric potential fields at $t=1.0$ are depicted in Figs.\hspace{1mm}8 and 9, respectively. Besides, Fig.\hspace{1mm}10 is plotted to display the evolutionary relative errors of temperature and electric potential fields in the $L^2$ norm and $H^1$ semi-norm senses.
\begin{figure}[!htb]
\centering
\begin{minipage}[c]{0.4\textwidth}
  \centering
  \includegraphics[width=45mm]{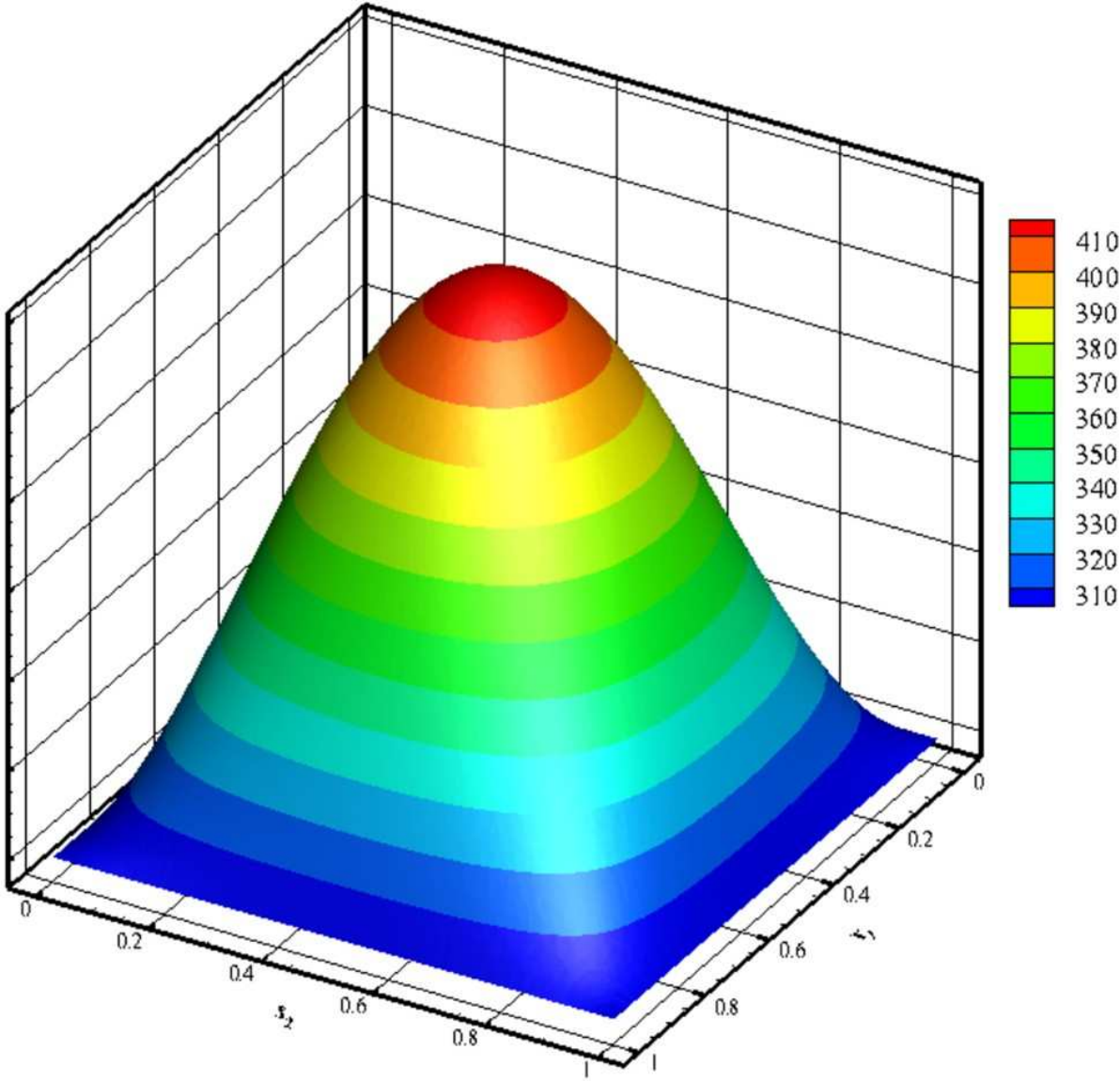}\\
  (a)
\end{minipage}
\begin{minipage}[c]{0.4\textwidth}
  \centering
  \includegraphics[width=45mm]{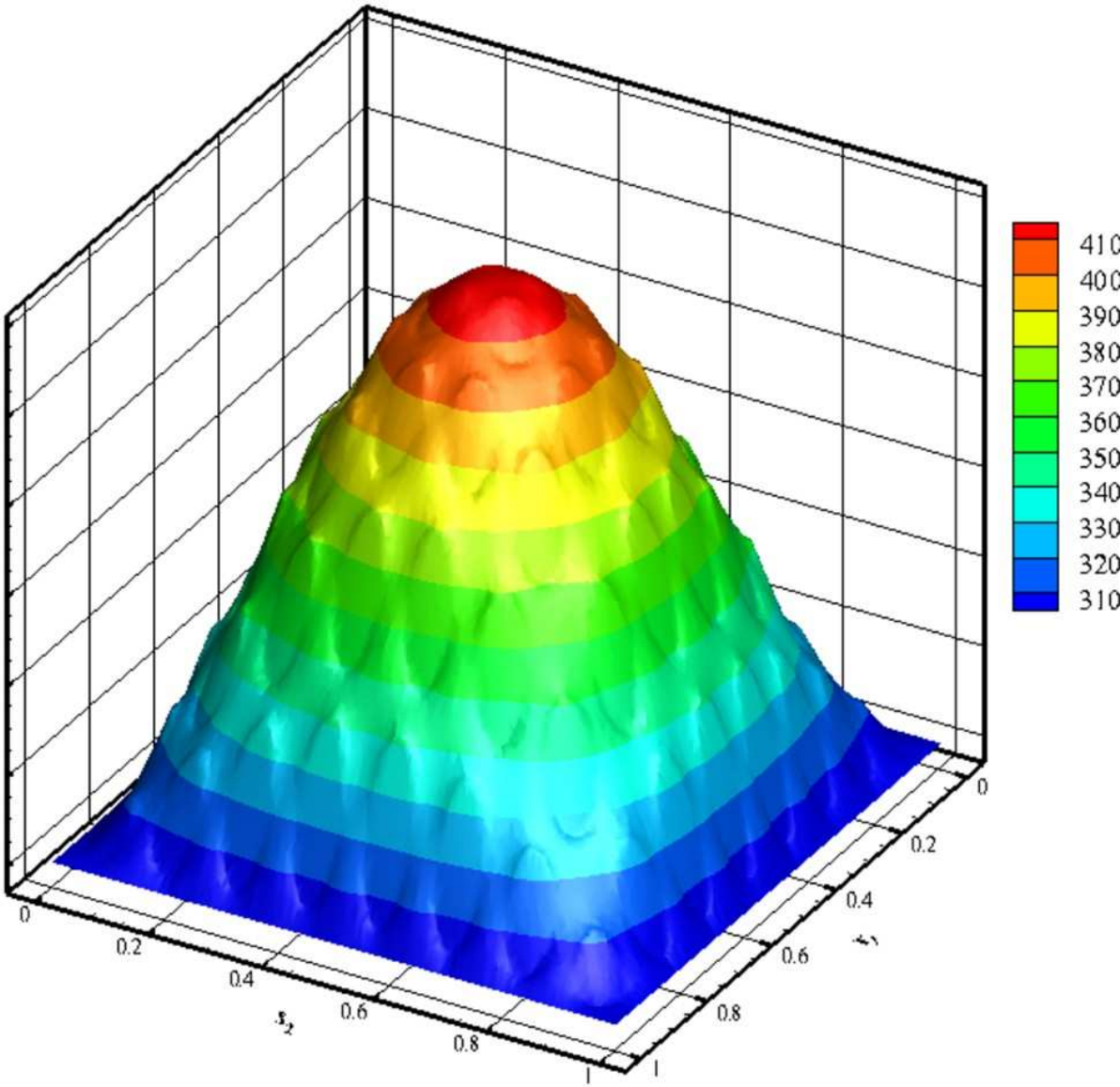}\\
  (b)
\end{minipage}
\begin{minipage}[c]{0.4\textwidth}
  \centering
  \includegraphics[width=45mm]{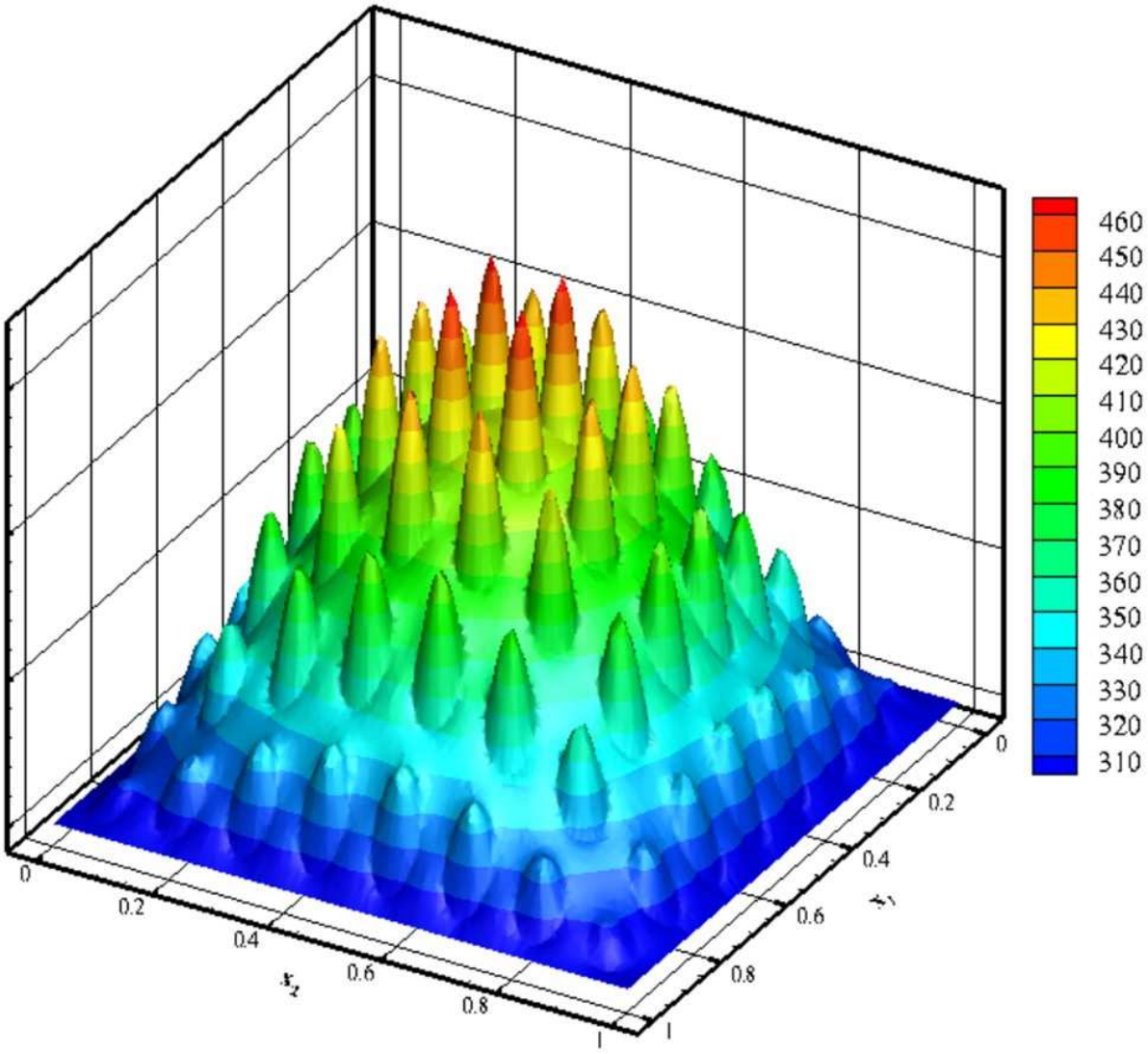}\\
  (c)
\end{minipage}
\begin{minipage}[c]{0.4\textwidth}
  \centering
  \includegraphics[width=45mm]{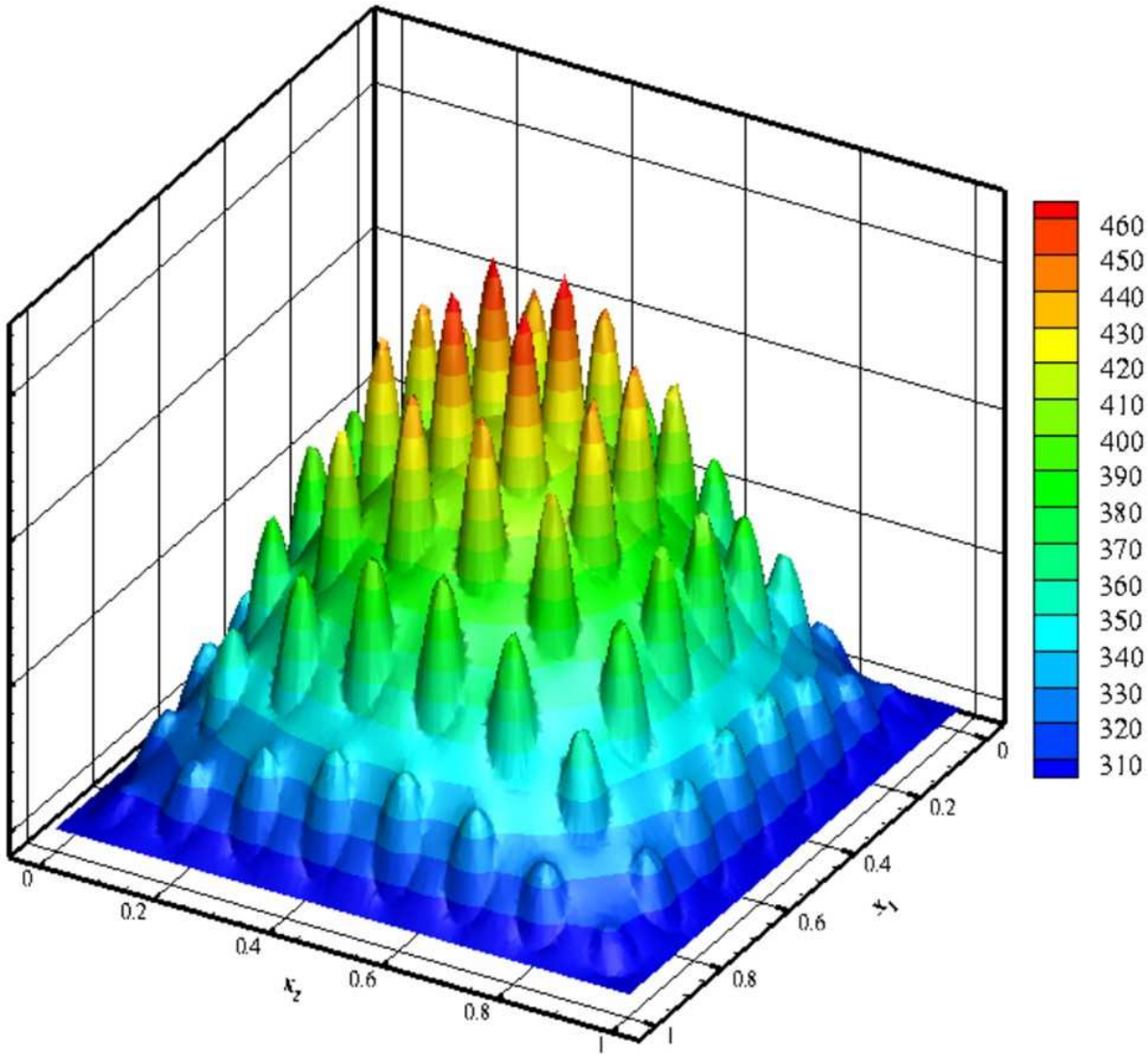}\\
  (d)
\end{minipage}
\caption{The temperature field in cross section $x_3=0.4375$ at $t=1.0$: (a) $u_{0}$; (b) $u^{(1\varepsilon)}$; (c) $u^{(2\varepsilon)}$; (d) $u_{\rm{DNS}}^\varepsilon$.}
\end{figure}
\begin{figure}[!htb]
\centering
\begin{minipage}[c]{0.4\textwidth}
  \centering
  \includegraphics[width=45mm]{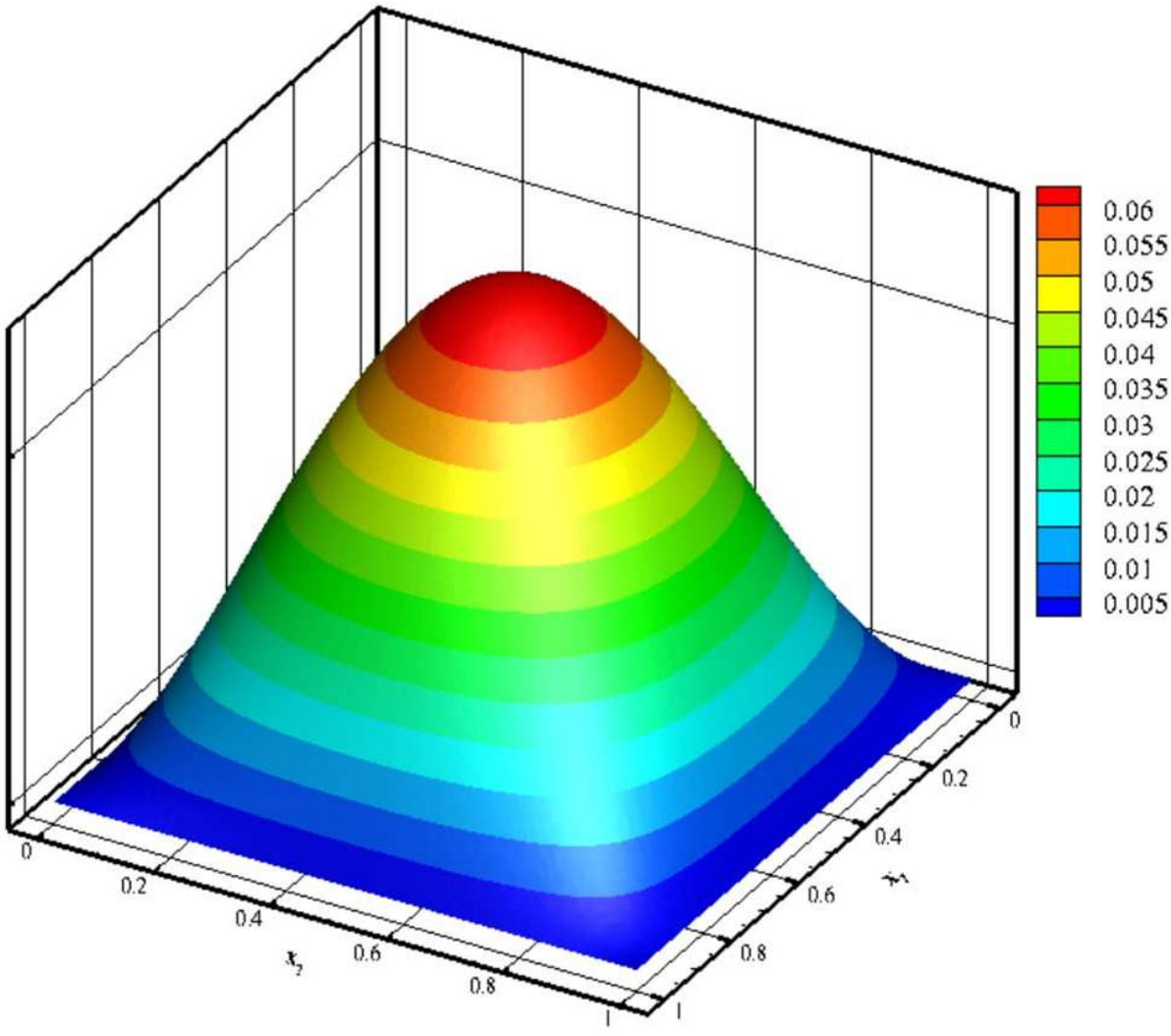}\\
  (a)
\end{minipage}
\begin{minipage}[c]{0.4\textwidth}
  \centering
  \includegraphics[width=45mm]{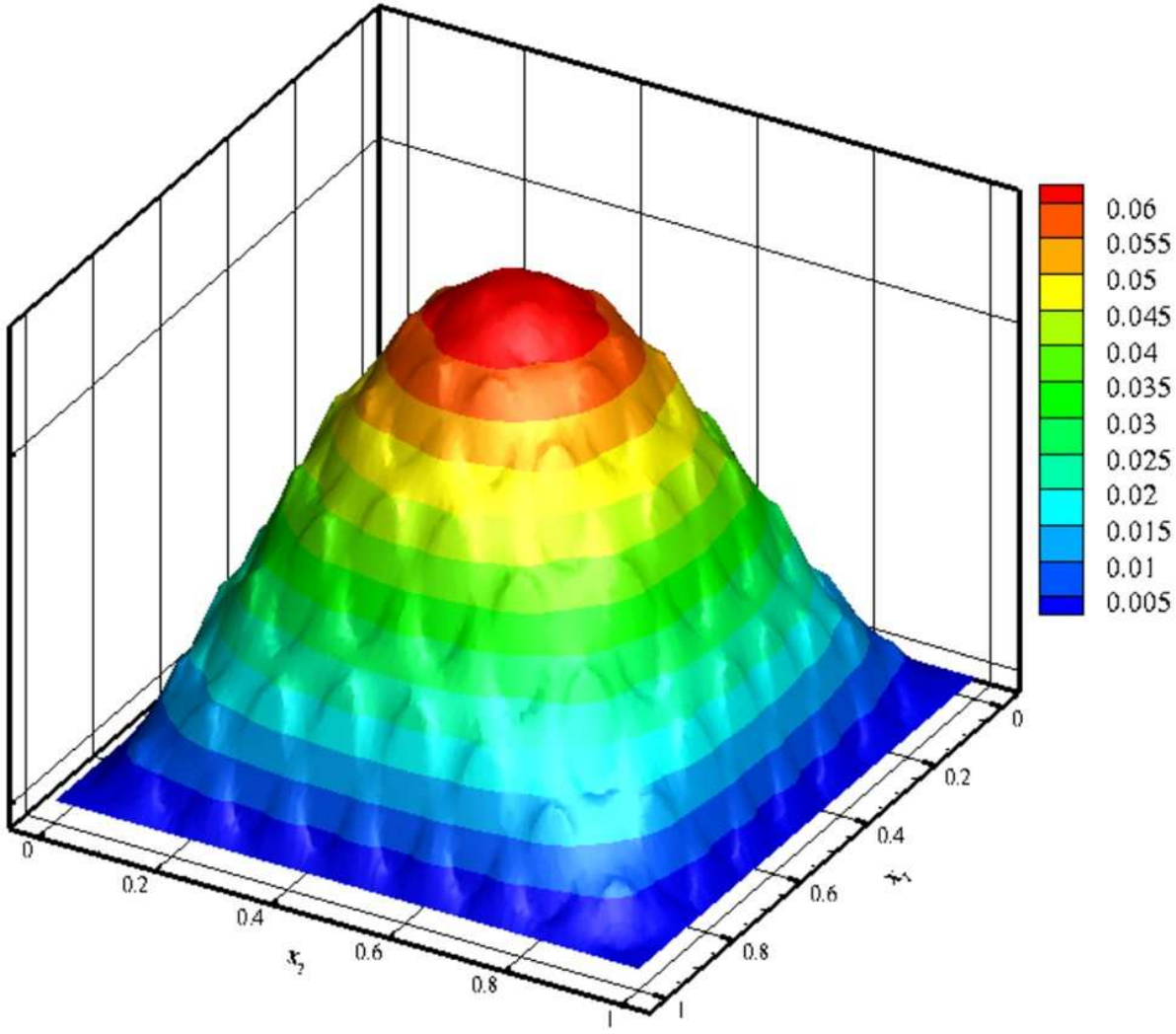}\\
  (b)
\end{minipage}
\begin{minipage}[c]{0.4\textwidth}
  \centering
  \includegraphics[width=45mm]{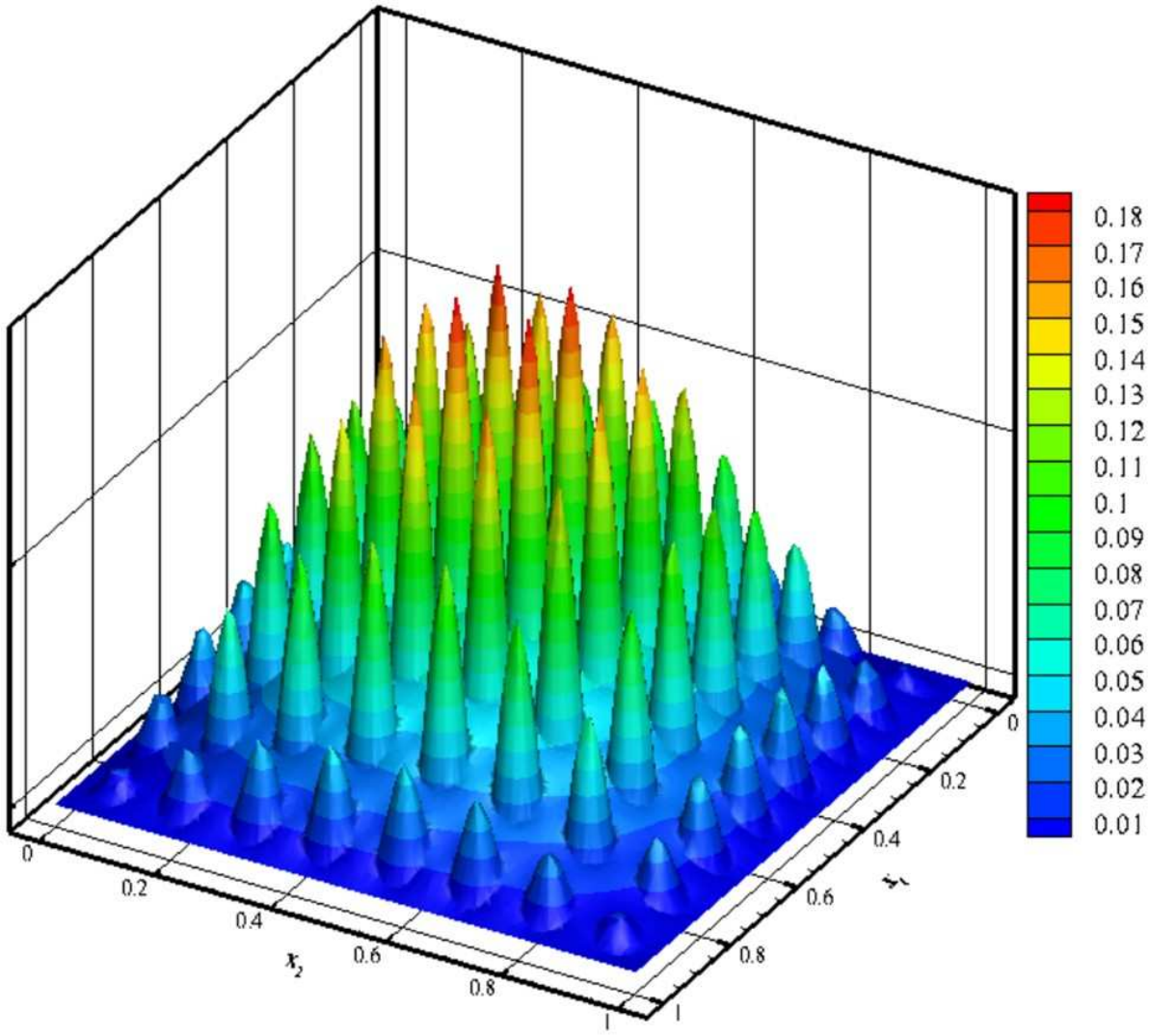}\\
  (c)
\end{minipage}
\begin{minipage}[c]{0.4\textwidth}
  \centering
  \includegraphics[width=45mm]{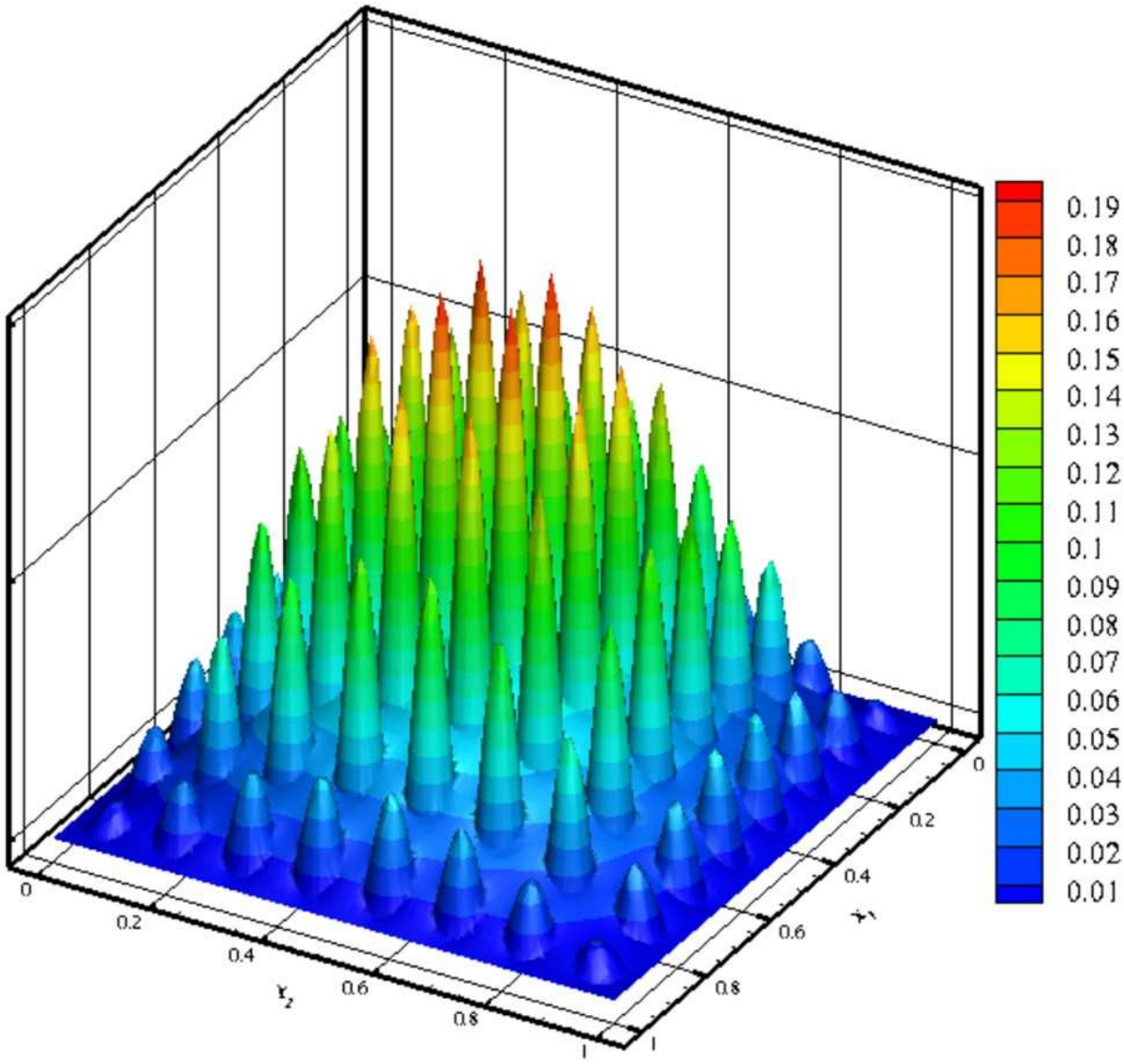}\\
  (d)
\end{minipage}
\caption{The electric potential field in cross section $x_3=0.4375$ at $t=1.0$: (a) $\phi_{0}$; (b) $\phi^{(1\varepsilon)}$; (c) $\phi^{(2\varepsilon)}$; (d) $\phi_{{\rm{DNS}}}^\varepsilon$.}
\end{figure}
\begin{figure}[!htb]
\centering
\begin{minipage}[c]{0.4\textwidth}
  \centering
  \includegraphics[width=45mm]{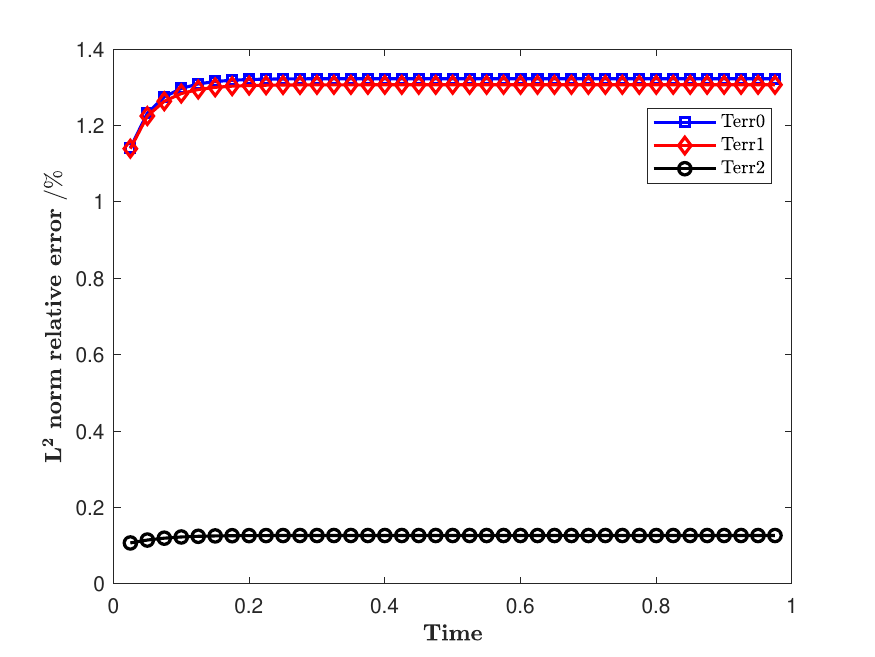}\\
  (a)
\end{minipage}
\begin{minipage}[c]{0.4\textwidth}
  \centering
  \includegraphics[width=45mm]{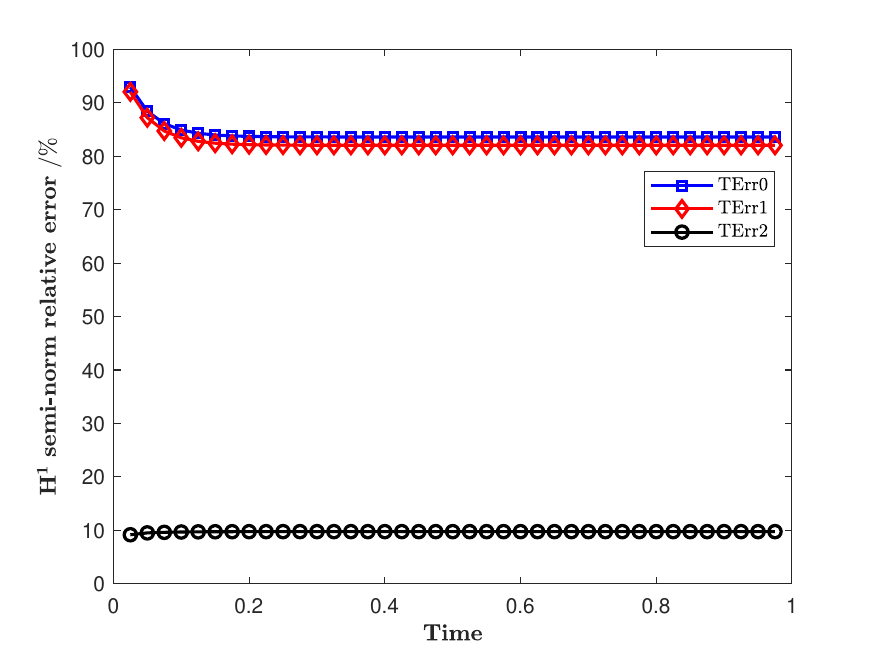}\\
  (b)
\end{minipage}
\begin{minipage}[c]{0.4\textwidth}
  \centering
  \includegraphics[width=45mm]{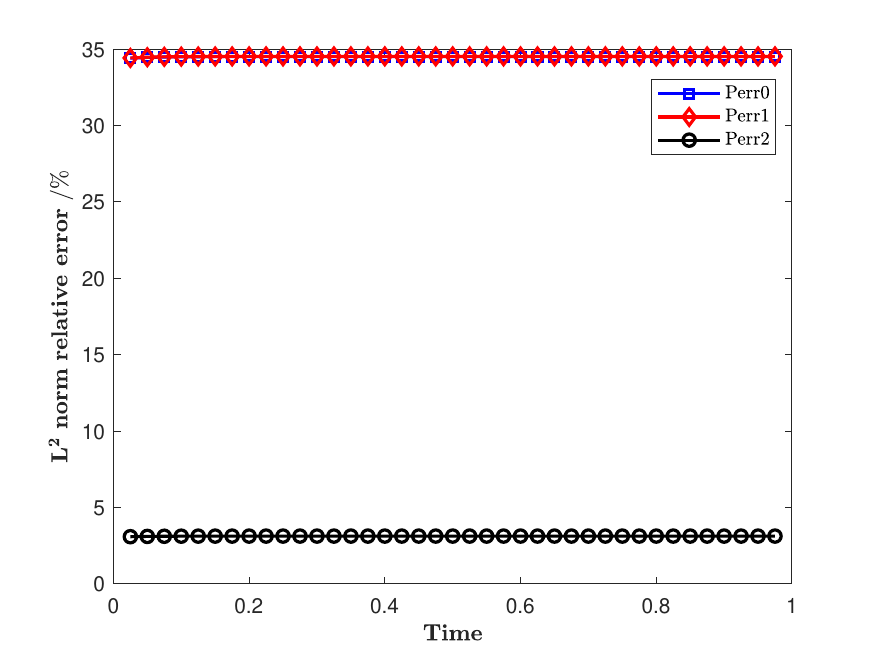}\\
  (c)
\end{minipage}
\begin{minipage}[c]{0.4\textwidth}
  \centering
  \includegraphics[width=45mm]{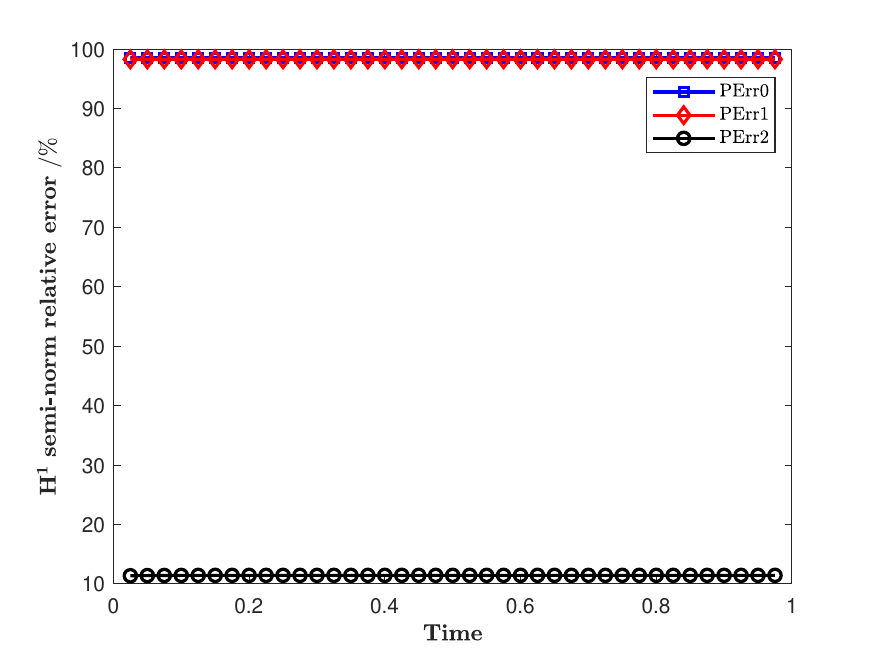}\\
  (d)
\end{minipage}
\caption{The evolutive relative errors of physical fields: (a) \rm{Terr}; (b) \rm{TErr}; (c) \rm{Perr}; (d) \rm{PErr}.}
\end{figure}

As illustrated in Table 3, the proposed HOMS approach can significantly economize computing resources by comparison with high-resolution DNS in terms of CPU memory and time. Specifically, the proposed HOMS approach can economize about $40.46\%$ computational time. It should be also mentioned that the longer multiscale simulation is implemented, the more computation time HOMS method save. From the computational results in Figs.\hspace{1mm}8 and 9, it can be clearly found that only HOMS solutions have the capacity to accurately simulate the nonlinear thermo-electric coupling behaviors of 3D heterogeneous structure, and the computational accuracy of homogenized and LOMS solutions is significantly inferior. Obviously, the homogenized method can only catch its macroscopic behaviors and the LOMS method can only capture its inadequate microscopic responses. Notably, as evidenced in Fig.\hspace{1mm}10 effectively validate that the proposed HOMS approach can retain long-time numerical stability because the blow-up phenomenon does not appear in the numerical solutions until $t=1.0$. Thus, the HOMS solutions would be utilized to catch the microscopic coupling responses of this 3D composite structure with massive unit cells. In real-world engineering applications, it is scarcely possible to obtain reference FEM solutions for large-scale composite structures. However, the HOMS method we proposed can effectively compute multiscale multiphysics problems of large-scale composite structures with minimal computational resource consumption.
\section{Conclusions}
In this work, a low-computational-cost HOMS computational method is presented for accurately simulating time-dependent nonlinear thermo-electric coupling problems of composite structures with microscopic heterogeneities. The pivotal contributions of this work are threefold: First, the novel macro-micro correlative formulations with the higher-order correction terms are established for composite structures with periodically microscopic configurations. Second, the local and global error analyses for multiscale solutions of composite structures are derived in detail. Third, an efficient two-stage algorithm with off-line and on-line stages is designed to resolve the limitation of intractable computational cost, and its convergence analysis is also derived. Numerical results testified that the proposed HOMS approach is high-accuracy and efficient for time-variant nonlinear thermo-electric coupling simulation of large-scale composite structures, and simultaneously economize the computational cost. Future work will focus on the extension of intricate nonlinear problems including thermal radiation and convection effects under high-temperature environment. Furthermore, intrinsic parallel advantage of the HOMS method in off-line stage will be exploited to further enhance its efficiency.
%\section*{Acknowledgments}
%This work was financially supported by the National Natural Science Foundation of China (No.\hspace{1mm}12001414), Young Talent Fund of Association for Science and Technology in Shaanxi, China (No.\hspace{1mm}20220506), Young Talent Fund of Association for Science and Technology in Xi'an, China (No.\hspace{1mm}095920221338), the Fundamental Research Funds for the Central Universities (No.\hspace{1mm}JB210702), the open foundation of Hubei Key Laboratory of Theory and Application of Advanced Materials Mechanics (Wuhan University of Technology) (No.\hspace{1mm}WUT-TAM202104), the National Natural Science Foundation of China (Nos.\hspace{1mm}11971386, 51739007 and 61971328). We also acknowledged the support by the Key Technology Research of FRP-Concrete Composite Structure and Xi'an Key Laboratory of Scientific Computation and Applied Statistics.
\appendix
\setcounter{equation}{0}
\renewcommand{\theequation}{A.\arabic{equation}}
\section*{Appendix A. Mathematical proof of ${\widehat \sigma _{ij}^*}(u_{0})={\widehat \sigma _{ij}}(u_{0})$}\\
As shown in (2.14), in order to prove ${\widehat \sigma _{ij}^*}(u_{0})={\widehat \sigma _{ij}}(u_{0})$, the key point lies in demonstrating that the superfluous part of ${\widehat \sigma _{ij}^*}(u_{0})$ compared to ${\widehat \sigma _{ij}}(u_{0})$ equals zero, namely
\begin{equation}
\frac{1}{|\Theta|}{\int_{\Theta}}\big( {\sigma_{\alpha_1j}^{(0)}({\bm{y}},{u_{0}}){\frac{\partial N_i({\bm{y}},{u_{0}})}{\partial y_{\alpha_1}}}}+ \sigma_{\alpha_1\alpha_2}^{(0)}({\bm{y}},{u_{0}}){\frac{\partial N_i({\bm{y}},{u_{0}})}{\partial y_{\alpha_1}}{\frac{\partial N_j({\bm{y}},{u_{0}})}{\partial y_{\alpha_2}}}}\big)d\Theta=0.
\end{equation}
Firstly, multiplying ${N_{{\alpha _2}}}({\bm{y}},{u_{0}})$ on both sides of first-order cell problem (2.12) and integrating over $\Theta$, we can derive following equality by virtue of Green's formula
\begin{equation}
\int_{\Theta}{\sigma}_{ij}^{(0)}({\bm{y}},{u_{0}})\frac{{\partial {N_{{\alpha _1}}}({\bm{y}},{u_{0}})}}{{\partial {y_j}}}\frac{\partial {N_{{\alpha _2}}}({\bm{y}},{u_{0}})}{\partial y_i}d\Theta=-\int_{\Theta} {\sigma}_{i\alpha_1}^{(0)}({\bm{y}},{u_{0}})\frac{\partial {N_{{\alpha _2}}}({\bm{y}},{u_{0}})}{\partial y_i}d\Theta.
\end{equation}
Next, some identity and index transforms are applied to equality (A.2) as below
\begin{equation}
\int_{\Theta} {\sigma}_{\alpha_1i}^{(0)}({\bm{y}},{u_{0}})\frac{\partial {N_{{\alpha _2}}}({\bm{y}},{u_{0}})}{\partial y_{\alpha_1}}d\Theta\!+\!\int_{\Theta}{\sigma}_{\alpha_1j}^{(0)}({\bm{y}},{u_{0}})\frac{{\partial {N_{i}}({\bm{y}},{u_{0}})}}{{\partial {y_j}}}\frac{\partial {N_{{\alpha _2}}}({\bm{y}},{u_{0}})}{\partial y_{\alpha_1}}d\Theta=0.
\end{equation}
\begin{equation}
\int_{\Theta} {\sigma}_{\alpha_1i}^{(0)}({\bm{y}},{u_{0}})\frac{\partial {N_{j}}({\bm{y}},{u_{0}})}{\partial y_{\alpha_1}}d\Theta+\int_{\Theta}{\sigma}_{\alpha_1\alpha_2}^{(0)}({\bm{y}},{u_{0}})\frac{{\partial {N_{i}}({\bm{y}},{u_{0}})}}{{\partial {y_{\alpha_2}}}}\frac{\partial {N_{j}}({\bm{y}},{u_{0}})}{\partial y_{\alpha_1}}d\Theta=0.
\end{equation}
\begin{equation}
\int_{\Theta} {\sigma}_{\alpha_1j}^{(0)}({\bm{y}},{u_{0}})\frac{\partial {N_{i}}({\bm{y}},{u_{0}})}{\partial y_{\alpha_1}}d\Theta+\int_{\Theta}{\sigma}_{\alpha_1\alpha_2}^{(0)}({\bm{y}},{u_{0}})\frac{{\partial {N_{j}}({\bm{y}},{u_{0}})}}{{\partial {y_{\alpha_2}}}}\frac{\partial {N_{i}}({\bm{y}},{u_{0}})}{\partial y_{\alpha_1}}d\Theta=0.
\end{equation}
Finally, dividing both sides of inequality (A.5) by a constant ${|\Theta|}$, we can successfully prove that equality (A.1) holds.
\setcounter{equation}{0}
\renewcommand{\theequation}{B.\arabic{equation}}\\
\section*{Appendix B. Detailed mathematical expressions of some functions}\\
The detailed expressions of ${{F}_{0}}({\bm{x}},{\bm{y}},t)$, ${{F}_{1}}({\bm{x}},{\bm{y}},t)$, ${{E}_{0}}({\bm{x}},{\bm{y}},t)$, ${{E}_{1}}({\bm{x}},{\bm{y}},t)$, ${{F}_{2}}({\bm{x}},{\bm{y}},t)$ and ${{E}_{2}}({\bm{x}},{\bm{y}},t)$ are presented as follows.
\begin{equation}
\begin{aligned}
&{{F}_{0}}({\bm{x}},{\bm{y}},t) = \Big[\widehat S(u_{0})}-{ \rho^{(0)}{c}^{(0)}\Big]\frac{{\partial {u_{0}}}}{{\partial t}}\\
&+ \Big[{k_{\alpha_1\alpha_2}^{(0)}}}-  { {\widehat k}_{\alpha_1\alpha_2}{+ \frac{\partial}{\partial y_i}\big( {k_{i\alpha_1}^{(0)}{M_{\alpha_2}}} \big)+{k_{\alpha_1j}\frac{\partial M_{\alpha_2}}{\partial y_j}}} \Big]\frac{\partial^2 u_{0}}{\partial x_{\alpha_1}\partial x_{\alpha_2}}\\
&+ \Big[ {\frac{\partial{k}_{i\alpha_1}^{(0)}}{\partial x_i}} - { \frac{\partial{\widehat k}_{i\alpha_1}}{\partial x_i}+ \frac{\partial}{\partial y_i}\big( {k_{ij}^{(0)}\frac{\partial M_{\alpha_1}}{\partial x_{j}}} \big)}{+\frac{\partial}{\partial x_i}\big( {k_{ij}^{(0)}\frac{\partial M_{\alpha_1}}{\partial y_{j}}} \big)} \Big]\frac{\partial u_{0}}{\partial x_{\alpha_1}}\\
&+\frac{\partial}{\partial y_i}\Big({{M_{\alpha_1}}\mathbf{D}^{(0,1)}{k_{i\alpha_2}^{(0)}}}+ {{M_{\alpha_1}}\mathbf{D}^{(0,1)}{k_{ij}^{(0)}}\frac{\partial M_{\alpha_2}}{\partial y_j}}\Big)\frac{\partial u_{0}}{\partial x_{\alpha_1}}\frac{\partial u_{0}}{\partial x_{\alpha_2}}\\
&+\Big[\sigma_{\alpha_1\alpha_2}^{(0)} - {\widehat \sigma}_{\alpha_1\alpha_2}^*
+\sigma_{i\alpha_2}^{(0)}\frac{\partial N_{\alpha_1}}{\partial y_i}+\sigma_{\alpha_1j}^{(0)}\frac{\partial N_{\alpha_2}}{\partial y_{j}}+\sigma_{ij}^{(0)}{\frac{\partial N_{\alpha_1}}{\partial y_{i}}{\frac{\partial N_{\alpha_2}}{\partial y_j}}}\Big]\frac{\partial \phi_{0}}{\partial x_{\alpha_1}}\frac{\partial \phi_{0}}{\partial x_{\alpha_2}}.
\end{aligned}
\end{equation}
\begin{equation}
\begin{aligned}
&{{F}_{1}}({\bm{x}},{\bm{y}},t) =  - {\rho^{(0)}c^{(0)}}\frac{\partial}{{\partial {t}}}\big(M_{\alpha_1}\frac{\partial u_{0}}{\partial x_{\alpha_1}}\big) {+\frac{\partial}{\partial x_i}\Big[ {k_{ij}^{(0)}\frac{\partial }{\partial x_{j}}}\big(M_{\alpha_1}\frac{\partial u_{0}}{\partial x_{\alpha_1}}\big) \Big]}\\
&+{+\frac{\partial}{\partial y_i}\Big[ {k_{ij}^{(1)}\frac{\partial }{\partial x_{j}}}\big(M_{\alpha_1}\frac{\partial u_{0}}{\partial x_{\alpha_1}}\big) \Big]}+{\frac{\partial}{\partial x_i}\Big[ {k_{ij}^{(1)}\frac{\partial }{\partial y_{j}}}\big(M_{\alpha_1}\frac{\partial u_{0}}{\partial x_{\alpha_1}}\big) \Big]}\\
&+\sigma_{ij}^{(0)}\frac{\partial \phi_{0}}{\partial x_{i}}\frac{\partial}{\partial x_{j}}\big(N_{\alpha_1}\frac{\partial \phi_{0}}{\partial x_{\alpha_1}}\big)+\sigma_{ij}^{(0)}\frac{\partial }{\partial x_{i}}\big(N_{\alpha_1}\frac{\partial \phi_{0}}{\partial x_{\alpha_1}}\big)\frac{\partial \phi_{0}}{\partial x_{j}}.
\end{aligned}
\end{equation}
\begin{equation}
\begin{aligned}
&{{E}_0}({\bm{x}},{\bm{y}},t) = \Big[{\sigma_{\alpha_1\alpha_2}^{(0)}}} -  { {\widehat \sigma}_{\alpha_1\alpha_2}{+\frac{\partial}{\partial y_i}\big( {\sigma_{i\alpha_1}^{(0)}{N_{\alpha_2}}} \big)+{\sigma_{\alpha_1j}\frac{\partial N_{\alpha_2}}{\partial y_j}}} \Big]\frac{\partial^2 \phi_{0}}{\partial x_{\alpha_1}\partial x_{\alpha_2}}\\
&+ \Big[{\frac{\partial{\sigma}_{i\alpha_1}^{(0)}}{\partial x_i}} -  { \frac{\partial{\widehat \sigma}_{i\alpha_1}}{\partial x_i}+ \frac{\partial}{\partial y_i}\big( {\sigma_{ij}^{(0)}\frac{\partial N_{\alpha_1}}{\partial x_{j}}} \big)}{+\frac{\partial}{\partial x_i}\big( {\sigma_{ij}^{(0)}\frac{\partial N_{\alpha_1}}{\partial y_{j}}} \big)} \Big]\frac{\partial \phi_{0}}{\partial x_{\alpha_1}}\\
&+\frac{\partial}{\partial y_i}\Big({{M_{\alpha_1}}\mathbf{D}^{(0,1)}{\sigma_{i\alpha_2}^{(0)}}}+ {{M_{\alpha_1}}\mathbf{D}^{(0,1)}{\sigma_{ij}^{(0)}}\frac{\partial N_{\alpha_2}}{\partial y_j}}\Big)\frac{\partial u_{0}}{\partial x_{\alpha_1}}\frac{\partial \phi_{0}}{\partial x_{\alpha_2}}.
\end{aligned}
\end{equation}
\begin{equation}
\begin{aligned}
&{{E}_1}({\bm{x}},{\bm{y}},t) =\frac{\partial}{\partial x_i}\Big[ {\sigma_{ij}^{(0)}\frac{\partial }{\partial x_{j}}}\big(N_{\alpha_1}\frac{\partial \phi_{0}}{\partial x_{\alpha_1}}\big) \Big]\\
&+\frac{\partial}{\partial y_i}\Big[ {\sigma_{ij}^{(1)}\frac{\partial }{\partial x_{j}}}\big(N_{\alpha_1}\frac{\partial \phi_{0}}{\partial x_{\alpha_1}}\big) \Big]+\frac{\partial}{\partial x_i}\Big[ {\sigma_{ij}^{(1)}\frac{\partial }{\partial y_{j}}}\big(N_{\alpha_1}\frac{\partial \phi_{0}}{\partial x_{\alpha_1}}\big) \Big].
\end{aligned}
\end{equation}
\begin{equation}
\begin{aligned}
&{F_2}({\bm{x}},{\bm{y}},t) =\frac{\partial}{\partial x_i}\Big( {k_{ij}\frac{\partial u_2}{\partial y_{j}}}\Big)+\frac{\partial}{\partial y_i}\Big( {k_{ij}\frac{\partial u_2}{\partial x_{j}}}\Big)+\varepsilon\frac{\partial}{\partial x_i}\Big(k_{ij}\frac{\partial u_2}{\partial x_{j}}\Big)-\varepsilon\rho c\frac{\partial u_2}{\partial t}.
\end{aligned}
\end{equation}
\begin{equation}
\begin{aligned}
&{E_2}({\bm{x}},{\bm{y}},t) = \frac{\partial}{\partial x_i}\Big( {\sigma_{ij}\frac{\partial \phi_2}{\partial y_{j}}}\Big)+\frac{\partial}{\partial y_i}\Big( {\sigma_{ij}\frac{\partial \phi_2}{\partial x_{j}}}\Big)+\varepsilon\frac{\partial}{\partial x_i}\Big(\sigma_{ij}\frac{\partial \phi_2}{\partial x_{j}}\Big).
\end{aligned}
\end{equation}
\setcounter{equation}{0}
\renewcommand{\theequation}{C.\arabic{equation}}\\
\section*{Appendix C. Mathematical proof of continuous property of microscopic cell functions}\\
Utilizing the similar technique provided in \cite{R37}, we can demonstrate that all auxiliary cell functions are continuous with respect to macroscopic temperature $u_0$. Taking the proof of microscopic cell function ${M_{\alpha_1}}$ as example, $M_{\alpha_1}(\bm{y},u_{0}^{s_1})$ and $M_{\alpha_1}(\bm{y},u_{0}^{s_2})$ are firstly denoted as the microscopic cell function ${M_{\alpha_1}}(\bm{y},u_{0})$ at temperature points $u_0^{s_1}$ and $u_0^{s_2}$ separately. Whereupon the variational equations for ${M_{\alpha_1}}(\bm{y},u_{0}^{s_1})$ and ${M_{\alpha_1}}(\bm{y},u_{0}^{s_2})$ are established by virtue of first-order cell problem (2.11).
\begin{equation}
\!\!-\!\!\int_{\Theta}{ k_{ij}^{(0)}(\bm{y},u_{0}^{s_1}){\frac{\partial M_{\alpha_1}(\bm{y},u_{0}^{s_1})}{\partial y_j}}}\frac{\partial \upsilon^{h_1}}{\partial y_i}d\Theta\!\!=\!\!\int_{\Theta}k_{i{\alpha_1}}^{(0)}(\bm{y},u_{0}^{s_1})\frac{\partial \upsilon^{h_1}}{\partial y_i}d\Theta,\;\forall\upsilon^{h_1}\!\!\in\!\! S_{h_1}(\Theta).
\end{equation}
\begin{equation}
\!\!-\!\!\int_{\Theta}{ k_{ij}^{(0)}(\bm{y},u_{0}^{s_2}){\frac{\partial M_{\alpha_1}(\bm{y},u_{0}^{s_2})}{\partial y_j}}}\frac{\partial \upsilon^{h_1}}{\partial y_i}d\Theta\!\!=\!\!\int_{\Theta}k_{i{\alpha_1}}^{(0)}(\bm{y},u_{0}^{s_2})\frac{\partial \upsilon^{h_1}}{\partial y_i}d\Theta,\;\forall\upsilon^{h_1}\!\!\in\!\! S_{h_1}(\Theta).
\end{equation}
Afterwards, exploiting the variational equation of ${M_{\alpha_1}}(\bm{y},u_{0}^{s_1})$ to subtract the variational equation of ${M_{\alpha_1}}(\bm{y},u_{0}^{s_2})$, this leads to a new variational equation as below.
\begin{equation}
\begin{aligned}
&\int_{\Theta} k_{ij}^{(0)}(\bm{y},u_{0}^{s_1}){{\frac{\partial}{\partial y_j}}}\big[M_{\alpha_1}(\bm{y},u_{0}^{s_1})-M_{\alpha_1}(\bm{y},u_{0}^{s_2})\big]\frac{\partial \upsilon^{h_1}}{\partial y_i}d\Theta\\
&=-\int_{\Theta}\big[k_{ij}^{(0)}(\bm{y},u_{0}^{s_1})-k_{ij}^{(0)}(\bm{y},u_{0}^{s_2})\big]{{\frac{\partial M_{\alpha_1}(\bm{y},u_{0}^{s_2})}{\partial y_j}}}\frac{\partial \upsilon^{h_1}}{\partial y_i}d\Theta\\
&-\int_{\Theta}\big[k_{ij}^{(0)}(\bm{y},u_{0}^{s_1})-k_{ij}^{(0)}(\bm{y},u_{0}^{s_2})\big]\frac{\partial \upsilon^{h_1}}{\partial y_i}d\Theta,\;\forall\upsilon^{h_1}\in S_{h_1}(\Theta).
\end{aligned}
\end{equation}
Then, assuming $|k_{ij}^{(0)}(\bm{y},u_{0}^{s_1})-k_{ij}^{(0)}(\bm{y},u_{0}^{s_2})|\!\leq\! C|u_{0}^{s_1}-u_{0}^{s_2}|$ and replacing $\upsilon^{h_1}$ in (C.3) with ${M_{\alpha_1}}(\bm{y},u_{0}^{s_1})-{M_{\alpha_1}}(\bm{y},u_{0}^{s_2})$, we have
\begin{equation}
\begin{aligned}
&\;\;\;\;\gamma_0\left \|{M_{\alpha_1}}(\bm{y},u_{0}^{s_1})-{M_{\alpha_1}}(\bm{y},u_{0}^{s_2})\right\|_{H^1_0(\Theta)}^2\\
&\!\!\leq\!\!\int_{\Theta} k_{ij}^{(0)}(\bm{y},\!u_{0}^{s_1}){{\frac{\partial\big[M_{\alpha_1}(\bm{y},u_{0}^{s_1})\!\!-\!\!M_{\alpha_1}(\bm{y},u_{0}^{s_2})\big]}{\partial y_j}}}\frac{\partial\big[M_{\alpha_1}(\bm{y},u_{0}^{s_1})\!\!-\!\!M_{\alpha_1}(\bm{y},u_{0}^{s_2})\big]}{\partial y_i}d\Theta\\
&\!\!=\!\!-\!\!\int_{\Theta}\big[k_{ij}^{(0)}(\bm{y},\!u_{0}^{s_1})\!\!-\!\!k_{ij}^{(0)}(\bm{y},u_{0}^{s_2})\big]{{\frac{\partial M_{\alpha_1}(\bm{y},u_{0}^{s_2})}{\partial y_j}}}\frac{\partial \big[M_{\alpha_1}(\bm{y},u_{0}^{s_1})\!\!-\!\!M_{\alpha_1}(\bm{y},u_{0}^{s_2})\big]}{\partial y_i}d\Theta\\
&\!\!-\!\!\int_{\Theta}\big[k_{ij}^{(0)}(\bm{y},u_{0}^{s_1})-k_{ij}^{(0)}(\bm{y},u_{0}^{s_2})\big]\frac{\partial }{\partial y_i}\big[M_{\alpha_1}(\bm{y},u_{0}^{s_1})-M_{\alpha_1}(\bm{y},u_{0}^{s_2})\big]d\Theta\\
&\!\!\leq\!\! C|u_{0}^{s_1}-u_{0}^{s_2}|\left \|M_{\alpha_1}(\bm{y},u_{0}^{s_1})-M_{\alpha_1}(\bm{y},u_{0}^{s_2})\right\|_{H^1_0({\Theta})}
\end{aligned}
\end{equation}
Consequently, in case $u_{0}^{s_1}\rightarrow u_{0}^{s_2}$, it is evident that ${M_{\alpha_1}}(\bm{y},u_{0}^{s_1})\rightarrow{M_{\alpha_1}}(\bm{y},u_{0}^{s_2})$ from (C.4), which means the continuous property of microscopic cell functions holds. In the future, we shall further excavate how to distribute representative macroscopic temperature points for acquiring optimally computing accuracy.

\bibliographystyle{siamplain}
\bibliography{references1}
\end{document}